\pdfoutput=1
  
\documentclass[10pt]{article}

\usepackage{amssymb,amsmath,psfrag,graphicx}
\usepackage{mathtools} 
\usepackage{verbatim}
\usepackage{color}
\usepackage{morefloats}
\usepackage{colortbl}
\usepackage{fullpage}
\usepackage{enumitem}
\usepackage{tabularx}
\usepackage{multicol}
\usepackage{algorithm}
\usepackage{algpseudocode}
\usepackage{trimclip} 
\usepackage{caption}  
\usepackage{subcaption} 

\usepackage{pst-grad}
\usepackage{pgfplots}
\usepackage[all]{xy}
\usetikzlibrary{calc}

\usepackage{soul}

\usepackage{marginnote}

\newtheorem{remark}{Remark}[section]

\definecolor{violachiaro}{rgb}{.8,0.69,1}

\definecolor{yel}{rgb}{1,1,0.554}
\newcommand{\nota}[1]{\ifodd\thepage{\marginnote{#1 }}\else{\reversemarginpar \marginnote{#1 }\reversemarginpar}\fi}

\definecolor{m}{rgb}{1,0,0}
\definecolor{s}{rgb}{0,0.666,0}
\definecolor{l}{rgb}{0,0,1}

\graphicspath{{}}

\begin{document}

\title{The INTERNODES method for the treatment of non-conforming multipatch
geometries in Isogeometric Analysis}

\author{Paola Gervasio \footnotemark[3] \and %
                   Federico Marini \footnotemark[4] }
\maketitle
\renewcommand{\thefootnote}{\fnsymbol{footnote}}
\footnotetext[3]{DICATAM, Universit\`a degli Studi di Brescia, via Branze 38, I-25123 Brescia, Italy ({\tt paola.gervasio@unibs.it})}

\footnotetext[4]{IMATI - CNR, Pavia, Italy ({\tt federico@imati.cnr.it})}

\begin{abstract}
In this paper we apply the INTERNODES method to solve second order elliptic
problems discretized by Isogeometric Analysis methods 
on non-conforming multiple patches in 2D and 3D geometries.
INTERNODES is an interpolation-based method 
that, on each interface of the
configuration, exploits two independent interpolation operators
to enforce the continuity of the traces and of the normal derivatives. 
INTERNODES easily handles both parametric and geometric NURBS
non-conformity.
We specify how to set up the interpolation matrices on
non-conforming interfaces, how to enforce
the continuity of the normal derivatives and we give special attention to
implementation aspects.
The numerical results show that INTERNODES exhibits optimal convergence rate
with respect to the mesh size of the NURBS spaces an that it is robust with
respect to jumping coefficients.
\end{abstract}


\section{Introduction}

Nowadays Isogeometric Analysis (IgA) \cite{chb_iga_book}
represents one of the most popular methods
for numerical simulations.
Its paradigm consists in expanding the 
Partial Differential Equations (PDE)
solution with respect to the basis functions of the same type of the ones 
(either B-splines or NURBS) used to
describe the geometry of the computational domain generated by CAD software.

Often, real-life problems are defined on complex geometries 
that usually consist of several patches.  Moreover these patches can feature
non-conformity. 

\begin{figure}
\begin{center}
\begin{subfigure}[t]{2in}
\includegraphics[width=\textwidth,height=3.cm]{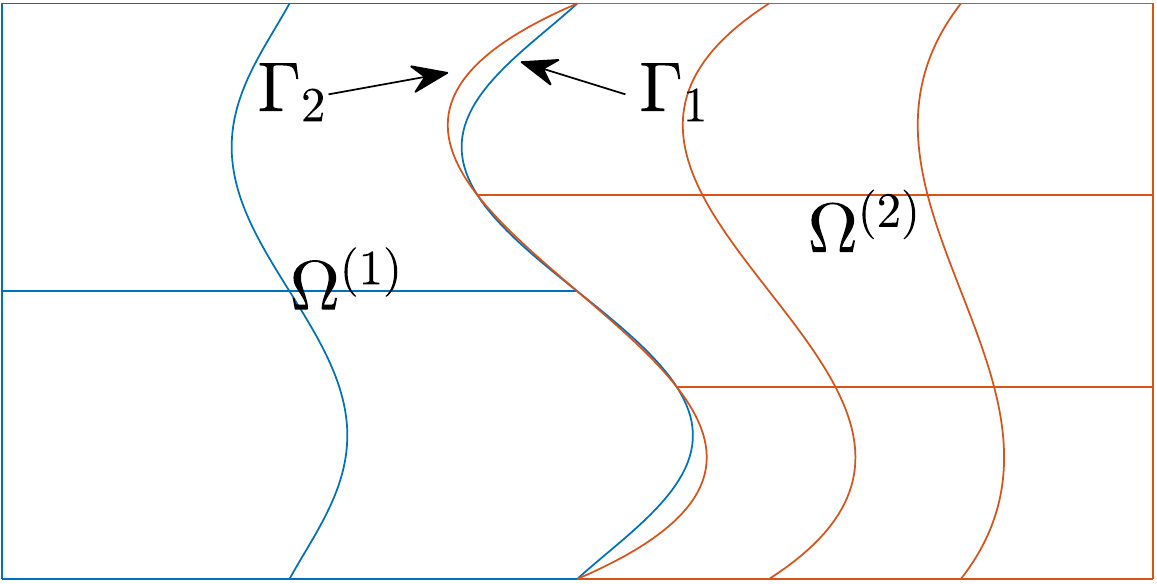}
\subcaption{}
\end{subfigure}
\begin{subfigure}[t]{2in}
\includegraphics[width=\textwidth,height=3.cm]{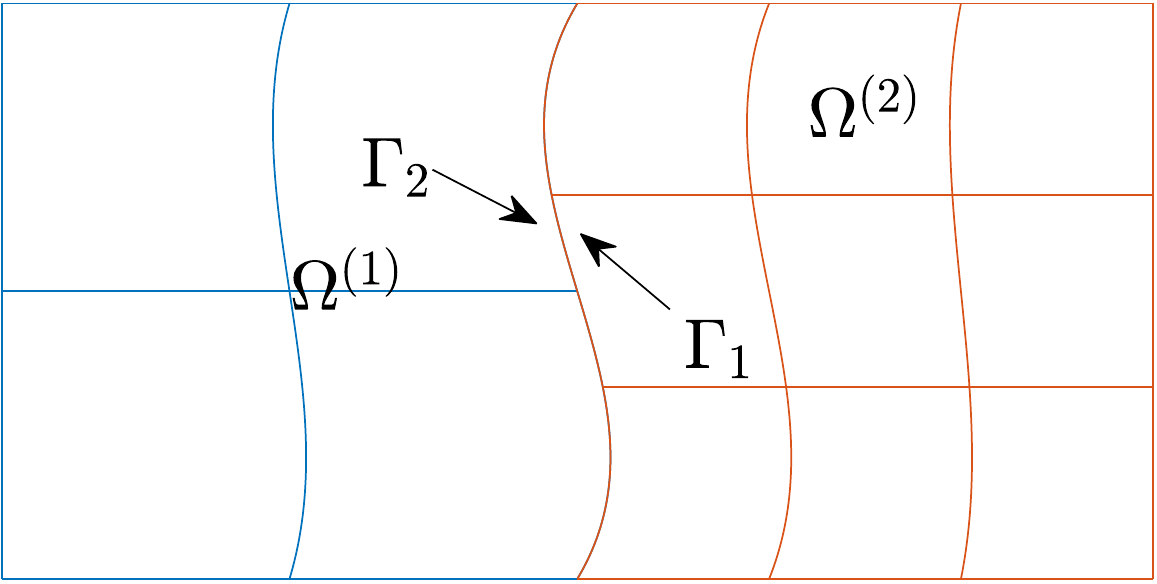}
\subcaption{}
\end{subfigure}
\begin{subfigure}[t]{2in}
\includegraphics[width=\textwidth,height=3.cm]{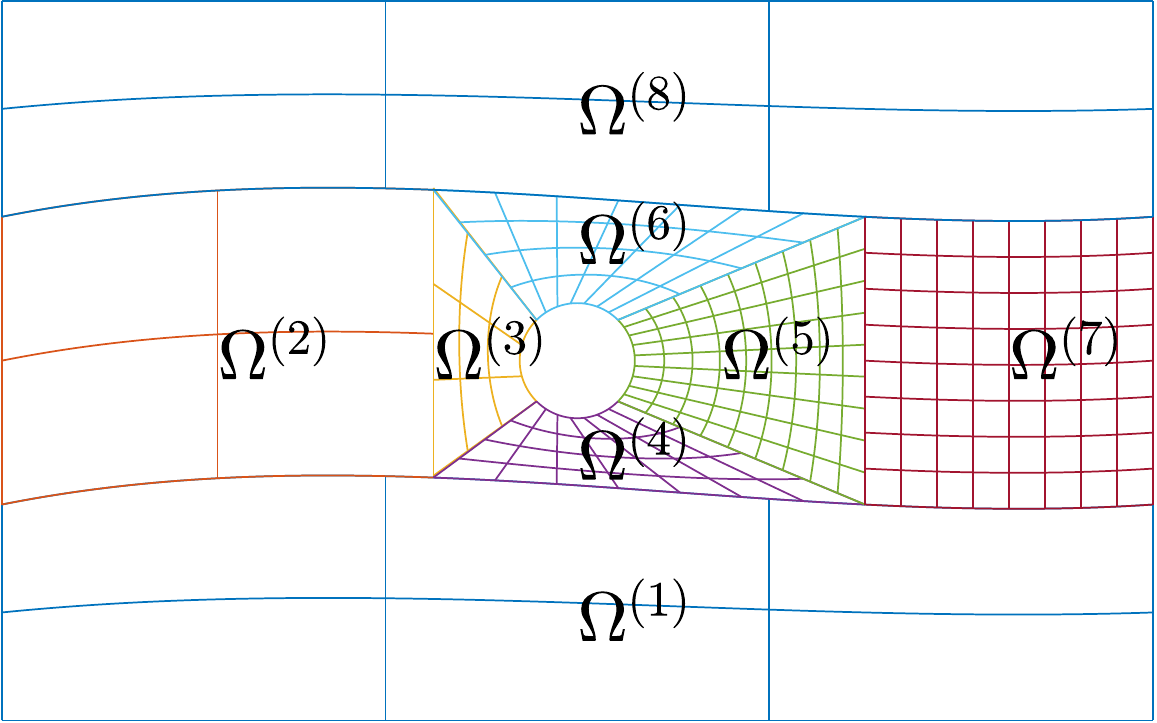}
\subcaption{}
\end{subfigure}
\end{center}
\caption{Non-conforming patches: (a) non-watertight patches, geometric
non-conformity; 
(b) and (c)  watertight patches but non-matching
parametrizations at the interface
}
\label{fig:geometrie}
\end{figure}

By non-conformity we mean either  \emph{geometrical
non-conformity} or \emph{parametric non-conformity}.
 The former one occurs when two adjacent patches share a common
boundary in the physical space only approximately (e.g., as 
result of CAD modeling operation) and we say 
that the interfaces are \emph{non-watertight}; an example is shown
in  Fig. \ref{fig:geometrie} (a).

In the case that two interfaces are watertight, however we can face 
parametric non-conformity that means that
different and totally unrelated discretizations 
(not necessarily the refinement one of the other) are considered inside 
the patches sharing the same
interface (or only a part of it); two examples are given in Fig. \ref{fig:geometrie} (b) and (c).

In the last years, the treatment of multipatch geometries has been investigated
in several papers, far from be exhaustive we mention
\cite{chr07,bbcs_tsplines,hesch_betsch_11,nguyen_etal_2014,ruess_etal_2014,zhu_2017,coox_etal_17,chan_etal_18,mi_zheng_18,
kleiss_ieti_12,hofer_langer_17,stavroulakis_ddiga_18,
brivadis_mortar_15,delorenzis_2011,de2012mortar,kargaran_2019}.
In this paper we propose to 
apply the INTERNODES method to solve elliptic problems within the IgA framework
on non-overlapping and non-conforming multipatch configurations.

Overlapping multi-patch domains without trimming have been faced 
in the Isogeometric Analysis context
in \cite{kargaran_2019}. In this very recent paper the authors 
propose the Overlapping Multi-Patch method, which is a non-iterative
reformulation of the Schwarz method (see, e.g. \cite{tw,qv_ddm}).

INTERNODES (INTERpolation for NOnconforming DEcompositionS) \cite{dfgq,
gq_internodes} is a general
purpose method to deal with non-conforming discretizations of PDEs in 2D and 3D
geometries split into non-overlapping subdomains (in fact 
IgA patches play the role of non-overlapping subdomains
in the domain decomposition  context \cite{tw,qv_ddm}).
The method was proposed in
\cite{dfgq} to solve elliptic PDEs by Finite Element Methods (FEM) and 
Spectral Element Methods (SEM) on two non-conforming subdomains, then 
its theoretical analysis, as well as its extension to decompositions with more
than two subdomains, has been carried out in
\cite{gq_internodes,gq_dd24,gq_camc}. 

INTERNODES has been successfully applied to solve Navier-Stokes equations
(\cite{dfq_fsi,forti_phd}) and multi-physics problems like the
Stokes-Darcy coupling  to simulate the filtration of fluid in porous domains
\cite{qg_dd24,gq_camc} and the Fluid-Structure Interaction problem
\cite{dfq_fsi,forti_phd}. 

It has been proved that, inside the FEM-SEM
context, INTERNODES exhibits optimal accuracy with respect
to the $H^1$-broken norm, i.e., the error between the global INTERNODES
solution and the exact one is proportional to
the best approximation error inside the subdomains.

Inside the subregions (or patches) of the decomposition we discretize the PDE
(here we have implemented the Galerkin formulation of IgA) by using 
NURBS spaces that are totally unrelated one each other.

To enforce the continuity of the traces and
the equilibration of the fluxes across the interfaces
between two adjacent patches, INTERNODES exploits two
independent interpolation operators: 
one for transferring the Dirichlet trace of the solution,
the other for the normal derivatives.
Like mortar methods, INTERNODES tags the opposite sides of an interface 
either master or slave: the
continuity of the traces is enforced on the slave side of the
interface (more precisely the Dirichlet trace is interpolated from the master
side to the slave one), while the equilibration of the fluxes is enforced 
on the master side of the same interface (the normal derivative is 
interpolated in a suitable way from the slave side to the master one).

In this paper we apply the interpolation at the 
Greville abscissae of the knot vectors \cite{piegl_nurbs_book,demko,abhrs_2010},
nevertheless other choices are possible (see, e.g.
\cite{abhrs_2010,montardini_st_cmame17}).

 One of the strength of
INTERNODES just consists  in working with totally unrelated
discretization spaces (with different sets of control points and
weights and different basis functions) without the need to enrich the
NURBS spaces or to insert necessary knots (as, e.g. T-splines do)
to ensure full compatibility at the interfaces.

When two interfaces are watertight (as in Fig. \ref{fig:geometrie} (b) and
(c))
 the point inversion from one interface to the adjacent one is
well defined also for non-matching parametrizations and standard 
NURBS interpolation is exploited to implement
INTERNODES.
When instead the interfaces are non-watertight (as in Fig. \ref{fig:geometrie}
(a)) we overcome the difficulty to
 project the Greville
abscissae from one face to the non-watertight corresponding one, by
exploiting the Radial Basis Functions interpolation
(see, e.g., \cite{dfq_rbf,dfq_fsi} in the FEM  context and
\cite{mi_zheng_18} in the IgA context).

To interpolate correctly the normal derivatives we need to assemble local 
interface mass matrices, but differently than in mortar methods,
no cross-mass matrix involving basis
functions living on the two opposite sides of the interface  and no
ad hoc  numerical quadrature (\cite{brivadis_mortar_quad})
are required by INTERNODES to build the inter-grid operators. 

To solve the  multipatch problem at the algebraic level,
the degrees of freedom internal to
the patches can be eliminated and the Schur complement system associated with
the degrees of freedom on the master skeleton can be solved by Krylov methods
(e.g., Bi-CGStab \cite{vander} or GMRES \cite{saad_schultz}),
as typical in domain decomposition methods of sub-structuring type.
In the case of only two patches, or when the decomposition is chessboard
like so that we can tag all the interfaces of a single patch as either master
or slave, we assemble the preconditioner  for the global 
Schur complement system starting from the 
 local Schur complement matrices associated with the master patches.

The numerical results of
 Sections \ref{sec:numres_2dom} and \ref{sec:numres_mdom} 
show that INTERNODES applied to IgA discretizations
exhibits optimal accuracy versus the mesh size $h$
for both 2D and 3D geometries and it is robust with respect to both
jumping coefficients and non-watertight interfaces.

This is the first paper that joins INTERNODES and IgA and a lot of questions
remain open:
 the analysis of the convergence rate in the IgA framework, 
the efficient 
solution of the Schur complement system by designing suitable
preconditioners in the case that a patch features both master and slave
edges, the formulation of the method on surfaces in 3D, its
application to contact mechanics problems and, last but not
least, the extension of the method to deal with multi-physics problems.
Even though these are indeed challenging tasks, the authors of this paper have
no reason to think that they
cannot be accomplished within INTERNODES, being the theoretical setting
presented herein clearly and the results promising. 
Compared to the mortar method,
the removal of the necessity of inter-grid quadrature is one of the most
attractive feature of INTERNODES. The authors
believe that it alone can be a sufficient reason to further develop INTERNODES
in the IgA framework.

The paper is organized as follows.
In Sect. \ref{sec:setting} we formulate the transmission problem; in Sect.
\ref{sec:internodes} we present INTERNODES for two patches; in Sect.
\ref{sec:discretization} we recall the definition of the NURBS basis functions,
we define the interpolation operators at the Greville nodes 
for both watertight and non-watertight configurations,
and we specify how to interpolate the
normal derivatives at the interface. In Sect. 
\ref{sec:alg} we give the algebraic formulation of the method on two patches,
while in Sect. \ref{sec:internodesM} we present INTERNODES on more general
configurations with $M>2$ patches. Finally, in Sect. \ref{sec:schurM} 
we provide the algorithms to implement INTERNODES and solve the Schur
complement system with respect to degrees of freedom on the master skeleton.
The numerical results are shown in Sect. \ref{sec:numres_2dom} and
\ref{sec:numres_mdom}.


\section{Problem setting}\label{sec:setting}

Let $\Omega\subset  {\mathbb R}^d$, with $d=2,3$,  be an open domain with
Lipschitz boundary $\partial\Omega$
$f\in L^2(\Omega)$, $\alpha\in L^\infty(\Omega)$ 
and $g\in H^{1/2}(\partial \Omega)$ be given functions.
We look for the solution $u$ of the self-adjoint second order elliptic problem 
\begin{eqnarray}\label{problem}
\left\{\begin{array}{ll}
-\Delta u+ \alpha u=f & \mbox{ in }\Omega\\
u=g & \mbox{ on }\partial\Omega.
\end{array}\right.
\end{eqnarray}

\noindent
For sake of simplicity, in the first part of the paper we deal with this simple
problem, then starting from Sect. \ref{sec:general_L} on wards we extend the method to
more general elliptic operators.

We denote by $\widetilde g$ a lifting of the Dirichlet datum $g$, i.e.
any function $\widetilde g\in H^1(\Omega)$ such that
$\widetilde g|_{\partial\Omega}=g$.
%

The weak form of problem (\ref{problem}) reads: find $u
\in H^1(\Omega)$ with $(u-\widetilde g)\in H^1_0(\Omega)$ such that
\begin{equation}\label{weak-problem}
a(u,v)=(f,v)_{L^2(\Omega)}\hskip 1.cm 
\forall v\in H^1_0(\Omega),
\end{equation}

\noindent
where 
$\displaystyle a(u,v)=\int_\Omega \nabla u \cdot \nabla v+ \alpha u v\ d\Omega$.

Under the assumption that
$\alpha\geq 0$ a.e. in $\Omega$, 
problem (\ref{weak-problem}) admits a unique solution (see, e.g., \cite{qvb})
that is stable w.r.t. the data $f$ and $g$.

In the next Section we introduce the INTERNODES method on 2-patches
decompositions. The more general case will be faced in Sect.
\ref{sec:internodesM}.

\section{The transmission problem for two subdomains}

We define a non-overlapping decomposition of  $\Omega$ into two 
subdomains
 $\Omega_1$ and $\Omega_2$ with Lipschitz boundary, such that
$$\overline\Omega=\overline{\Omega^{(1)}}\cup\overline{\Omega^{(2)}}, \qquad 
\Omega^{(1)}\cap\Omega^{(2)}=\emptyset,$$ 
while
$\Gamma_{12}=\overline{\Omega^{(1)}}\cap\overline{\Omega^{(2)}}$ is the common
interface  that we 
assume be of class $C^{1,1}$  
(see \cite[Def. 1.2.1.2]{grisvard85}) to allow the normal derivative of $u$ on
it to be well defined.

Then, for $k=1,2$ we define:
$\partial\Omega^{(k)}_D=\partial{\Omega}^{(k)}\cap\partial\Omega_D$. Let
$u^{(k)}$ be the restriction of $u$ to $\Omega^{(k)}$, then $u^{(1)}$ and 
$u^{(2)}$ are the solutions of
the transmission problem 
(see \cite{qv_ddm})
\begin{eqnarray}\label{transmission_strong}
\left\{ \begin{array}{ll}
-\Delta u^{(k)}+\alpha u^{(k)} = f  
&\textrm{ in }\Omega^{(k)}, \quad k=1,2 \\
u^{(k)} = g  &\textrm{ on }\partial\Omega^{(k)}_D, \quad k=1,2 \\
u^{(1)} = u^{(2)}  &\textrm{ on }\Gamma_{12} \\
\displaystyle \frac{\partial u^{(1)}}{\partial {\bf n}_1}
+\frac{\partial u^{(2)}}{\partial {\bf n}_2}= 0 &\textrm{ on }\Gamma_{12},
      \end{array}
    \right.
\end{eqnarray}
where ${\bf n}_k$ is the outward unit normal vector to $\partial\Omega^{(k)}$
(on $\Gamma_{12}$ it holds ${\bf n}_1=-{\bf n}_2$).

For $k=1,2$, 
we define the functional spaces
\begin{eqnarray*}
V^{(k)}=\left\{ v \in H^1(\Omega) |\  v=0 \mbox{ on } 
\partial\Omega^{(k)}_D\right\}, \qquad
V^{(k)}_0=H^1_0(\Omega^{(k)}),\\
\Lambda = \left\{\lambda \in H^{1/2}(\Gamma_{12})\,\,|\,\, \exists v\in 
H^1(\Omega) \mbox{ such that } v|_{\Gamma_{12}}=\lambda \right\},
\end{eqnarray*}
noticing that $\Lambda = H^{1/2}_{00}(\Gamma_{12})$ if
$\Gamma_{12}\cap\partial\Omega\neq\emptyset$.

We denote by 
$\displaystyle a^{(k)}(u,v)=\int_{\Omega^{(k)}}\nabla u\nabla v+\alpha uv\,d\Omega$
 the restriction of the bilinear form $a(\cdot,\cdot)$ to $\Omega^{(k)}$ and
we set
$\widetilde g^{(k)}=\widetilde g_{|\Omega^{(k)}}$.

The weak form of the transmission problem (\ref{transmission_strong})
reads  (see \cite{qv_ddm}): for $k=1,2$ look for
$u^{(k)}\in H^1(\Omega^{(k)})$ with $(u^{(k)}-\widetilde g^{(k)})\in V^{(k)}$
such that
\begin{eqnarray}\label{transmission_weak}
\left\{
\begin{array}{ll}
a^{(k)}(u^{(k)},v^{(k)}) = {\cal F}^{(k)}(v^{(k)})
&\forall v^{(k)}\in V^{(k)}_0 \\[2mm]
u^{(1)}= u^{(2)} &\textrm{on }\Gamma_{12} \\[2mm]
\displaystyle
\sum_{k=1,2}a^{(k)}(u^{(k)},{\cal L}^{(k)}\eta)=
\sum_{k=1,2}{\cal F}^{(k)}({\cal L}^{(k)}\eta)
&\forall \eta\in\Lambda
 \end{array}
\right.
\end{eqnarray}
where 
\begin{equation}\label{functionalF}
{\cal F}^{(k)}(v^{(k)})=(f,v^{(k)})_{L^2(\Omega^{(k)})}
\qquad \forall v^{(k)}\in
V^{(k)},
\end{equation}
while
\begin{equation}\label{R}
{\cal L}^{(k)}:\Lambda\to V^{(k)},\quad s.t. \quad 
({\cal L}^{(k)}\eta)|_{\Gamma}=\eta \qquad \forall \eta \in \Lambda
\end{equation}
denotes any possible linear and continuous \emph{lifting operator}
from $\Gamma_{12}$ to $\Omega^{(k)}$.

\begin{remark}\label{rem:flussi}
{\rm
Denoting by $\langle\cdot,\cdot\rangle$ the duality pairing between $\Lambda$
and its dual space $\Lambda'$, the
distributional form of the interface condition 
(\ref{transmission_strong})$_4$
reads
\begin{equation}\label{flussi}
\langle \frac{\partial u^{(1)}}{\partial {\bf n}_1}
+\frac{\partial u^{(2)}}{\partial
{\bf n}_2},\eta\rangle = 0 \qquad \forall \eta\in\Lambda
\end{equation}
and it is equivalent to (\ref{transmission_weak})$_3$. 
(In the case that $f\in L^2(\Gamma)$, we have  $\langle f,
v\rangle=\int_\Gamma f v\, d\Gamma$.)
}

\end{remark}
\subsection{Formulation   of INTERNODES}\label{sec:internodes}

Bearing in mind the Isogeometric Analysis framework, 
the two subdomains $\Omega^{(1)}$ and $\Omega^{(2)}$ 
introduced in the previous section play the role of two
disjoint patches of a suitable multipatch decomposition of $\Omega$.

For $k=1,2$, let ${\cal N}^{(k)}_{h_k}$ be two finite dimensional spaces
arising from Isogeometric Analysis discretization
that can be totally unrelated to each other and
set 
\begin{equation}\label{spaces1}
V^{(k)}_{h_k}={\cal N}^{(k)}_{h_k}\cap V^{(k)},\qquad
V^{(k)}_{0,h_k}={\cal N}^{(k)}_{h_k}\cap V^{(k)}_0.
\end{equation}
Then we denote by 
$u_{h_k}^{(k)}\in {\cal N}^{(k)}_{h_k}$ the approximation of $u^{(k)}$ we are
looking for.

\begin{figure}[b!]
\begin{center}
\scalebox{0.4}{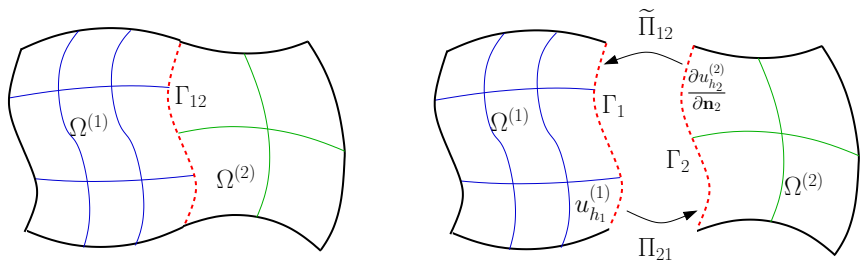}
\end{center}
\caption{The interface
$\Gamma_{12}=\partial\Omega^{(1)}\cap \partial\Omega^{(2)}$ and the two sides
$\Gamma_1$ and $\Gamma_2$ for a multipatch geometry when $d=2$.
The core idea of INTERNODES: $\Pi_{21}$ interpolates the trace from 
$\Gamma_1\subset\partial\Omega^{(1)}$ to $\Gamma_2\subset\partial\Omega^{(2)}$, 
$\widetilde\Pi_{12}$ interpolates the normal derivative
 from $\Gamma_2$ to $\Gamma_1$}
\label{fig:Pi}
\end{figure}

Let us denote by $\Gamma_1$ and $\Gamma_2$ 
the two sides of $\Gamma_{12}$ as part of the boundary of either $\Omega^{(1)}$
or  $\Omega^{(2)}$ (see Fig. \ref{fig:Pi}), and 
by $Y^{(k)}_{h_k}$ the space of the trace on $\Gamma_k$ of the functions of
${\cal N}^{(k)}_{h_k}$, for $k=1,2$ (see Fig. \ref{fig:Pi}).

Even if $\Gamma_1$ and $\Gamma_2$ may represent the same geometric
curve (when $d=2$) or surface (when $d=3$), 
we distinguish them to underline  on which side of the interface we are
working.

To use non-conforming discretizations in $\Omega^{(1)}$ and
$\Omega^{(2)}$  implies that the trace spaces
$Y^{(1)}_{h_1}$ and $Y^{(2)}_{h_2}$ may not match.
In such a case,
to enforce the continuity of the trace (i.e., the interface condition 
(\ref{transmission_weak})$_2$) and the equilibration of normal derivatives
 (i.e., the interface condition (\ref{flussi}) or, equivalently,
(\ref{transmission_weak})$_3$), we introduce two independent interpolation
 operators:
the first one, named $\Pi_{21}$, is designed
to interpolate the trace of $u_{h_1}^{(1)}$ from $\Gamma_1$ to
$\Gamma_2$, while
the second one, named $\widetilde\Pi_{12}$, is used to interpolate in a
suitable way the normal derivative $\frac{\partial u^{(2)}_{h_2}}{\partial
{\bf n}_2}$ from $\Gamma_2$ to $\Gamma_1$ (see Fig. \ref{fig:Pi}).

We give here
the basic idea of the INTERNODES method when it is applied  to
the weak transmission problem
(\ref{transmission_weak}), and we postpone the rigorous description of 
the method to the next sections, after defining
the interpolation operators and after explaining how to transfer
the normal derivative across the interface.

The INTERNODES method applied to (\ref{transmission_weak}) reads as follows.
For $k=1,2$, let $\widetilde g_{h_k}^{(k)}\in{\cal N}^{(k)}_{h_k}$ be 
a suitable approximation of $\widetilde g^{(k)}$.
 Then, for $k=1,2$ we look for
$u^{(k)}_{h_k}\in{\cal N}^{(k)}_{h_k}$ such that $(u_{h_k}^{(k)}-\widetilde 
g_{h_k}^{(k)})\in V^{(k)}_{h_k}$ and
\begin{eqnarray}\label{internodes_weak}
\left\{
\begin{array}{ll}
a^{(k)}(u_{h_k}^{(k)},v_{h_k}^{(k)}) = {\cal F}^{(k)}(v_{h_k}^{(k)})
&\ \forall v_{h_k}^{(k)}\in V^{(k)}_{0,h_k},\quad k=1,2 \\[2mm]
u^{(2)}_{h_2}=\Pi_{21}u^{(1)}_{h_1} & \mbox{ on }\Gamma_2\\
\displaystyle\langle \frac{\partial u^{(1)}_{h_1}}{\partial {\bf n}_1}
+\widetilde\Pi_{12}\frac{\partial u^{(2)}_{h_2}}{\partial
{\bf n}_2},\eta^{(1)}_{h_1}\rangle = 0 &\  \forall \eta^{(1)}_{h_1}\in
Y^{(1)}_{h_1}.
\end{array}\right.
\end{eqnarray}

The interface condition
(\ref{internodes_weak})$_2$ characterizes the role of the 
interfaces $\Gamma_1$ and $\Gamma_2$. Since 
 the trace on $\Gamma_2$ depends on the trace on $\Gamma_1$,
following the terminology typical of mortar methods, 
the interface  $\Gamma_1$ is named
\emph{master}, while $\Gamma_2$ is named \emph{slave}.

\begin{remark}[Analysis of INTERNODES]
{\rm 
The INTERNODES method has been analyzed in \cite{gq_internodes} in the Finite
Element framework.  More precisely, if quasi-uniform and affine triangulations
are considered inside each subdomain and Lagrange interpolation is applied 
to enforce the interface conditions, 
it has been proved (\cite{gq_internodes}) that INTERNODES
yields a solution that is unique, stable, and
convergent with an \emph{optimal}
rate of convergence (i.e., that of the best approximation error
in every subdomain).

Two interpolation operators are needed to guarantee the optimal convergence
rate of the method with respect to the discretization parameters.
As a matter of fact, it is well known that using a single
interpolation operator (jointly with its transpose) instead of two different
operators is not optimal.
The approach using a single interpolation operator is also known as 
\emph{point-wise approach}, see
\cite{bmp94,bbdmkmp} and \cite[Sect. 6]{dfgq}.

The same arguments used in \cite{gq_internodes} can be used in
the Isogeometric Analysis framework too, to prove the
existence, the uniqueness, and  the stability of the solution of problem
(\ref{internodes_weak}). 

A convergence theorem, establishing the error bound for
the INTERNODES method with respect to the mesh size $h=\max_k h_k$
in the framework of Isogeometric Analysis
is an open problem at the moment of writing the present paper.

Nevertheless, the numerical results provided in
 the next Sections show that INTERNODES
exhibits optimal accuracy versus the mesh size $h$
for both 2D and 3D geometries.
}
\end{remark}

\begin{remark}\label{rem:strong_weak}
{\rm INTERNODES could be applied to the strong form (\ref{transmission_strong}) of the
transmission problem. In this case it is sufficient to replace
(\ref{internodes_weak})$_1$ with the discrete counterpart of
(\ref{transmission_strong})$_1$.
}
\end{remark}


\section{Discretization by Isogeometric Analysis}
\label{sec:discretization}

Let $Z=\{0=\zeta_0,\zeta_1,\ldots,\zeta_{n-1}, \zeta_{n_{el}}=1\}$ 
be  the set of $(n_{el}+1)$ distinct
knot values in the one-dimensional
patch $[0,1]$. Given the integer $p\ge 1$, let
\begin{equation}
\Xi=\{\xi_1,\xi_2,\ldots,\xi_q\}=\{\underbrace{\zeta_0,\ldots, \zeta_0}_{p+1}, 
\underbrace{\zeta_1,\ldots,\zeta_1}_{m_1}, \ldots,
\underbrace{  \zeta_{n_{el}-1},  \ldots,  \zeta_{n_{el}-1}}_{m_{n_{el}-1}},
\underbrace{\zeta_{n_{el}},\ldots, \zeta_{n_{el}}}_{p+1} \}
\end{equation}

\noindent be a $p$-open knot vector, whose internal knots are repeated
at most $p$ times, their
multiplicity being denoted $m_i$.
If $q$ is the cardinality of $\Xi$, we consider
the number $n=n(\Xi) = q - p - 1$.

Starting from the knot vector $\Xi$
we define the $n(\Xi)$ uni-variate B-spline functions of degree $p$ and of global regularity
$C^{p-\textrm{max}_i\{m_i\}}$ in the patch $[0,1]$ by means of the Cox-de Boor recursion
formula as follows
(\cite{chb_iga_book}). 
For $i=1,\ldots,n(\Xi)$,  set
\begin{eqnarray}\label{N0}
\widehat B_{i,0}(\hat x)=\left\{\begin{array}{ll}
1 & \mbox{if }\xi_i\leq \hat x< \xi_{i+1}\\
0 &\mbox{otherwise,}
\end{array}\right.
\end{eqnarray}
 and, for $\ell=1,\ldots, p$ and $i=1,\ldots,q-\ell-1$, set
\begin{eqnarray}\label{bsplines}
\widehat B_{i,\ell}(\hat x)=\frac{\hat x-\xi_i}{\xi_{i+\ell}-\xi_i}\widehat
B_{i,\ell-1}(\hat x)
+ \frac{\xi_{i+\ell+1}-\hat x}{\xi_{i+\ell+1}-\xi_{i+1}}\widehat
B_{i+1,\ell-1}(\hat x).
\end{eqnarray}

The $d$-times tensor product of the set $Z$ induces a Cartesian
grid in the parameter domain $\widehat \Omega=[0,1]^d$. 
Then we exploit the tensor product rule for the
construction of multivariate B-splines functions:
\begin{equation}
  \widehat B_{i_1,\ldots,i_d,\,p_1,\ldots,p_d}(\hat x_1,\ldots,\hat x_d)=
  \widehat B_{i_1,\,p_1}(\hat x_1)\cdots\widehat B_{i_d,\,p_d}(\hat x_d).
\end{equation}

We assume for sake of simplicity that the  knots $\zeta_i$
are equally spaced (i.e. the resulting knot vector is \emph{uniform})
 along all the parameter directions, and we define 
the mesh size $h=1/n_{el}$. 
We also assume that the multiplicities of the internal knots are all equal to
1, thus the resulting B-spline functions belong to $C^{p-1}$.
Finally we assume that the knots vectors $\Xi_1,\ldots,\Xi_d$
and
the polynomial degrees $p_1,\ldots,p_d$ are the same along any 
direction of the parameter domain, bearing in
mind that what we are going to formulate applies as well
to more general situations for which either
different knot vectors (uniform or non-uniform) 
or different polynomial degrees  or different global 
regularities are considered along the
directions of the parameter domain.
With these assumptions, the number of uni-variate basis functions along each
direction is equal to $n = n_{el} + p$ and, by
tensor product means, the number of multivariate basis functions is $N=n^d$.

B-splines are the building blocks for the parametrization of geometries
$\Omega\subset\mathbb R^d$ of interest. Given a
set of $N$ so-called \emph{control points} ${\bf P}_i\in{\mathbb R}^d$, the
geometrical map ${\bf F}:[0,1]^d \longrightarrow \Omega$
defined as\footnote{For the sake of simplicity, we index the B-spline basis
functions with a uni-variate index $i=1,\ldots,N$
instead of a more precise multi-index ${\bf i}=(i_1,\ldots,i_d)$. The same
simplification holds for $p$.}
\begin{equation}
  {\bf F}(\hat x_1,\ldots,\hat x_d) = \sum_{i=1}^N {\bf P}_i \widehat B_{i,p}
(\hat x_1,\ldots,\hat x_d)
\end{equation}
is a parametrization of $\Omega$, its shape being 
governed by the control points.
This can be seen as the starting point of the
techniques typically adopted by the CAD community for the representation of
geometries.

Even though they can be used to parametrize a wide variety of shapes,
B-splines do not allow to exactly represent objects
such as conic sections and many others typical of the engineering
design. To overcome this drawback the CAD community and hence 
Isogeometric Analysis exploits NURBS (Non-Uniform Rational B-Splines). 

Given a set of positive weights
$\{w_1,\ w_2,\ldots, w_{N}\}$ associated with the control points ${\bf P}_i$,
multivariate NURBS basis functions are 
\begin{equation}\label{multivariate_nurbs}
\widehat \varphi_{i,p}({\bf x})=
\frac{\widehat B_{i,p}(\widehat x_1,\ldots,\widehat x_d)\,
w_{i}} 
{\sum_{j=1}^{n^d} \widehat B_{j,p}(\widehat x_1,\ldots,\widehat x_d)
w_{j}}, \qquad
i=1,\ldots,n^d.
\end{equation}

Notice that, 
by the definition of the knot vectors, the basis functions
associated with the corners of $\widehat\Omega$ are interpolatory.

Then we denote by 
\begin{equation}\label{NURBS_space_parameter}
\widehat {\cal N}_h=\widehat {\cal N}_h(\Xi,p)=span\{\widehat \varphi_{i,p}
(\widehat{\bf x}),\
i=1,\ldots,N\}
\end{equation}
the space spanned by the multivariate NURBS 
basis functions (\ref{multivariate_nurbs})
 on the parameter domain
$\widehat\Omega$. The sub-index $h$ is an abridged notation that expresses the
dependence of the space on both
the number of elements $n_{el}$ induced by the knot vector $\Xi$ and the
polynomial degree $p$ of the B-spline.

NURBS are in fact piecewise rational B-splines and inherit the  global
continuity in the patch by the B-spline $\widehat B_{i,p}$.
The index $p$ used in the definition of $\widehat \varphi_{i,p}$ represents the
polynomial degree of the originating B-splines, but it is evident that
$\widehat  \varphi_{i,p}$ in general are not piecewise  polynomials.

For a deeper analysis of NURBS basis functions and their
practical use in CAD frameworks, we refer to \cite{piegl_nurbs_book}.

Even if in \cite{bbrs_acta} it is shown that the isoparametric paradigm
 can be relaxed 
(i.e. by using a NURBS space for the parametrization of $\Omega$ via 
${\bf F}$ and an 
unrelated B-spline space for the discretization of the PDE), 
in this paper we will follow this concept and we will use 
the same NURBS space for both the parametrization of the subdomain and for the
discrete space.

\medskip
Since the discretizations in the two patches $\Omega^{(1)}$ and $\Omega^{(2)}$
are independent of each-other, 
for any $k=1,2$, we consider $p-$open multivariate knot vectors $\Xi^{(k)}$ 
(with  $\Xi^{(1)}$ and  $\Xi^{(2)}$ independent of
each other) and 
polynomial degrees $p^{(k)}$ (again with $p^{(1)}$ and $p^{(2)}$ 
independent of each other). 
Again, the cardinality of the associated NURBS spaces are $n^{(k)}=n_{el}^{(k)}
+ p^{(k)}$ along each direction so that
the global cardinality of the multivariate NURBS  space is
$N^{(k)}=(n^{(k)})^d$.

Given a set of real positive weights in each patch,
 we define the (parameter-)multivariate NURBS spaces
\begin{equation}\label{NURBS_space_parameter_k}
\widehat {\cal N}_{h_k}^{(k)}=\widehat {\cal N}_{h_k}^{(k)}(\Xi^{(k)},p^{(k)})
=span\{\widehat \varphi^{(k)}_{i,p^{(k)}},
\ i=1,\ldots,N^{(k)}\},\qquad k=1,2.
\end{equation}

We assume that each physical patch $\Omega^{(k)}$
is given through a NURBS transformation of the parameter domain
$\widehat \Omega$.
Such transformation, denoted by ${\bf F}^{(k)}:\widehat \Omega\to \Omega^{(k)}$,
 is defined by a set of control points
${\bf P}_i^{(k)}\in{\mathbb R}^d$ for $i=1,\ldots, N^{(k)}$, 
thus every point ${\bf x}\in\Omega^{(k)}$ is given by
\begin{equation}\label{map}
{\bf x} = {\bf F}^{(k)}(\widehat{\bf x})=\sum_{i=1}^{N^{(k)}} {\bf P}_i^{(k)}
\, \widehat \varphi^{(k)}_{i,p^{(k)}}
(\widehat{\bf x}).
\end{equation}

Throughout the paper we assume that the 
the mappings ${\bf F}^{(k)}$ are invertible, of class $C^1$ and their 
inverses are of class $C^1$.
After setting  $\widehat\Gamma_k=({\bf F}^{(k)})^{-1}(\Gamma_k)$, we define 
 the space of the traces on $\widehat\Gamma_k$ 
\begin{equation*}
\widehat Y^{(k)}_{h_k}=\left\{\widehat \lambda =\widehat
v|_{\widehat \Gamma_k},\,\,\widehat v\in\widehat {\cal N}^{(k)}_{h_k},
\right\},
\end{equation*}
whose dimension is
$n^{(k)}_\Gamma= (n^{(k)})^{d-1}$, then we denote by
${\bf F}^{(k,\widehat\Gamma_k)}:{\mathbb R}^{d-1}\to {\mathbb R}^{d}$ the
restriction of ${\bf F}^{(k)}$ to $\widehat\Gamma_k$.

The basis functions $\widehat\mu^{(k)}_j$ (with $j=1,\ldots,n^{(k)}_\Gamma$)
of $\widehat Y^{(k)}_{h_k}$  are defined starting from those of
 $\widehat {\cal N}^{(k)}_{h_k}$, more precisely they are the restriction to
$\widehat \Gamma_k$
of those basis functions of $\widehat{\cal N}^{(k)}_{h_k}$ that are not
identically
null on $\widehat \Gamma_k$. Thus, for any basis function
$\widehat\mu^{(k)}_j$ of $\widehat Y^{(k)}_{h_k}$ there exists a unique 
 basis function
$\widehat \varphi^{(k)}_{i_j }\in \widehat {\cal N}^{(k)}_{h_k}$, such that 
\begin{equation}\label{mu-phi}
\widehat \mu_j^{(k)}=(\widehat \varphi^{(k)}_{i_j })|_{\widehat \Gamma_k}.
\end{equation}

Now we define the NURBS function space over the physical domain
$\Omega^{(k)}$ as the \emph{push-forward} of the NURBS function space over
$\widehat \Omega$ through $\bf F^{(k)}$:
\begin{equation}\label{spaces_Nk}
{\cal N}_{h_k}^{(k)}={\cal N}_{h_k}^{(k)}(\Xi^{(k)},p^{(k)})=span
\left\{ \varphi_{i,p^{(k)}}^{(k)}=
\widehat
\varphi_{i,p^{(k)}}^{(k)}\circ({\bf F}^{(k)})^{-1},\ i=1,\ldots,
N^{(k)}\right\}.
\end{equation}

From now on, for sake of clearness, 
we denote the basis functions of ${\cal N}^{(k)}_{h_k}$
 by $\varphi^{(k)}_i $ (instead of $\varphi_{i,p^{(k)}}^{(k)}$).

Starting from (\ref{spaces_Nk}), we define 
the finite dimensional spaces $V^{(k)}_{h_k}$ and $V^{(k)}_{0,h_k}$ as in
(\ref{spaces1}) and the trace space
\begin{equation}
\label{trace_space}
Y^{(k)}_{h_k}=span\{ \mu^{(k)}_j=\widehat\mu^{(k)}_j\circ
({\bf F}^{(k,\widehat\Gamma_k)})^{-1},\ j=1,\ldots,n_\Gamma^{(k)}\}.
\end{equation}

Relations (\ref{mu-phi}) suggest us how to define the discrete counterpart of
the lifting operators ${\cal L}^{(k)}$ invoked in the weak transmission
problem (\ref{transmission_weak}).

For $k=1,2$  we define the linear and continuous discrete lifting operator
$\overline{\cal L}^{(k)}:Y^{(k)}_{h_k}\longrightarrow
{\cal N}^{(k)}_{h_k}$ such that 
$ \overline{\cal L}^{(k)}\mu_j^{(k)}=\varphi^{(k)}_{i_j }$
for any basis function  $\mu_j^{(k)}\in
Y^{(k)}_{h_k}$,
i.e. the lifting of $\mu_j^{(k)}$ is the NURBS basis function
$\varphi^{(k)}_{i_j }$ whose trace on $\Gamma_k$ is $\mu_j^{(k)}$.

\subsection{Interpolation operators}\label{sec:interpolation}

Contrary to mortar methods that are based on $L^2-$projection operators, 
INTERNODES takes advantage of two interpolation operators 
to exchange information between the interfaces $\Gamma_{1}$ and
$\Gamma_2$.
To implement the interpolation process, in this paper
 we use the \emph{Greville abscissae}
(\cite{deboor,abhrs_2010}), but other families of
interpolation nodes could be considered as well (see, e.g.,
\cite{abhrs_2010,montardini_st_cmame17}).

Starting  from the $p-$open knot vector
$\Xi=\left\{\xi_i\right\}_{i=1}^q$ (with 
$q=n+p+1$) in the parameter domain $[0,1]$, 
the \emph{Greville abscissae} (also known as \emph{averaged knot vector}
\cite[Ch. 9]{piegl_nurbs_book}) are defined by
\begin{equation}\label{greville}
\xi_{i,G}= \frac{1}{p}\sum_{j=i+1}^{i+p} \xi_j,\qquad i=1,...,n.
\end{equation}

The assumption that the knot vector $\Xi$ is $p-$open implies that
$\xi_{1,G}=\xi_1=0$ and $\xi_{n,G}=\xi_n=1$.

The Greville abscissae interpolation is  proved to be stable up to degree 3
(\cite{abhrs_2010}), while there are examples of instability for degrees higher
than 19 on particular non-uniform meshes (more precisely, meshes with geometric 
refinement \cite{abhrs_2010,manni_etal_2015}).

For $k=1,2$, we define $\xi_{i,G}^{(k)}\in[0,1]$ as in (\ref{greville}) and,
by tensor product, we build the Greville nodes
\begin{equation}\label{greville_parameter}
\widehat{\bf x}_{i,G}^{(k)}=(\xi_{i_1,G}^{(k)},\ldots,\xi_{i_d,G}^{(k)})\in
\widehat\Omega,
\qquad
\mbox{for }i_1,\ldots,i_d\in\{1,\ldots,n^{(k)}\}.
\end{equation}
The points
\begin{equation}\label{greville_physical}
{\bf x}_{i,G}^{(k)}={\bf F}^{(k)}(\widehat{\bf x}_{i,G}^{(k)})
\end{equation}
are the images  of the Greville nodes in 
the physical patch $\Omega^{(k)}$.
In fact, only the Greville nodes
laying on $\Gamma_k$ will be used during the interpolation process,
these points are denoted by
${\bf x}_{i,G}^{(\Gamma_k)}$. Finally, we denote by
$\widehat{\bf x}_{i,G}^{(\Gamma_k)}=
({\bf F}^{(k,\widehat\Gamma_k)})^{-1}({\bf x}_{i,G}^{(\Gamma_k)})$
the counter-images on $\widehat\Gamma_k$ of the Greville nodes.
Notice that $\widehat{\bf x}_{i,G}^{(\Gamma_k)}\in{\mathbb R}^{d-1}$, while 
${\bf x}_{i,G}^{(\Gamma_k)}\in {\mathbb R}^d$.

We define two different classes of interpolation operators, depending on
the fact that the two interfaces are watertight or not.  The
interpolation we introduce for the watertight case will be generalized to
multipatch decomposition with more than 2 watertight patches,
as in Fig. \ref{fig:geometrie} (c).

\subsection{Interpolation on watertight
interfaces}\label{sec:nurbs_interpolation}

Here the interface $\Gamma_1$ and $\Gamma_2$ 
describe either the same curve in ${\mathbb R}^2$ or the same surface in 
${\mathbb R}^3$, i.e. they coincide at the geometric level with the interface
$\Gamma_{12}$,  but the NURBS spaces $Y_{h_1}^{(1)}$ and $Y_{h_2}^{(2)}$
do not match, that is they are 
unrelated and not necessarily the refinement one of the
other. The set of the Greville
abscissae on $\Gamma_1$ differs a-priori from that of the Greville abscissae
 on $\Gamma_2$.

Given $\lambda^{(1)}\in Y^{(1)}_{h_1}$ and $\lambda^{(2)}\in Y^{(2)}_{h_2}$,
we define
$$
\Pi_{21}:Y^{(1)}_{h_1}\to Y^{(2)}_{h_2},\quad \mbox{ and }\quad
\Pi_{12}:Y^{(2)}_{h_2}\to Y^{(1)}_{h_1}
$$
by the 
interpolation conditions at the Greville nodes laying on $\Gamma_2$ and
$\Gamma_1$, respectively, i.e.,
\begin{eqnarray}\label{Pikl}
\begin{array}{lll}
(\Pi_{21}\lambda^{(1)})({\bf x}_{i,G}^{(\Gamma_2)})=\lambda^{(1)}({\bf
x}_{i,G}^{(\Gamma_2)}), & \mbox{for any }i=1,\ldots, n^{(2)}_\Gamma,\\[2mm]
(\Pi_{12}\lambda^{(2)})({\bf x}_{i,G}^{(\Gamma_1)})=\lambda^{(2)}({\bf
x}_{i,G}^{(\Gamma_1)}), & \mbox{for any }i=1,\ldots, n^{(1)}_\Gamma.
\end{array}
\end{eqnarray}

These interpolation operators are particular instances of those analyzed, e.g.,
 in \cite{demko,piegl_nurbs_book,abhrs_2010}.

Let us suppose that $\lambda^{(1)}\in Y_{h_1}^{(1)}$ is known
 and we want to compute 
the function 
\begin{equation}\label{interpolation}
Y_{h_2}^{(2)}\ni\psi=\Pi_{21}\lambda^{(1)},
\end{equation}
we proceed as
follows. First of all, for $k\in\{1,2\}$ we define the matrices
\begin{equation}\label{matrix_Gkk}
(G_{kk})_{ij}=\widehat\mu_j^{(k)}
(\widehat{\bf x}_{i,G}^{(\Gamma_k)}),\qquad i,j=1,\ldots, n_\Gamma^{(k)}
\end{equation}
and for $\ell\in\{1,2\}$ with $\ell\neq k$ we define the matrices
$G_{k\ell}$ by the relations

\begin{equation}\label{matriciG}
(G_{k\ell})_{ij}=\widehat\mu_j^{(\ell)}
(({\bf F}^{(\ell,\widehat\Gamma_\ell)})^{-1}({\bf x}_{i,G}^{(\Gamma_k)})),
\qquad i=1,\ldots,n_\Gamma^{(k)},\ j=1,\ldots,n_\Gamma^{(\ell)},
\end{equation}

\noindent
i.e., we evaluate the basis functions $\widehat
\mu_j^{(\ell)}$ of the trace space 
$\widehat Y_{h_\ell}^{(\ell)}$ (for $\ell=1,2$) at 
the counter-image w.r.t. the map ${\bf F}^{(\ell,\widehat\Gamma_\ell)}$, 
of the Greville nodes of $\Omega^{(k)}$ laying on $\Gamma_k$ (for $k=1,2$). 
Since $\Gamma_1$ and $\Gamma_2$ describe the same set in the physical
space,  the point-inversion problem
\begin{equation}\label{point-inversion}
find\quad \widehat{\bf x}\in \widehat{\Gamma}_\ell:\ 
({\bf F}^{(\ell,\widehat\Gamma_\ell)})(\widehat{\bf x})=
{\bf x}_{i,G}^{(\Gamma_k)}
\end{equation}
 with 
$\ell\neq k$ has a unique solution (that, e.g., can be numerically computed by 
the Newton method), thus
the entries of the matrices
$G_{k\ell}$ are well defined.

Then, we denote by $\lambda_j^{(1)}$ (for $j=1,\ldots,n^{(1)}_\Gamma$) the
known coefficients of the expansion of
$\lambda^{(1)}$ with respect to the basis function $\mu_j^{(1)}$ of
$Y^{(1)}_{h_1}$ and by $\psi_j$ the unknown
coefficients of the expansion of $\psi$ with respect
to the basis functions of $Y_{h_2}^{(2)}$. In view of (\ref{trace_space}),
for any ${\bf x}\in \Gamma_2$
it holds
\begin{eqnarray}
\begin{array}{l}
\lambda^{(1)}({\bf x})=
\displaystyle\sum_{j=1}^{n^{(1)}_\Gamma}\lambda_j^{(1)}\mu_j^{(1)}({\bf x})
= \sum_{j=1}^{n^{(1)}_\Gamma}\lambda_j^{(1)}
\widehat\mu_j^{(1)}(({\bf F}^{(1,\widehat\Gamma_1)})^{-1}({\bf x})),\\
\psi({\bf x})=
\displaystyle\sum_{j=1}^{n^{(2)}_\Gamma}\psi_j\mu_j^{(2)}({\bf x})
= \sum_{j=1}^{n^{(2)}_\Gamma}\psi_j
\widehat\mu_j^{(2)}(({\bf F}^{(2,\widehat\Gamma_2)})^{-1}({\bf x}))
.
\end{array}
\end{eqnarray}
From now on, 
the expansion of a function  of $Y_{h_k}^{(k)}$
with respect to the basis functions $\mu_j^{(k)}$ is named \emph{primal}
and the associated coefficients are named \emph{primal coefficients}.

Thanks to both (\ref{matrix_Gkk}) and (\ref{matriciG}),
the interpolation conditions $\psi({\bf
x}_{i,G}^{(\Gamma_2)})=\lambda^{(1)}({\bf x}_{i,G}^{(\Gamma_2)})$ 
 (i.e. (\ref{Pikl})$_1$)
 read:
\begin{equation}\label{interp_bis}
\sum_{j=1}^{n^{(2)}_\Gamma}\psi_j\underbrace{\widehat\mu_j^{(2)}
(({\bf F}^{(2,\widehat\Gamma_2)})^{-1}({\bf x}_{i,G}^{(\Gamma_2)})) }_{(G_{22})_{ij}}=
\sum_{j=1}^{n^{(1)}_\Gamma}\lambda_j^{(1)}
\underbrace{\widehat\mu_j^{(1)}
(({\bf F}^{(1,\widehat\Gamma_1)})^{-1}({\bf x}_{i,G}^{(\Gamma_2)}))}_{(G_{21})_{ij}}, \qquad
i=1,\ldots,n^{(2)}_\Gamma.
\end{equation}

Denoting by $\boldsymbol\lambda^{(1)}$ ($\boldsymbol\psi$, resp.)
the array whose components are 
the values $\lambda_j^{(1)}$ ($\psi_j$, resp.), 
(\ref{interp_bis}) becomes
\begin{equation}\label{G22_nurbs}
G_{22}\boldsymbol\psi=G_{21}\boldsymbol\lambda^{(1)}.
\end{equation}

In conclusion, 
given $\boldsymbol\lambda^{(1)}$, we compute $\boldsymbol\psi$ by
\begin{equation}\label{matrix_P21}
\boldsymbol\psi=P_{21}\boldsymbol\lambda^{(1)}, 
\quad\mbox{ with }\quad P_{21}=G_{22}^{-1}G_{21}.
\end{equation}

The matrix $P_{21}$ (with $n^{(2)}_\Gamma$ rows and 
$n^{(1)}_\Gamma$ columns) implements the interpolation
operator $\Pi_{21}$. 

Proceeding in a similar way, we define the matrix 
\begin{equation}\label{matrix_P12}
P_{12}=G_{11}^{-1}G_{12}
\end{equation}
that is the algebraic counterpart of the interpolation operator 
$\Pi_{12}$.

The matrices $G_{kk}$ (for $k=1,2$) are non-singular 
(see \cite[Ch. 9.2]{piegl_nurbs_book}), thus we can compute $P_{k\ell}$ by 
solving the linear system $G_{kk} P_{k\ell}=G_{k\ell}$ by any suitable method. Notice
that the size of $G_{kk}$ is equal to the number of degrees of freedom on the
interface $\Gamma_k$, then it is considerably less than the number of degrees of freedom in
$\Omega^{(k)}$. Moreover, the computation of $P_{k\ell}$ is done only once 
during
the initialization step of INTERNODES, and only matrix-vector products
between $P_{k\ell}$ and the array of the degrees of freedom on $\Gamma_{\ell}$ 
are required at the successive steps of INTERNODES (see Algorithm
\ref{alg:Slambda}).

\begin{remark}\label{rem:Piu}
Notice that 
$\Pi_{21}u_{h_1}^{(1)}$ in (\ref{internodes_weak})$_2$
stands for  $\Pi_{21}(u_{h_1}^{(1)})|_{\Gamma_1}$.
\end{remark}

\subsection{RL-RBF interpolation on non-watertight
interfaces}\label{sec:rbf_interpolation}

In this Section we face the case for which
 $\Gamma_1$ and $\Gamma_2$ are only approximation of the common boundary
$\Gamma_{12}$,
 as e.g.  in Fig. \ref{fig:geometrie} (a) and for which the point-inversion
problem (\ref{point-inversion}) may have no solution. Typically in mortar
method this drawback is faced by projecting the point ${\bf
x}_{i,G}^{(\Gamma_k)}$ onto $\widehat \Gamma_\ell$. Here we overcome the
problem by using the \emph{Rescaled Localized  Radial Basis
Function (RL-RBF)} interpolation operators introduced in
formula (3.1) of \cite{dfq_rbf}.  

\begin{figure}
\begin{center}
\includegraphics[trim={0 40 0 0},width=0.25\textwidth]{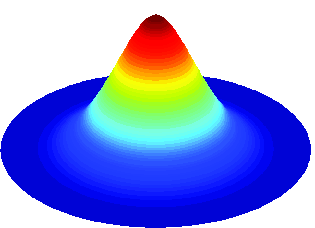}
\qquad 
\scalebox{0.6}{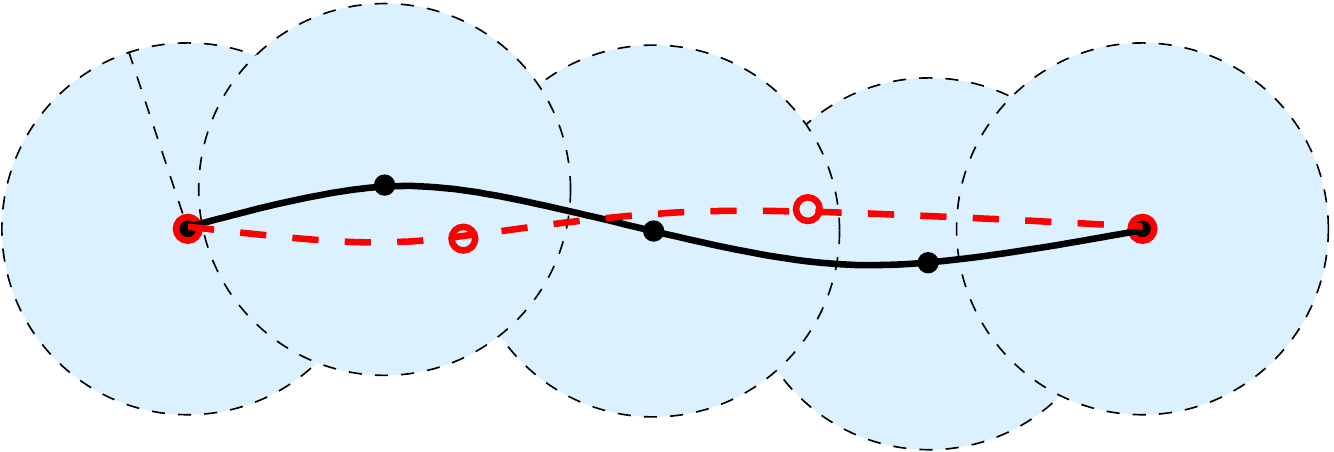}
\end{center}
\caption{At left, an example of locally supported $C^2-$ Wendland radial basis function
$\tilde\phi^{(k)}_j({\bf x})$ when $d=2$. At right, the interfaces 
$\Gamma_1$ (continuous line) and $\Gamma_2$
(dashed line), the Greville abscissae ${\bf x}_{i,G}^{(\Gamma_1)}$ (full
circles), 
and  ${\bf x}_{j,G}^{(\Gamma_2)}$ (empty circles) and the support
(the union of the light blue circles) of the RL-RBF interpolant 
$\Pi_{RL-RBF}^{(1)}\lambda$}
\label{fig:Wendland}
\end{figure}

More precisely, for $k=1,2$,
for $i=1,\ldots,n^{(k)}_\Gamma$ and for any ${\bf x}\in {\mathbb R}^d$
($d=2,3$) let
$$\tilde\phi^{(k)}_j({\bf x})=\phi(\|{\bf x}-{\bf x}_{j,G}^{(\Gamma_k)}\|,r)=
\max\left\{0,\left(1-\frac{\|{\bf x}-{\bf x}_{j,G}^{(\Gamma_k)}\|}{r}\right)^4\right\}
\left(1+4\frac{\|{\bf x}-{\bf x}_{j,G}^{(\Gamma_k)}\|}{r}\right)$$
 be the locally supported $C^2-$ Wendland radial basis function
(\cite{wendland1995piecewise}) 
centered at ${\bf x}_{j,G}^{(\Gamma_k)}$ with radius
$r>0$ (see Fig. \ref{fig:Wendland}), where  $\|\cdot\|$ denotes the
euclidean norm in ${\mathbb R}^d$.

If $\lambda$ is any continuous function defined on $\Gamma_k$, 
the Radial Basis Function
(RBF) interpolating $\lambda$ at the nodes ${\bf x}_{i,G}^{(\Gamma_k)}$, with 
$i=1,\ldots, n_\Gamma^{(k)}$, reads
\begin{equation} \label{rbf}
(\Pi^{(k)}_{RBF}\lambda)({\bf x})=
\sum_{j=1}^{n_\Gamma^{(k)}}(\gamma_\lambda^{(k)})_j\tilde\phi^{(k)}_j({\bf x}),
\end{equation}

\noindent
where the real values $(\gamma_\lambda^{(k)})_j$ are the solutions of the 
$n_\Gamma^{(k)}\times n_\Gamma^{(k)}$ linear system 
\begin{equation}\label{rbf-coeff}
(\Pi^{(k)}_{RBF}\lambda)({\bf x}_{i,G}^{(\Gamma_k)})=
\sum_{j=1}^{n_\Gamma^{(k)}}(\gamma_\lambda^{(k)})_j\tilde\phi^{(k)}_j({\bf x}_{i,G}^{(\Gamma_k)})
=\lambda({\bf x}_{i,G}^{(\Gamma_k)}), \qquad i=1,\ldots, n_\Gamma^{(k)}.
\end{equation}

After setting $g(x)\equiv 1$, the RL-RBF interpolant 
of $\lambda$ at the nodes
${\bf x}_{i,G}^{(\Gamma_k)}$  reads (\cite{dfq_rbf}):
\begin{equation}\label{rl-rbf}
(\Pi_{RL-RBF}^{(k)} \lambda)({\bf x})=\frac{(\Pi_{RBF}^{(k)} \lambda)({\bf x})}
{(\Pi_{RBF}^{(k)} g)({\bf x})}=
\frac{\sum_{j=1}^{n_\Gamma^{(k)}}(\gamma_\lambda^{(k)})_j\tilde\phi^{(k)}_j({\bf x})}
{\sum_{j=1}^{n_\Gamma^{(k)}}(\gamma_g^{(k)})_j\tilde\phi^{(k)}_j({\bf x})}.
\end{equation}

The advantage of RL-RBF interpolant (\ref{rl-rbf}) with respect to 
 (\ref{rbf}) is that (\ref{rl-rbf}) reproduces 
exactly constant functions, while  (\ref{rbf}) does not (\cite{dfq_rbf}).

Notice that $\Pi_{RL-RBF}^{(k)} \lambda$ is defined in ${\mathbb R}^d$
and not only on the  $(d-1)-$dimensional manifold $\Gamma_k$ of ${\mathbb
R}^d$, moreover
its support depends on the chosen radius $r$ (see, Fig. \ref{fig:Wendland},
right). We refer to \cite{dfq_rbf} for the discussion about the optimal choice
of the radius $r$ and the accuracy of RL-RBF interpolation.

For $k,\ell\in\{1,2\}$ we define the matrices
$$(\Phi_{k\ell})_{ij}=\tilde\phi_j^{(\ell)}({\bf x}_{i,G}^{(\Gamma_k)})\quad
\mbox{for }k,\ell\in\{1,2\},\ i=1,\ldots, n_\Gamma^{(k)},
\ j=1,\ldots,  n_\Gamma^{(\ell)}$$
and the RL-RBF interpolation matrices 
\begin{equation}\label{P-RBF}
(P^{RBF}_{\ell k})_{ij}=
\frac{(\Phi_{\ell k}\Phi_{kk}^{-1})_{ij}}
{(\Phi_{\ell k}\Phi_{kk}^{-1}{\bf 1})_{i}}, \qquad  i=1,\ldots, n_\Gamma^{(\ell)},
\ j=1,\ldots, n_\Gamma^{(k)},
\end{equation}
where ${\bf 1}$ denotes a column array with all entries  equal to 1.

The evaluation of the RL-RBF interpolant
$\Pi_{RL-RBF}^{(k)}\lambda$ at the Greville nodes ${\bf
x}_{i,G}^{(\Gamma_\ell)}$ is given by the matrix vector product
\begin{eqnarray}\label{rbf_alg}
(\Pi_{RL-RBF}^{(k)} \lambda)({\bf x}_{i,G}^{(\Gamma_\ell)}) =
\sum_{j=1}^{n_\Gamma^{(k)}}P_{\ell k}^{RBF}\lambda({\bf x}_{j,G}^{(\Gamma_k)}),
\qquad i=1,\ldots, n_\Gamma^{(\ell)}
\end{eqnarray}
Notice that $P_{\ell k}^{RBF}$ is applied to a vector of nodal values and
produces a vector of nodal values. 

\medskip
Now we exploit the RL-RBF matrices (\ref{P-RBF}) 
to build the INTERNODES interpolation operators for non-watertight interfaces.
Given $\lambda^{(1)}\in Y_{h_1}^{(1)}$ and 
$\lambda^{(2)}\in Y_{h_2}^{(2)}$  we define the interpolation operators
$$
\Pi_{21}:Y^{(1)}_{h_1}\to Y^{(2)}_{h_2},\quad \mbox{ and }\quad
\Pi_{12}:Y^{(2)}_{h_2}\to Y^{(1)}_{h_1}
$$
as follows:

given $\lambda^{(1)}\in Y_{h_1}^{(1)}$, the function 
$\Pi_{21}\lambda^{(1)}\in Y^{(2)}_{h_2}$ is the NURBS defined on 
$\Gamma_2$ that interpolates the RL-RBF interpolant $\Pi_{RL-RBF}^{(1)}
\lambda^{(1)}$ at the Greville abscissae ${\bf x}_{i,G}^{(\Gamma_2)}$, that is
\begin{eqnarray}\label{PI21_rbf}
(\Pi_{21}\lambda^{(1)})({\bf x}_{i,G}^{(\Gamma_2)})=
(\Pi_{RL-RBF}^{(1)} \lambda^{(1)})({\bf x}_{i,G}^{(\Gamma_2)}), \quad
i=1,\ldots,n_\Gamma^{(2)},
\end{eqnarray}
and similarly, given $\lambda^{(2)}\in Y_{h_2}^{(2)}$, the function
$\Pi_{12}\lambda^{(2)}\in Y^{(1)}_{h_1}$ is the NURBS defined on 
$\Gamma_1$ that interpolates the RL-RBF interpolant $\Pi_{RL-RBF}^{(2)}
\lambda^{(2)}$ at the Greville abscissae ${\bf x}_{i,G}^{(\Gamma_1)}$, that is
\begin{eqnarray}\label{PI12_rbf}
(\Pi_{12}\lambda^{(2)})({\bf x}_{i,G}^{(\Gamma_1)})=
(\Pi_{RL-RBF}^{(2)} \lambda^{(2)})({\bf x}_{i,G}^{(\Gamma_1)}), \quad
i=1,\ldots,n_\Gamma^{(1)}.
\end{eqnarray}

We show how to compute 
$Y_{h_2}^{(2)}\ni\psi=\Pi_{21}\lambda^{(1)}=\Pi_{RL-RBF}^{(1)} \lambda^{(1)}$
by  (\ref{PI21_rbf}).
Denoting by $\{\lambda^{(1)}_j\}$ the coefficients of $\lambda^{(1)}$ with
respect to the NURBS basis functions of $Y_{h_1}^{(1)}$ and by 
$\psi_i$ the coefficients of $\psi$ with
respect to the NURBS basis functions of $Y_{h_2}^{(2)}$, 
the idea surrounding the interpolation operator $\Pi_{21}$ 
can be resumed by the following diagram:
\begin{eqnarray*}
\displaystyle\{\lambda^{(1)}_j\}
\xrightarrow[G_{11}]{}
\{\lambda^{(1)}({\bf x}_{j,G}^{(\Gamma_1)})\}
\xrightarrow[P^{RBF}_{21}]{}
\{(\Pi_{RL-RBF}^{(1)}\lambda^{(1)})({\bf x}_{i,G}^{(\Gamma_2)})
=\psi({\bf x}_{i,G}^{(\Gamma_2)})\}
\xrightarrow[G_{22}^{-1}]{}
\{\psi_i\}.
\end{eqnarray*}

Since
 $\lambda^{(1)}({\bf x}_{j,G}^{(\Gamma_1)})=(G_{11} \boldsymbol\lambda^{(1)})_j$
and $\psi({\bf x}_{i,G}^{(\Gamma_2)})=(G_{22}\boldsymbol\psi)_i$, the
algebraic counterpart of (\ref{PI21_rbf}) reads
\begin{equation}\label{G22_rbf}
G_{22}\boldsymbol\psi=P^{RBF}_{21} G_{11}\boldsymbol\lambda^{(1)}.
\end{equation}
Thus, given $\boldsymbol\lambda^{(1)}$, we compute 
$\boldsymbol\psi$ by
\begin{equation}\label{matrix_P21_RBF}
\boldsymbol\psi=P_{21}\boldsymbol\lambda^{(1)}, 
\quad\mbox{ with }\quad P_{21}=G_{22}^{-1}P^{RBF}_{21}G_{11}.
\end{equation}

We notice the difference between (\ref{G22_nurbs})  and (\ref{G22_rbf}): 
the matrix $G_{21}$ in (\ref{G22_nurbs}) is replaced by the product
$P^{RBF}_{21}G_{11}$ in (\ref{G22_rbf}) 
that allows the transfer of information between the two
non-watertight interfaces, by circumventing the projection task.

Proceeding in a similar way, we define the RL-RBF interpolation matrix
from $\Gamma_2$ to $\Gamma_1$:
\begin{equation}\label{matrix_P12_RBF}
P_{12}=G_{11}^{-1}P^{RBF}_{12}G_{22}.
\end{equation}

\begin{remark} 
From now on, when the interfaces will be watertight we will implement
INTERNODES by using the interpolation matrices (\ref{matrix_P21}) and
(\ref{matrix_P12}), while in presence of non-watertight interfaces we will use
(\ref{matrix_P21_RBF}) and (\ref{matrix_P12_RBF}).
\end{remark}

\subsection{Transferring the normal derivatives}
\label{sec:normalder}

We start by underlying that the interpolation matrix 
$P_{\ell k}$ 
defined in either (\ref{matrix_P21})--(\ref{matrix_P12}) or
 (\ref{matrix_P21_RBF})--(\ref{matrix_P12_RBF})
 applies to the primal coefficients of the functions
belonging to $Y_{h_k}^{(k)}$.
Despite this, the normal derivative
$ \frac{\partial u^{(k)}_{h_k}}{\partial 
{\bf n}_k}$ is a functional belonging to the space $(Y_{h_k}^{(k)})'$ dual of 
$Y_{h_k}^{(k)}$, thus we have to understand how to apply the interpolation to
 the normal
derivatives and we have to explain which is the meaning of
$\widetilde\Pi_{12}\frac{\partial u^{(2)}_{h_2}}{\partial 
{\bf n}_2}$ as it appears 
in the interface condition (\ref{internodes_weak})$_3$.

For $k=1,2$ and 
for any basis function $\mu_i^{(k)}$ of $Y_{h_k}^{(k)}$
we define the real values
\begin{equation}\label{residual}
r^{(k)}_i=\langle \frac{\partial u^{(k)}_{h_k}}{\partial 
{\bf n}_k},\mu_i^{(k)} \rangle
=\int_{\Gamma_k}\frac{\partial u^{(k)}_{h_k}}{\partial {\bf n}_k}\,
\mu_i^{(k)}d\Gamma , \qquad i=1,\ldots,n^{(k)}_\Gamma.
\end{equation}

When dealing with the weak form (\ref{transmission_weak}) of the transmission
problem, to compute $r^{(k)}_i$ with a small effort, 
we can work as follows.
Let $u^{(k)}_{h_k}$ be any function in ${\cal N}^{(k)}_{h_k}$ such that
\begin{equation}\label{eqweak}
a^{(k)}(u^{(k)}_{h_k},v_{h_k})={\cal F}^{(k)} (v_{h_k})+
\int_{\partial\Omega^{(k)}}\frac{\partial u^{(k)}_{h_k}}{\partial 
{\bf n}_k}v_{h_k}\,d\Gamma \qquad \forall v_{h_k}\in {\cal N}^{(k)}_{h_k},
\end{equation}
and 
let us set $G_k=\partial\Omega^{(k)}\setminus\Gamma_k$.
Then the values $r_i^{(k)}$ read also
\begin{equation}\label{flux}
r^{(k)}_i=a^{(k)}(u^{(k)}_{h_k},\overline{\cal L}^{(k)}\mu_i^{(k)})-{\cal F}^{(k)}
(\overline{\cal L}^{(k)}\mu_i^{(k)})
-\int_{G_k} \frac{\partial u^{(k)}_{h_k}}{\partial {\bf n}_k}
\,\overline{\cal L}^{(k)}\mu_i^{(k)}\,d\Gamma.
\end{equation}

The algebraic implementation of (\ref{flux}) is in fact a
matrix-vector product between the stiffness matrix and the array of the degrees
of freedom associated with the interface (see
(\ref{matrix_corr})--(\ref{res_modified})). Obviously, in 
the case that the strong
form of the transmission problems
 is preferred instead of the Galerkin one, the formula
(\ref{residual}) must be considered instead of (\ref{flux}).

\begin{figure}
\begin{center}
\includegraphics[trim={20 80 20
40},width=0.43\textwidth]{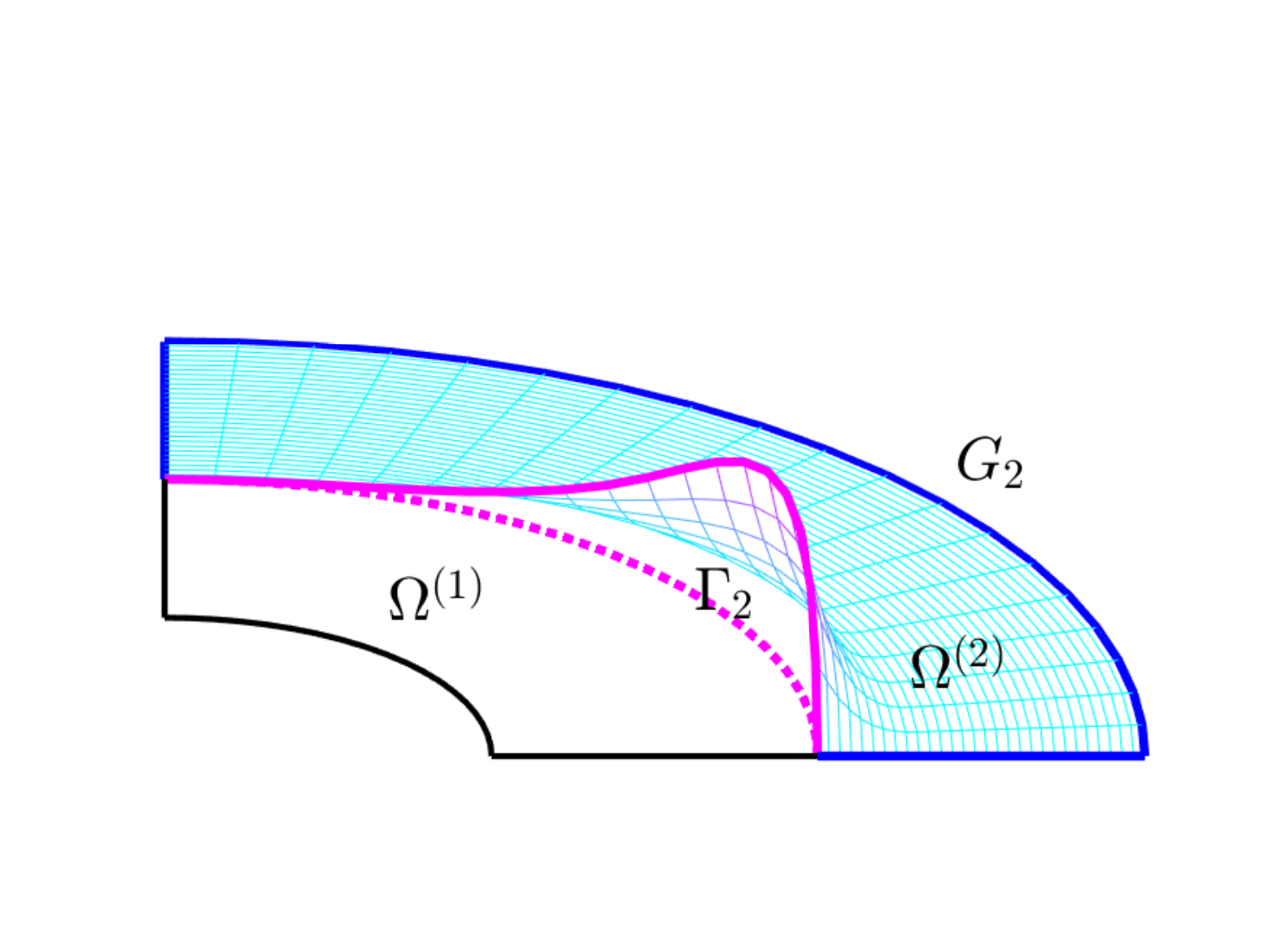}\quad
\includegraphics[trim={20 60 20
40},width=0.43\textwidth]{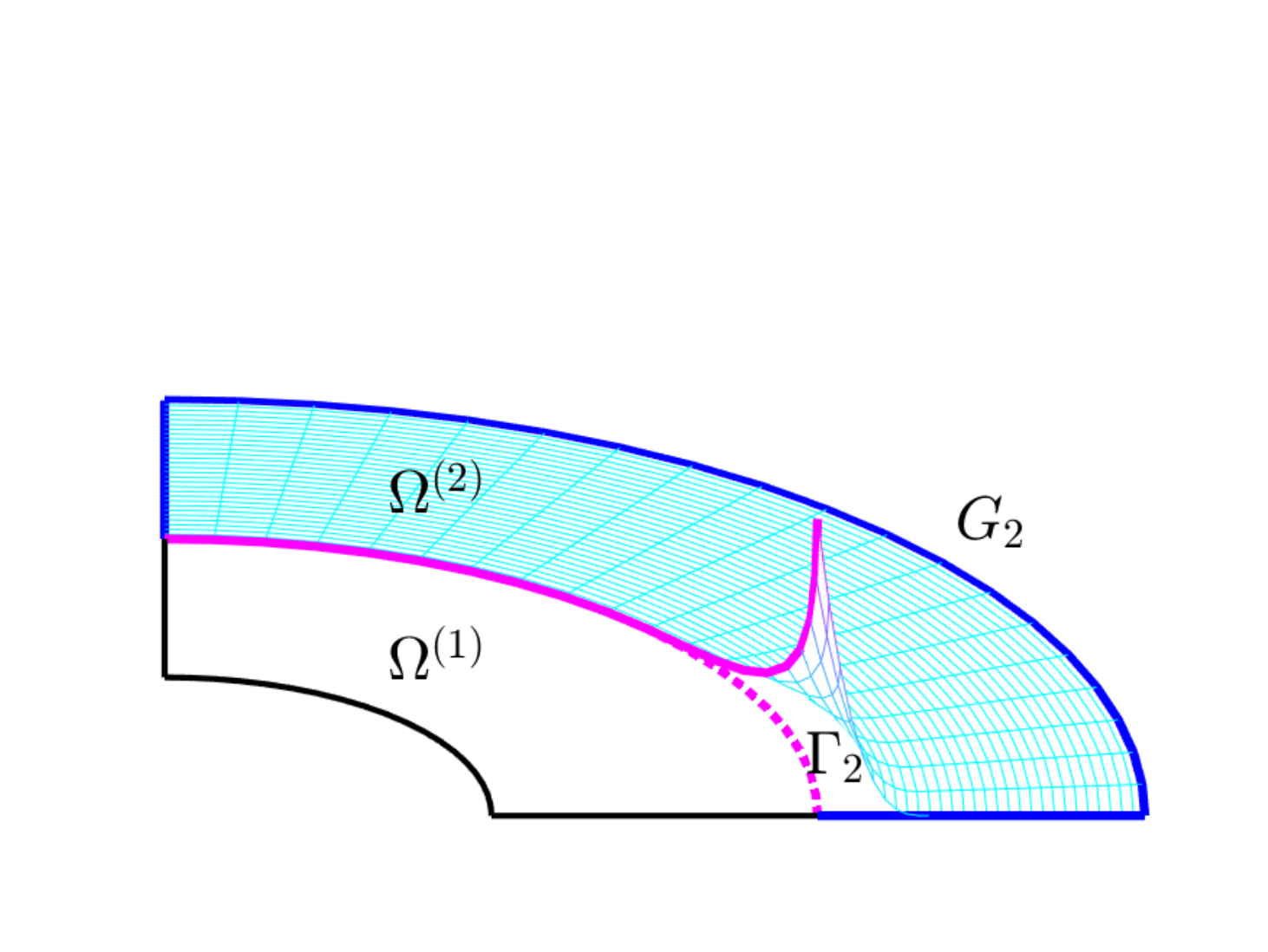}\quad
\end{center}
\caption{The lifting $\overline{\cal L}^{(2)}\mu_j^{(2)}$. At left,
$\mu_j^{(2)}$ (the purple function defined on $\Gamma_2$)
is identically null on $\partial\Gamma_2$ and $\overline{\cal
L}^{(2)}\mu_j^{(2)}$ is
identically null on $G_2=\partial\Omega^{(2)}\setminus\Gamma_2$ ($G_2$
is drawn in blue);
at right,  $\mu_j^{(2)}$ is not identically null on $\partial\Gamma_2$ and
 $\overline{\cal
L}^{(2)}\mu_j^{(2)}$
 is not identically null on $G_2$. }
\label{fig:lifting_R}
\end{figure}

The presence of the 
last term in (\ref{flux}) is justified as follows. We notice that when
$\mu_j^{(k)}$ is not identically null on $\partial\Gamma_k$,
then  $\overline{\cal L}^{(k)}\mu_j^{(k)}$ is not identically null on
the set $G_k$ (see Fig. \ref{fig:lifting_R}).

For $k=1,2$ we
 denote by ${\bf r}_{\overline\Gamma_k}$ the array whose entries are
the real values $r_i^{(k)}$, for $i=1,\ldots, n^{(k)}_\Gamma$ and,
following the nomenclature used in linear algebra,
${\bf r}_{\overline\Gamma_k}$ is named \emph{residual vector}.

The values $r^{(k)}_i$ are not the coefficients of the primal expansion of
$\frac{\partial u_{h_k}^{(k)}}{\partial{\bf n}_k}$,
so we cannot apply the interpolation matrix $P_{\ell k}$ to the array ${\bf
r}_{\overline\Gamma_k}$. Rather they
are the coefficients of $\frac{\partial u_{h_k}^{(k)}}{\partial{\bf n}_k}$
with respect to the 
basis $\{\Phi_j^{(k)}\}_{j=1}^{n^{(k)}_\Gamma}$  in
$(Y_{h_k}^{(k)})'$ that is \emph{dual} to 
$\{\mu_j^{(k)}\}_{j=1}^{n^{(k)}_\Gamma}$ (see, e.g.
\cite{gq_internodes,brauchli_oden}), i.e. satisfying
\begin{equation*}
\langle \Phi_j^{(k)},\mu_i^{(k)}\rangle=\int_{\Gamma_k}\Phi_j^{(k)}\,
\mu_i^{(k)}d\Gamma=\delta_{ij}, \qquad i,j=1,\ldots,n^{(k)}_\Gamma
\end{equation*}
(where $\delta_{ij}$ is the Kronecker delta), and it holds
\begin{equation}\label{residual_expansion_dual}
\frac{\partial u_{h_k}^{(k)}}{\partial{\bf n}_k} 
=\sum_{j=1}^{n^{(k)}_\Gamma}r^{(k)}_j\Phi_j^{(k)}.
\end{equation}


Nevertheless, $Y_{h_k}^{(k)}$ and $(Y_{h_k}^{(k)})'$ are the same 
(finite dimensional)
algebraic space (\cite{brauchli_oden}) and we denote by ${\cal J}_k$ 
the canonical isomorphism between 
 $Y_{h_k}^{(k)}$ and its dual $(Y_{h_k}^{(k)})'$.

To transfer the normal derivative 
$\frac{\partial u^{(2)}_{h_2}}{\partial {\bf n}_2}$ from
$\Gamma_2$ to $\Gamma_1$,
 we define 
the operator $\widetilde \Pi_{12}:(Y_{h_2}^{(2)})'\to
(Y_{h_1}^{(1)})'$ such that
$\widetilde \Pi_{12}={\cal J}_{1}\Pi_{12}{\cal J}_2^{-1},$
i.e.
\begin{displaymath}
\xymatrix{
\widetilde\Pi_{12}:\  (Y_{h_2}^{(2)})' \ar[r]_{\ \ {\cal J}_2^{-1}}
& Y_{h_2}^{(2)} \ar[r]_{\Pi_{12}}  
& Y_{h_1}^{(1)} \ar[r]_{{\cal J}_1}
& (Y_{h_1}^{(1)})'.}
\end{displaymath}

The \emph{interface mass matrix} 
$M_{\Gamma_k}$ on $\Gamma^{(k)}$, whose entries are
\begin{equation}\label{local_mass}
(M_{\Gamma_k})_{ij}=(\mu_j^{(k)},\mu_i^{(k)})_{L^2(\Gamma_k)}, 
\qquad i,j=1,\ldots, n^{(k)}_\Gamma,
\end{equation}
is the matrix corresponding to the isomorphism ${\cal J}_k$.

Assume that $u^{(2)}_{h_2}$ is known, 
the computation of $\widetilde\Pi_{12}\frac{\partial
u_{h_2}^{(2)}}{\partial{\bf n}_2}$  is carried out as
follows:
\begin{enumerate}
\item compute the entries of the array ${\bf r}_{\overline\Gamma_2}$ by
(\ref{flux});
\item compute the array
$
{\bf z}_{\overline\Gamma_2}=M_{\Gamma_2}^{-1}{\bf r}_{\overline\Gamma_2},
$
the entries of ${\bf z}_{\overline\Gamma_2}$ 
are in fact the coefficients of the expansion (named $z_{h_2}^{(2)}$) of  
$\frac{\partial u_{h_2}^{(2)}}{\partial{\bf n}_2}$
with respect to  the primal basis $\mu_j^{(2)}$ of $Y_{h_2}^{(2)} $
 (see, e.g., \cite{gq_internodes,brauchli_oden});
\item compute  the array
${\bf s}_{\overline\Gamma_1}=P_{12}{\bf z}_{\overline\Gamma_2},$
the entries of ${\bf s}_{\overline\Gamma_1}$ 
are the primal coefficients of the function $s=\Pi_{12}z_{h_2}^{(2)}\in
Y_{h_1}^{(1)}$;
\item compute the array
$
\breve{\bf r}_{\overline\Gamma_1}=M_{\overline\Gamma_1}{\bf
s}_{\overline\Gamma_1},
$
i.e., come back to the dual expansion.
\end{enumerate}

The entries of the array
\begin{equation}\label{interp_z}
\breve{\bf r}_{\overline\Gamma_1}=
M_{\Gamma_1}P_{12}M_{\Gamma_2}^{-1}{\bf r}_{\overline\Gamma_2}
\end{equation}
are the coefficients of the expansion of $\widetilde\Pi_{12}
\frac{\partial u^{(2)}_{h_2}}{\partial {\bf n}_2}$ with respect
to the dual basis $\{\Phi_j^{(1)}\}$, i.e.,
\begin{equation}\label{Pitilde}
\widetilde\Pi_{12}\displaystyle \frac{\partial u^{(2)}_{h_2}}{\partial 
{\bf n}_2}=\sum_{j=1}^{n_\Gamma^{(1)}}\breve r_j\Phi_j^{(1)}, \quad
\mbox{ and }\quad 
\langle \widetilde\Pi_{12}\displaystyle \frac{\partial u^{(2)}_{h_2}}{\partial 
{\bf n}_2},\mu_i^{(1)}\rangle=\breve r_i.
\end{equation}

In conclusion, in view of (\ref{residual_expansion_dual}) and (\ref{Pitilde}),
the algebraic counterpart of the interface condition
(\ref{internodes_weak})$_3$ reads

\begin{equation}\label{residual_equil_alg}
{\bf r}_{\overline\Gamma_1}+M_{\Gamma_1}P_{12}M_{\Gamma_2}^{-1}{\bf
r}_{\overline\Gamma_2}={\bf 0}.
\end{equation}

 \section{The algebraic form of INTERNODES}\label{sec:alg}

For $k=1,2$ we define the following sets of indices:
\begin{description}[noitemsep]
 \item[--] ${\cal I}_{\overline\Omega^{(k)}}=\{1,\ldots,N^{(k)}\}$;
 \item[--] ${\cal I}_k$ 
the subset of the indices of ${\cal I}_{\overline\Omega^{(k)}}$
associated with the basis functions of ${\cal N}^{(k)}_{h_k}$ that
 are identically null on $\partial \Omega^{(k)}$;
\item[--] ${\cal I}_{\overline \Gamma_k}$ 
the subset of the indices of ${\cal I}_{\overline\Omega^{(k)}}$
associated with the basis functions of ${\cal N}^{(k)}_{h_k}$ that
are not identically null
on $\Gamma_k$ (even if $\Gamma_k=\overline\Gamma_k$, the 
bar over $\Gamma_k$ stresses the fact that we are taking into account for
all the basis functions that are not identically null on $\overline\Gamma_k$);
\item[--] ${\cal I}_{\Gamma_k}$ 
the subset of the indices of ${\cal I}_{\overline\Omega^{(k)}}$
associated with the basis functions of ${\cal N}^{(k)}_{h_k}$ that
are not identically null on the interior of $\Gamma_k$ and are identically null
 on $\partial\Gamma_k$;
\item[--] ${\cal I}_{D_k}$ 
the subset of the indices of ${\cal I}_{\overline\Omega^{(k)}}$
associated with the Dirichlet degrees
of freedom;
\item[--]
${\cal I}_{\partial\Gamma_k}={\cal I}_{\overline\Gamma_k}\setminus{\cal
I}_{\Gamma_k}$.
\end{description}

We define the local stiffness matrices  $A^{(k)}$ whose entries are
\begin{equation}\label{matrices}
A^{(k)}_{ij}=a^{(k)}(\varphi_j,\varphi_i), \qquad i,j\in {\cal
I}_{\overline\Omega^{(k)}},
\end{equation}
then let 
$$
A^{(k,k)}=A^{(k)}({\cal I}_k,{\cal I}_k)$$
be the submatrix of $A^{(k)}$ obtained by taking both rows and columns of 
$A^{(k)}$ whose indices belong to ${\cal I}_k$.
Similarly, we define the submatrices
$A^{(k,\Gamma_k)}=A^{(k)}({\cal I}_{k},{\cal I}_{\Gamma_k})$,
$A^{(\Gamma_k,\Gamma_k)}=A^{(k)}({\cal I}_{\Gamma_k},{\cal I}_{\Gamma_k}),$
$A^{(\overline\Gamma_k,k)}=A^{(k)}({\cal I}_{\overline\Gamma_k},{\cal I}_{k})$,
$A^{(\overline\Gamma_k,D_k)}=A^{(k)}({\cal I}_{\overline\Gamma_k},{\cal
I}_{D_k})$, 
$A^{(\overline\Gamma_k,\overline\Omega^{(k)})}=A^{(k)}({\cal I}_{\overline\Gamma_k},
{\cal I}_{\overline\Omega^{(k)}})$, 
and so on.

Moreover, we define the array ${\bf f}^{(k)}$ 
whose entries are
\begin{eqnarray*}
f^{(k)}_i={\cal F}_k(\varphi_i^{(k)}),\qquad  i\in {\cal
I}_{\overline\Omega^{(k)}},
\end{eqnarray*}
the array ${\bf u}^{(k)}$ of the degrees of freedom in
$\overline\Omega^{(k)}$, and using the same notation as above, the subarrays
\begin{eqnarray*}
\begin{array}{lll}
{\bf f}^{(k)}_0={\bf f}^{(k)}({\cal I}_k), & 
{\bf f}_{\overline\Gamma_k}= {\bf f}^{(k)}({\cal I}_{\overline \Gamma_k}), &
{\bf f}_{\Gamma_k}= {\bf f}^{(k)}({\cal I}_{\Gamma_k}), \\[1mm]
{\bf u}^{(k)}_0={\bf u}^{(k)}({\cal I}_k), & 
{\bf u}_{\overline\Gamma_k}= {\bf u}^{(k)}({\cal I}_{\overline \Gamma_k}), &
{\bf u}_{\Gamma_k}= {\bf u}^{(k)}({\cal I}_{\Gamma_k}).
\end{array} 
\end{eqnarray*}

Finally, we denote by ${\bf g}_{k}$ and  ${\bf g}_{\partial\Gamma_k}$
the arrays of all the Dirichlet degrees of
freedom associated with $\partial\Omega^{(k)}_D$ and
$\partial\Gamma_k$, respectively.

To evaluate the last integral of (\ref{flux})
 we define the matrix $C^{(k)}$ whose non-null entries are
\begin{eqnarray}\label{matrix_corr}
\begin{array}{lll}
C^{(k)}_{ij}&
 \displaystyle =-\int_{G_k} \frac{\partial\varphi_j ^{(k)}}
{\partial{\bf n}_k}
\varphi_i^{(k)}, & \mbox{ for }i\in{\cal I}_{\partial\Gamma_k},\ j\in{\cal
I}_{\overline\Omega^{(k)}}
\end{array}
\end{eqnarray}
and, as done for the stiffness matrix $A^{(k)}$, we set
$C^{(\Gamma_k,k)}=C^{(k)}({\cal I}_{\Gamma_k},{\cal I}_{k}),$
$C^{(\Gamma_k,\Gamma_k)}=C^{(k)}({\cal I}_{\Gamma_k},{\cal I}_{\Gamma_k}),$
$C^{(\Gamma_k,D_k)}=C^{(k)}({\cal I}_{\Gamma_k},{\cal I}_{D_k})$,
and so on.

The integrals in (\ref{matrix_corr}) can be easily computed by exploiting the
definition of the NURBS basis functions, moreover
the rows of $C^{(k)}$ associated with all the degrees of freedom not
belonging to ${\cal I}_{\partial\Gamma_k}$ are null. Thus the computation of
$C^{(k)}$ is very cheap.

Then we define
\begin{equation}\label{def:AC}
\widehat{A}^{(\Gamma_k,X)}=A^{(\Gamma_k,X)}+C^{(\Gamma_k,X)} \quad \mbox{and}
\quad
\widehat{A}^{(\overline\Gamma_k,X)}=A^{(\overline\Gamma_k,X)}+C^{(\overline
\Gamma_k,X)},
\end{equation}
where $X\in\{\overline\Omega^{(k)},\ k,\ \Gamma_k,\ \overline\Gamma_k,\ D_k\}$, so
that the algebraic implementation of (\ref{flux}) reads, for $k=1,2$,
\begin{equation}\label{res_modified}
{\bf
r}_{\overline\Gamma_k}=\widehat{A}^{(\overline\Gamma_k,\overline\Omega^{(k)})}{\bf u}^{(k)}
-{\bf f}_{\overline\Gamma_k}.
\end{equation}

By defining the two \emph{intergrid matrices}
\begin{equation}\label{matrices_Q}
Q_{21}=P_{21},\qquad Q_{12}=M_{\Gamma_1}P_{12} M_{\Gamma_2}^{-1},
\end{equation}

the algebraic counterpart of (\ref{internodes_weak})$_{2,3}$ read
(see (\ref{matrix_P21}) and (\ref{residual_equil_alg}))
\begin{eqnarray}\label{eq_trace_nc_alg}
\begin{array}{l}
 {\bf u}_{\overline\Gamma_2}=Q_{21}  {\bf u}_{\overline\Gamma_1},\qquad \qquad
{\bf r}_{\overline\Gamma_1}+Q_{12}{\bf r}_{\overline\Gamma_2}={\bf 0}.
\end{array}
\end{eqnarray}

\noindent
By introducing the following submatrices:
\begin{equation*}
Q_{21}^{(\overline\Gamma_2,\Gamma_1)}=Q_{21}({\cal I}_{\overline\Gamma_2},{\cal
I}_{\Gamma_1}),\quad
Q_{21}^{(\Gamma_2,\partial\Gamma_1)}=Q_{21}({\cal I}_{\Gamma_2},{\cal
I}_{\partial\Gamma_1}),\quad
Q_{12}^{(\Gamma_1,\overline\Gamma_2)}=Q_{12}({\cal I}_{\Gamma_1},{\cal I}_{\overline\Gamma_2}),
\end{equation*}
and by using (\ref{eq_trace_nc_alg}),
the algebraic form of (\ref{internodes_weak}) reads
\begin{eqnarray}\label{alg_2dom_nc_3x3}
\underbrace{\left[\begin{array}{ccc}
A^{(1,1)} &  0 & A^{(1,\Gamma_1)} \\[1mm]
0 & A^{(2,2)}&   A^{(2,\overline\Gamma_2)}Q_{21}^{(\overline\Gamma_2,\Gamma_1)} \\[1mm]
\widehat{A}^{(\Gamma_1,1)} 
&Q_{12}^{(\Gamma_1,\overline\Gamma_2)}\widehat{A}^{(\overline\Gamma_2,2)}
& \widehat{A}^{(\Gamma_1,\Gamma_1)}+
 Q_{12}^{(\Gamma_1,\overline\Gamma_2)}\widehat{A}^{(\overline\Gamma_2,
\overline\Gamma_2)}Q_{21}^{(\overline\Gamma_2,\Gamma_1)} 
\end{array}
\right]}_{\displaystyle\mathbb A}
\left[
\!\!
\begin{array}{l}
{\bf u}^{(1)}_0\\[1mm]
{\bf u}^{(2)}_0\\[1mm]
{\bf u}_{\Gamma_1}
\end{array}
\!\!
\right] 
=
\left[
\!\!
\begin{array}{c}
{\bf f}^{(1)}_0\\[1mm]
{\bf f}^{(2)}_0\\[1mm]
{\bf f}_{\Gamma_1}+Q_{12}^{(\Gamma_1,\overline\Gamma_2)}
{\bf f}_{\overline\Gamma_2}
\end{array} 
\!\!
\right]-{\bf G},
\end{eqnarray}
where  the array
\begin{eqnarray}\label{G_dirichlet}
{\bf G}=
\left[
\begin{array}{c}
{\bf G}_1\\
{\bf G}_2\\
{\bf G}_{\Gamma_1}
\end{array}
\right]=
\left[
\begin{array}{c}
A^{(1,D_1)}{\bf g}_1\\[1mm]
A^{(2,D_2)}{\bf g}_2+A^{(2,\Gamma_2)}Q_{21}^{(\Gamma_2,
\partial\Gamma_1)}{\bf g}_{\partial\Gamma_1}\\[1mm]
\widehat{A}^{(\Gamma_1,D_1)}{\bf g}_1+
 Q_{12}^{(\Gamma_1,\overline\Gamma_2)}(\widehat{A}^{(\overline\Gamma_2,D_2)}
{\bf g}_2+\widehat{A}^{(\overline\Gamma_2,\Gamma_2)}Q_{21}^{(\Gamma_2,
\partial\Gamma_1)}{\bf g}_{\partial\Gamma_1})
\end{array}
\right]
\end{eqnarray}
is non null only when non-homogeneous Dirichlet conditions are 
given on $\partial\Omega_D$ and implements the lifting of the Dirichlet datum.

Finally the degrees of freedom in $\Omega^{(1)}$ are given by
${\bf u}^{(1)}=[{\bf u}^{(1)}_0,\ {\bf u}_{\Gamma_1},\ {\bf g}_1]$
while the one in  $\Omega^{(2)}$ are given by
${\bf u}^{(2)}=[{\bf u}^{(2)}_0,\ Q_{21}^{(\Gamma_2,\overline\Gamma_1)}
{\bf u}_{\overline\Gamma_1}, {\bf g}_2]$.

The presence of the terms
$\widehat{A}^{(\overline\Gamma_2,\Gamma_2)}Q_{21}^{(\Gamma_2,
\partial\Gamma_1)}{\bf g}_{\partial\Gamma_1}$
and $A^{(2,\Gamma_2)}Q_{21}^{(\Gamma_2,
\partial\Gamma_1)}{\bf g}_{\partial\Gamma_1}$  in
the last two rows of (\ref{G_dirichlet}) is motivated by the fact that
the trace of $u^{(2)}_{h_2}$ on
the interface $\Gamma_2$ is the interpolation through $\Pi_{21}$
of the trace of $u^{(1)}_{h_1}$ on $\Gamma_1$.

System (\ref{alg_2dom_nc_3x3}) represents the algebraic form
of INTERNODES implemented
in practice.
By taking $Q_{12}=Q_{21}=I$ we recover the algebraic system
associated with classical
conforming domain decomposition (see, e.g., \cite{qv_ddm,tw}).

Notice that, even though the residuals are defined up to the boundary
of
$\Gamma_k$, the algebraic counterpart of condition
(\ref{internodes_weak})$_3$ is imposed only on the degrees of freedom 
internal to $\Gamma_1$.
In this way the number of equations and the number of unknowns in
(\ref{alg_2dom_nc_3x3}) do coincide.

\subsection{An iterative method to solve 
 (\ref{alg_2dom_nc_3x3})}\label{sec:schur}

An alternative way to solve system (\ref{alg_2dom_nc_3x3}) consists in
eliminating the variables ${\bf u}^{(1)}_0$ and
${\bf u}^{(2)}_0$ and in solving the Schur complement system (\cite{tw,qv_ddm})
\begin{equation}\label{schur}
S{\bf u}_{\Gamma_1}={\bf b},
\end{equation}
where
\begin{eqnarray}
\label{schur_matrix}
S=S_{\Gamma_1}+Q_{12}^{(\Gamma_1,\overline\Gamma_2)}S_{\overline\Gamma_2}
Q_{21}^{(\overline\Gamma_2,\Gamma_1)}, &&
{\bf b}={\bf b}_{\Gamma_1}+Q_{12}^{(\Gamma_1,\overline\Gamma_2)}
{\bf b}_{\overline\Gamma_2}-{\bf G}_{\Gamma_1},\\
\label{local_schur}
S_{\Gamma_1}=\widehat{A}^{(\Gamma_1,\Gamma_1)}-
\widehat{A}^{(\Gamma_1,1)}(A^{(1,1)})^{-1}A^{(1,\Gamma_1)},
&&
S_{\overline\Gamma_2}=\widehat{A}^{(\overline\Gamma_2,\overline\Gamma_2)}-
\widehat{A}^{(\overline\Gamma_2,2)}(A^{(2,2)})^{-1}A^{(2,\overline\Gamma_2)}\\
\label{local_rhs}
{\bf b}_{\Gamma_1}={\bf f}_{\Gamma_1}-\widehat{A}^{(\Gamma_1,1)}
(A^{(1,1)})^{-1} ({\bf f}^{(1)}_0-{\bf G}_1), 
&&
{\bf b}_{\overline\Gamma_2}={\bf f}_{\overline\Gamma_2}
- \widehat{A}^{(\overline\Gamma_2,2)}(A^{(2,2)})^{-1}({\bf f}^{(2)}_0-{\bf G}_2).
\end{eqnarray}

$S_{\Gamma_1}$ and $S_{\overline\Gamma_2}$ are the local Schur complement
matrices, while ${\bf b}_{\Gamma_1}$ and ${\bf b}_{\overline\Gamma_2}$
are the local right hand sides.

System (\ref{schur}) can be solved, e.g., by a preconditioned Krylov method
(Bi-CGStab or GMRES) with $S_1$ as preconditioner.
Notice that the matrix $Q_{12}^{(\Gamma_1,\overline\Gamma_2)}
S_{\overline\Gamma_2}Q_{21}^{(\overline\Gamma_2,\Gamma_1)}$
 is not a good candidate to play the role of
preconditioner since it may be singular.

Since $Q_{12}$ is not the
transpose of $Q_{21}$, even if the differential operator is symmetric, the
Schur complement system $S$ is not. 
Nevertheless the local systems continue to be symmetric and can be solved
either by standard Cholesky factorization or the PCG method.
For example, we can compute and store
the Cholesky factorization of the
matrices $A^{(k,k)}$ and dispose of a
function that implements
the action of $S$ on a given array ${\boldsymbol\lambda}$
whose entries are the degrees of freedom associated with $\Gamma_1$.

 We will describe the iterative approach in 
Sect. \ref{sec:schurM}  for general 
multipatch configurations.

Once ${\bf u}_{\Gamma_1}$ is known,
the variables ${\bf u}^{(1)}_0$ and ${\bf u}^{(2)}_0$ are recovered by
solving the local subsystems
\begin{eqnarray*}
A^{(1,1)}{\bf u}^{(1)}_0&=&{\bf f}^{(1)}_0-{\bf G}_1-A^{(1,\Gamma_1)}{\bf
u}_{\Gamma_1},\\
A^{(2,2)}{\bf u}^{(2)}_0&=&{\bf f}^{(2)}_0-{\bf G}_2-A^{(2,\overline\Gamma_2)}
Q_{21}^{(\overline\Gamma_2,\Gamma_1)}{\bf u}_{\Gamma_1}.
\end{eqnarray*}

Finally, ${\bf u}_{\overline\Gamma_1}$ is recovered by assembling ${\bf
u}_{\Gamma_1}$ and ${\bf g}_{\partial\Gamma_1}$ and 
the numerical solution on $\overline\Gamma_2$ is reconstructed by the
interpolation formula ${\bf u}_{\overline \Gamma_2}=Q_{21}{\bf
u}_{\overline\Gamma_1}$.

\section{Numerical results for 2 patches}\label{sec:numres_2dom}

The aim of this section is twofold. From one hand we show that
INTERNODES does not deteriorate the accuracy of IGA discretization (we say that
the method exhibits optimal accuracy).
Then we show that, if the
interfaces  are non-watertight, the accuracy of the
INTERNODES solution depends on the maximum size $d_\Gamma$
of gaps and overlaps between 
$\Gamma_1$ and $\Gamma_2$, and smaller $d_\Gamma$, smaller the error.

To highlight these features of INTERNODES we consider here very regular
solutions.
Other test case with less regular solutions will be taken into account in Sect.
\ref{sec:numres_mdom} in the case of $M>2$ patches.

 \subsection{Test case \#1.}

Let us consider the differential problem (\ref{problem}) in $\Omega=\{(x,y)\in
{\mathbb R}^2:\ x\geq0, \ y\geq 0, \ 1\leq x^2+y^2\leq 4\}$
 with  $\alpha=0$, and 
$f$ and  $g$ such that the exact solution is
$u(x,y)=\sin(1.5\pi x)\sin(3\pi y)$.

\begin{figure}
\begin{center}
\begin{subfigure}[t]{2in}
\includegraphics[width=0.7\textwidth]{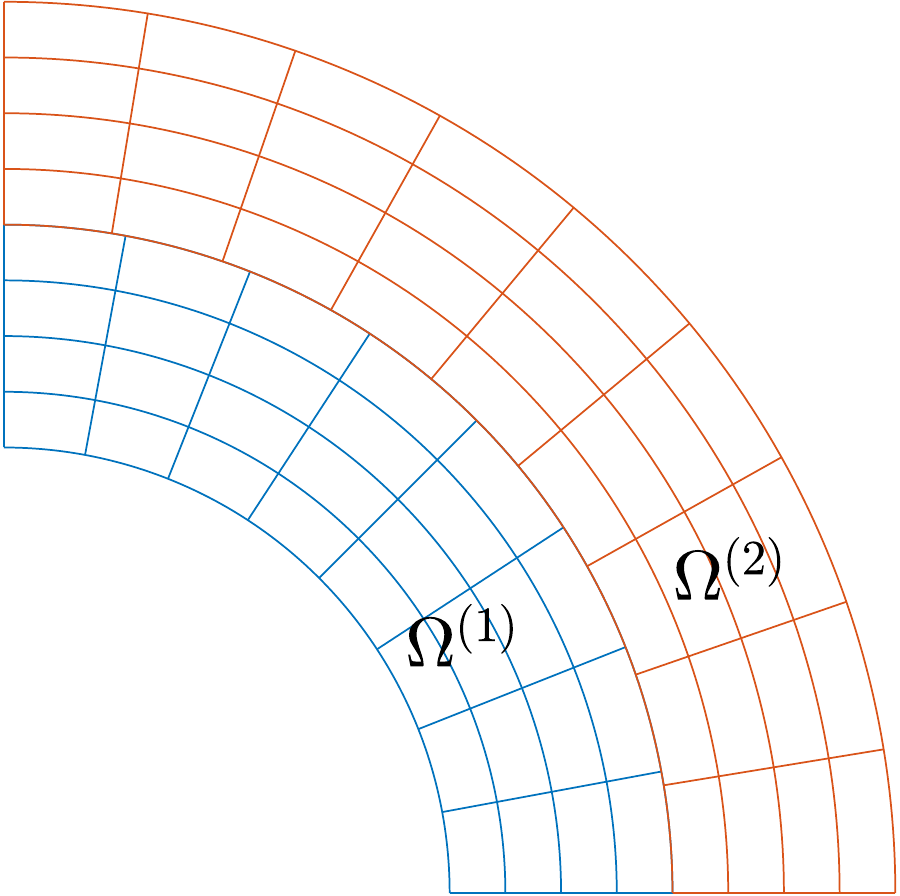}
\subcaption{}
\end{subfigure}
\begin{subfigure}[t]{2in}
\includegraphics[width=0.7\textwidth]{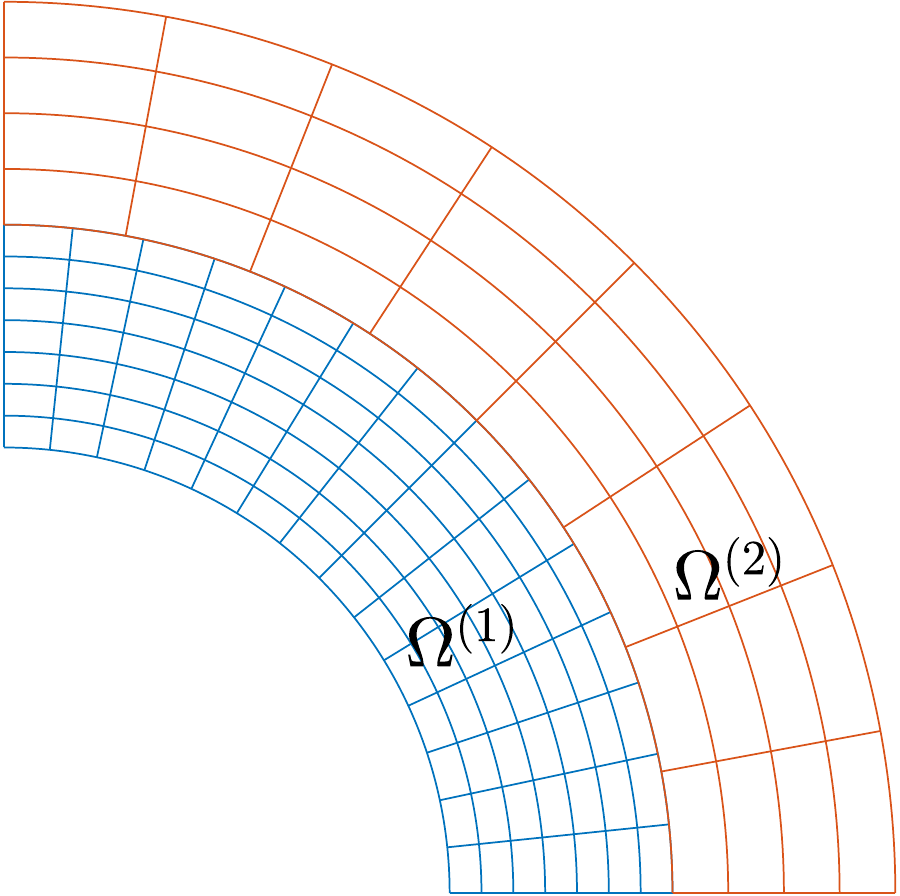}
\subcaption{}
\end{subfigure}
\begin{subfigure}[t]{2in}
\includegraphics[width=0.7\textwidth]{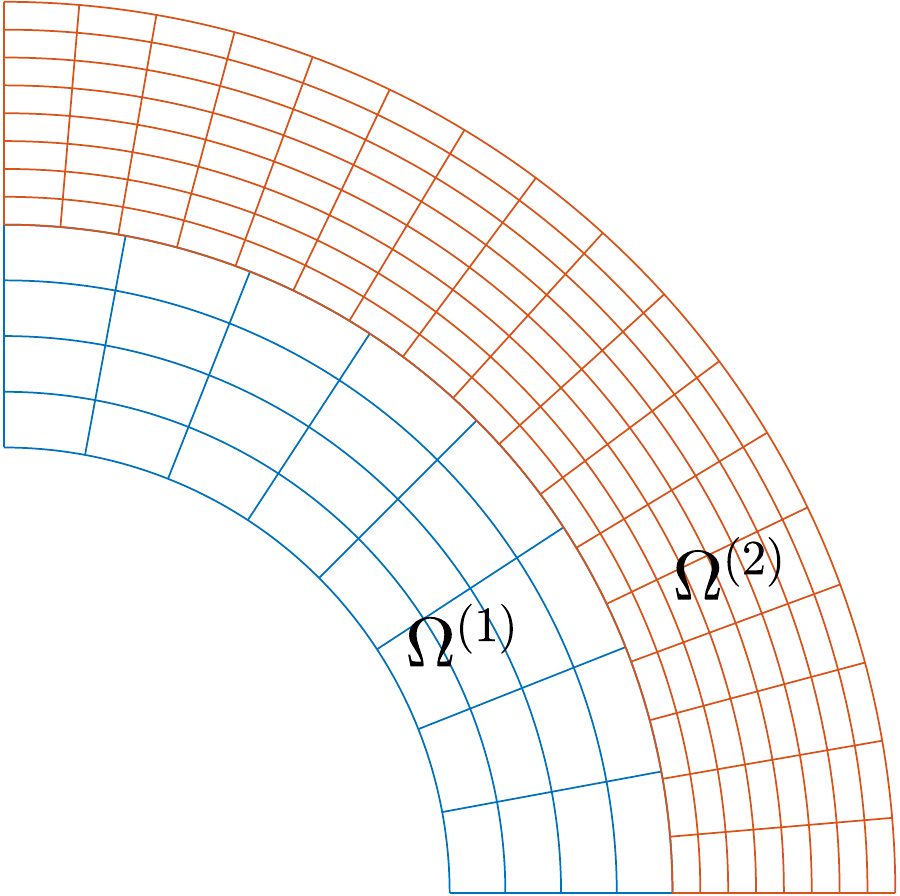}
\subcaption{}
\end{subfigure}
\end{center}
\caption{\emph{Test case \#1.} 
Balanced  (a), 
master-refined  (b), slave-refined (c) discretization sets with $\overline
n=8$}
\label{fig:ring_2dom}
\end{figure}

The computational domain $\Omega$ is split into the patches
$\Omega^{(1)}=\{(x,y)\in\Omega:\ x^2+y^2\leq (1.5)^2\}$ and 
$\Omega^{(2)}=\{(x,y)\in\Omega:\ (1.5)^2\leq x^2+y^2\leq 4\}$ (see Fig.
\ref{fig:ring_2dom}). Each patch is parameterized by NURBS as 
described in Sect. \ref{sec:discretization}.
The weights and
the control points of the circular arcs are chosen as described in \cite[Sect.
2.4.1.1]{chb_iga_book}; each patch is built first as a 
single element with 6 control points, more precisely:
\begin{center}
\begin{tabular}{l|cccccc}
${\bf P}_i^{(1)}$ & (1,0) & (1.5,0) & (0,1) & (1,1) & (0,1.5) &
(1.5,1.5)\\[2mm]
${\bf P}_i^{(2)}$ & (1.5,0) & (2,0) & (0,1.5) & (1.5,1.5) & (0,2) &
(2,2)\\[2mm]
$w_i^{(1)}=w_i^{(2)}$&  1 & 1 & 1 & $\sqrt{2}/2$  &  1 & $\sqrt{2}/2$,
\end{tabular}
\end{center}

then it is ${\sf k}-$refined (see 
\cite[Sect. 2.1.4.3]{chb_iga_book}) up to polynomial degree $p^{(k)}$  and 
continuity order $p^{(k)}-1$ along both the coordinates,
and finally it is uniformly $h-$refined leaving the underlying geometry and its
parametrization intact (see \cite[Sect. 2.1.4.1]{chb_iga_book}).

Then, let
\begin{eqnarray}\label{numerical_sol}
u_h=\left\{\begin{array}{ll}
u_{h_1}^{(1)} & \mbox{ in }\Omega^{(1)}\\
u_{h_2}^{(2)} & \mbox{ in }\Omega^{(2)}
\end{array}\right.
\end{eqnarray}
denote the numerical solution computed with INTERNODES.

We consider three non-matching parametrizations named:
  balanced, master-refined and slave-refined (in fact
these are $h-$refinements).
We take equal polynomial degrees $p^{(1)}=p^{(2)}=p\in\{2,\ldots, 5\}$ 
and, for $\overline n\in\{8,16,24,32\}$,
 we define the number of elements inside the patches as follows:
\begin{center}
\begin{tabular}{l|ccc}
patch & balanced & master-refined & slave-refined\\
\hline
$\Omega^{(1)}$ & $(\overline n/2)\times \overline n$ &
             $(\overline n-1)\times 2(\overline n-1)$ &
             $(\overline n/2)\times \overline n$\\
$\Omega^{(2)}$ & $\overline n/2 \times (\overline n+1)$ &
             $(\overline n/2)\times \overline n$ &
             $\overline n\times (2\overline n+1)$.
\end{tabular}
\end{center}
The first (second, resp.)
parameter coordinate is mapped onto the physical radial (angular, resp.) 
coordinate, and the non-conformity is a consequence of the different number of
elements along the second coordinate.

\begin{figure}
\begin{center}
\includegraphics[width=0.3\textwidth]{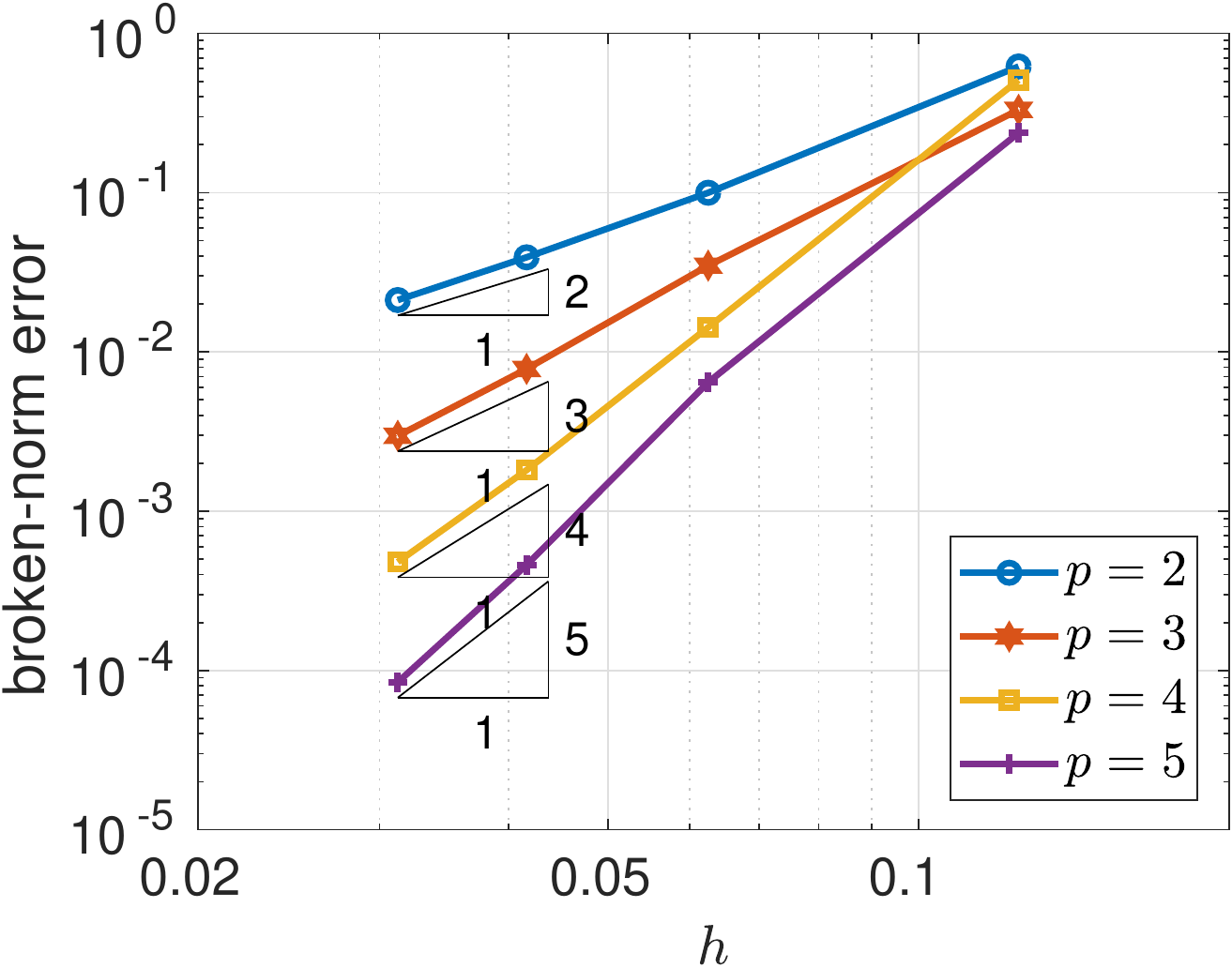}\quad
\includegraphics[width=0.3\textwidth]{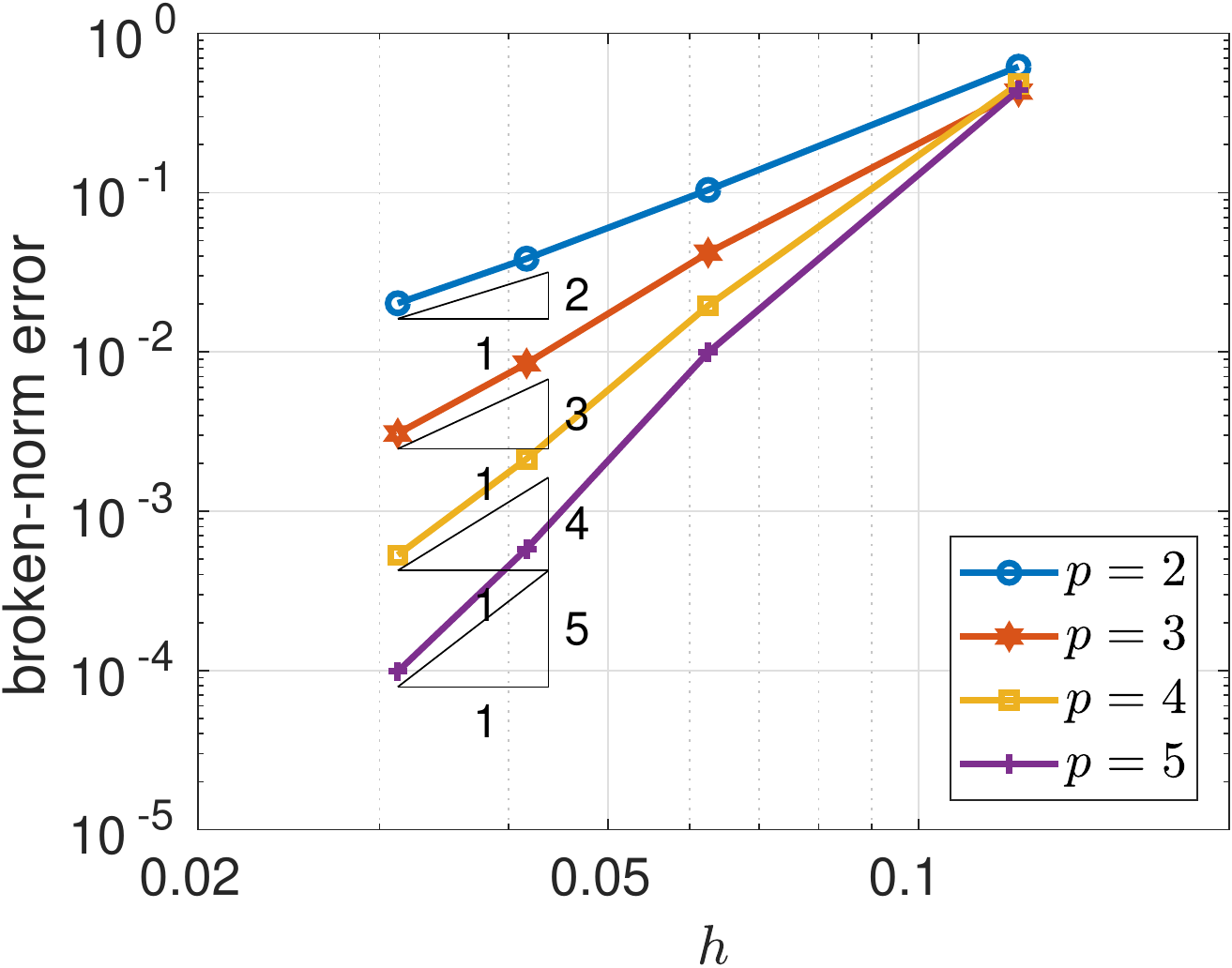}\quad
\includegraphics[width=0.3\textwidth]{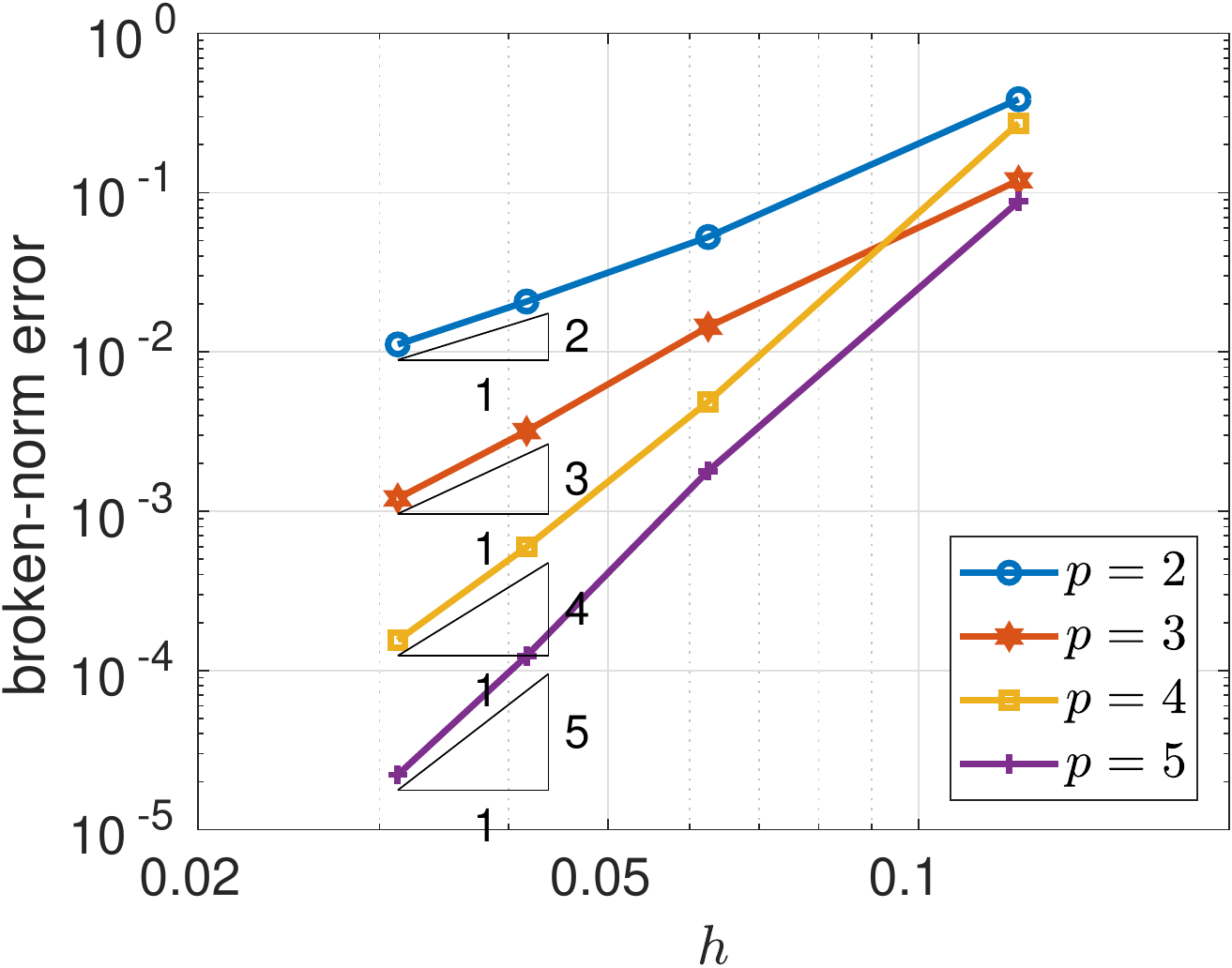}
\end{center}
\caption{\emph{Test case \#1.}
The broken-norm error (\ref{err_broken}) versus the mesh size.
At left the error for the balanced configuration, in the middle the error for
the master-refined configuration, at right the error for the slave-refined
configuration}
\label{fig:errori_2dom_h}
\end{figure}

In Fig. \ref{fig:ring_2dom},
 we show the the three discretization sets when $\overline n=8$,
while in Fig. \ref{fig:errori_2dom_h} we show the broken-norm errors 
\begin{equation}\label{err_broken}
\|u_h-u\|_*=\left(\sum_{k}
\frac{\|u_h-u\|^2_{H^1(\Omega^{(k)})}}{\|u\|^2_{H^1(\Omega^{(k)})}}
\right)^{1/2}
\end{equation}
with respect to the maximum mesh size $h=\max_k h_k$
($h_k$ is the mesh-size in $\Omega^{(k)}$),
for the three discretization sets. 

The convergence of INTERNODES is optimal versus the mesh-size $h$,
 in the sense that the broken-norm
errors behave like $h^{p}$ when $h\to0$, 
exactly as the error in $H^1-$norm of the 
Galerkin Isogeometric methods (see, e.g., \cite[Thm. 3.4 and Cor. 4.16]{bbsv}).

\begin{figure}
\begin{center}
\includegraphics[width=0.35\textwidth]{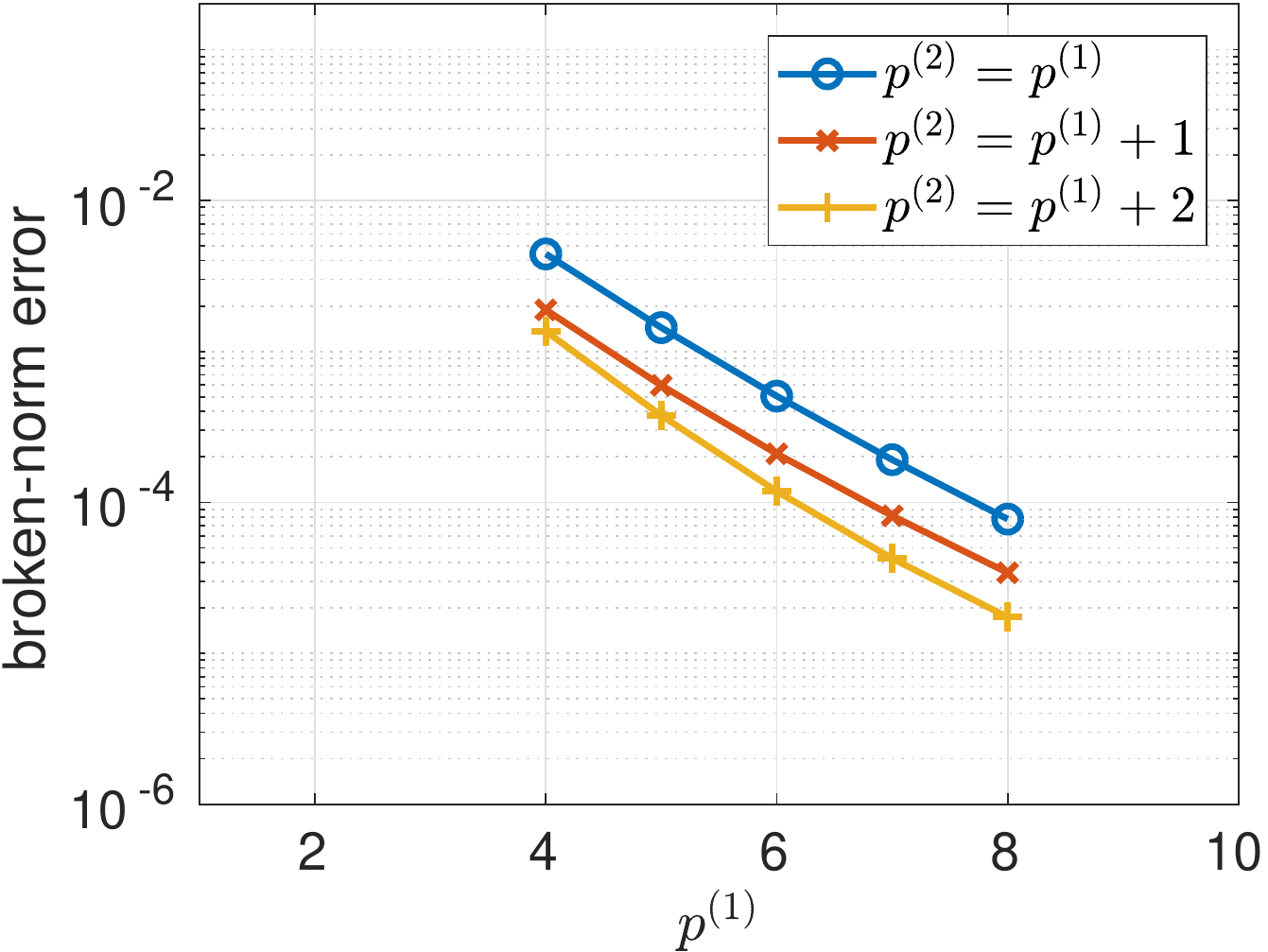}\quad
\includegraphics[width=0.34\textwidth]{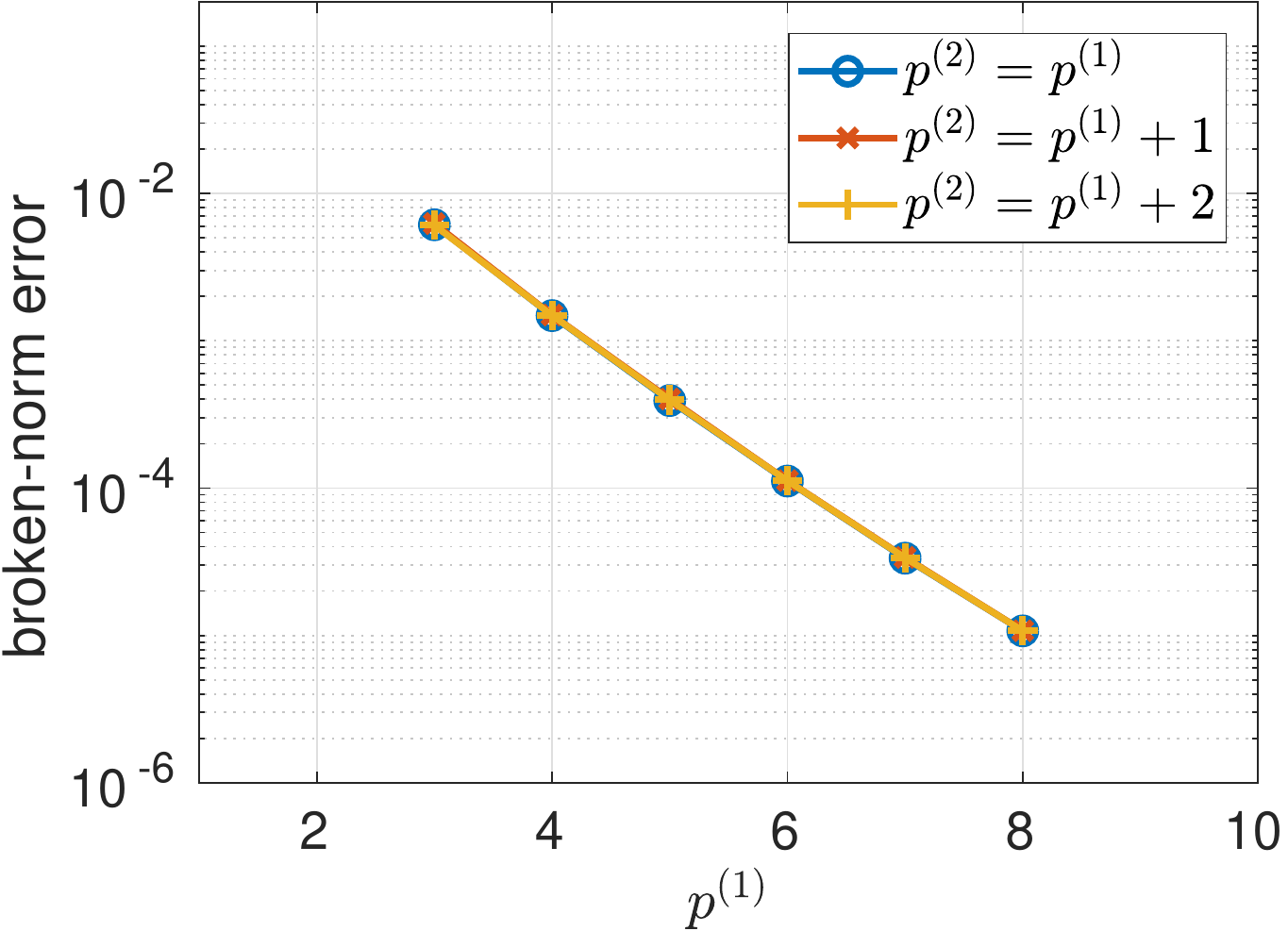}
\end{center}
\caption{\emph{Test case \#1.}
The broken-norm error versus the polynomial degree $p^{(1)}$ in the patch
$\Omega^{(1)}$.
At left the error for the balanced configuration 
 at right the error for
the slave refined configuration, $\overline n=20$ 
for both cases}
\label{fig:errori_2dom_p}
\end{figure}

\begin{figure}[t!]
\begin{center}
\includegraphics[width=0.4\textwidth]{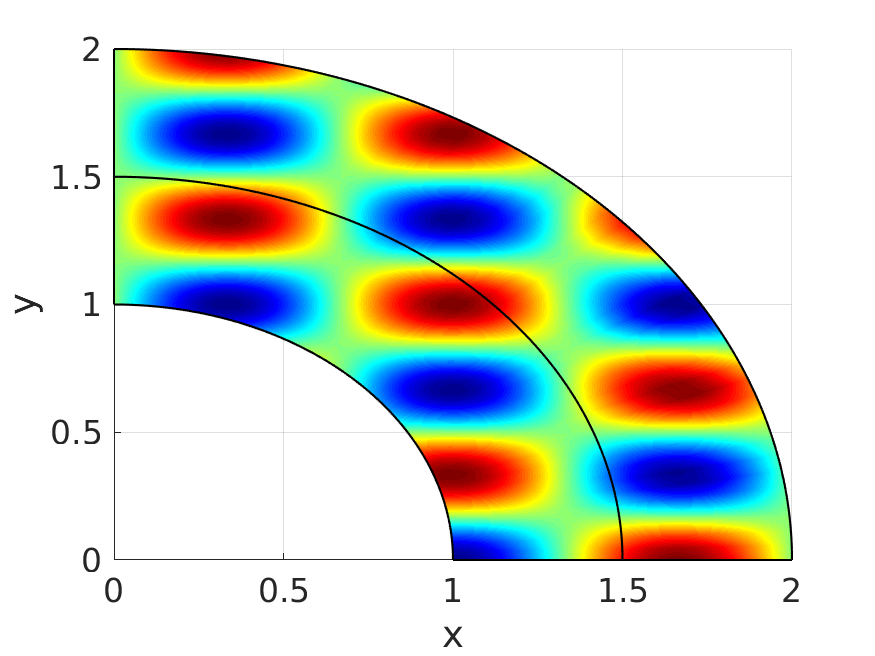}\quad
\includegraphics[width=0.5\textwidth]{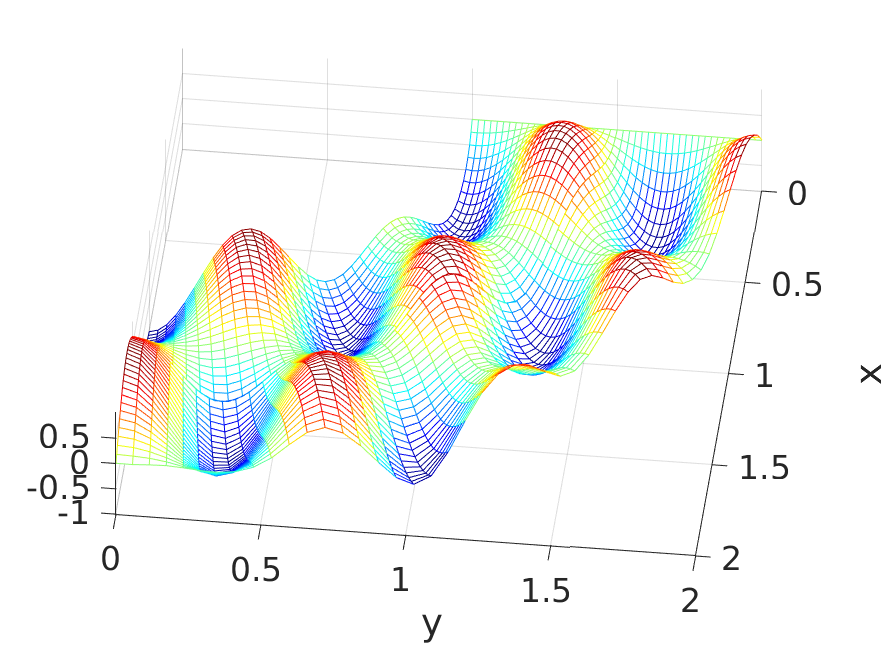}
\end{center}
\caption{\emph{Test case \#1.}
The numerical solution obtained with
$p^{(1)}=p^{(2)}=3$, $(8\times 16)$ elements in $\Omega^{(1)}$, and 
$(8\times 17)$ elements in $\Omega^{(2)}$. The first (second, resp.)
parameter coordinate is mapped
on the physical radial (angular, resp.) coordinate. 
The number of grid points used in the right plot is about
twice the number of degrees of freedom}
\label{fig:sol_2dom}
\end{figure}

In Fig. \ref{fig:errori_2dom_p} we show the broken-norm errors versus the
polynomial degree $p^{(1)}$, with three different choices for $p^{(2)}$:
$p^{(2)}=p^{(1)}$, $p^{(2)}=p^{(1)}+1$ and $p^{(2)}=p^{(1)}+2$ and the
discretization sets: balanced and slave refined. In all the cases 
$\overline n=20$. Finally, in Fig. \ref{fig:sol_2dom},
we present two qualitative pictures of the INTERNODES solution.

\subsection{Test case \#2. Non-watertight patches}\label{sec:test2}
Let us consider the differential problem (\ref{problem}) in
$\Omega=(0,2)\times(0,1)$
 with  $\alpha=0$, and
$f$ and  $g$ such that the exact solution is
$u(x,y)=e^{-3(x-1)^2-4(y-0.6)^2}(1+\sin(3\pi x)\cos(3\pi y)).$

\begin{figure}
\begin{center}
\includegraphics[width=0.4\textwidth]{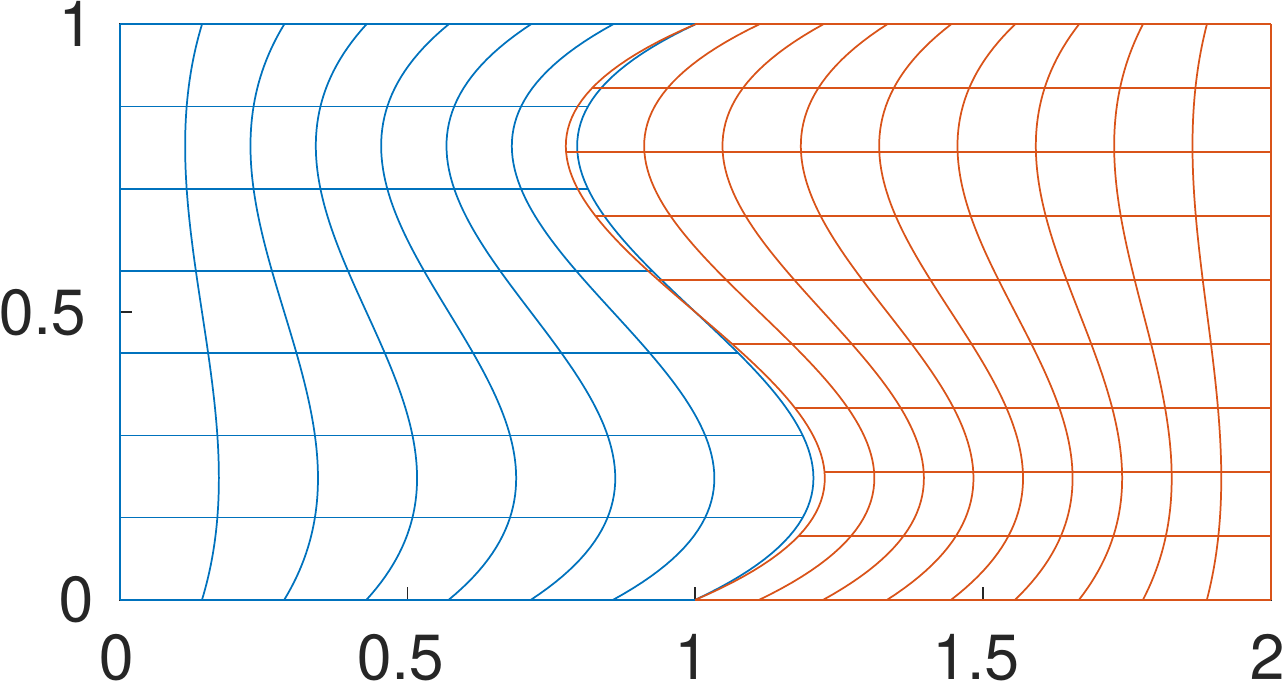}\quad
\includegraphics[trim={0 60 0 0 },width=0.5\textwidth]{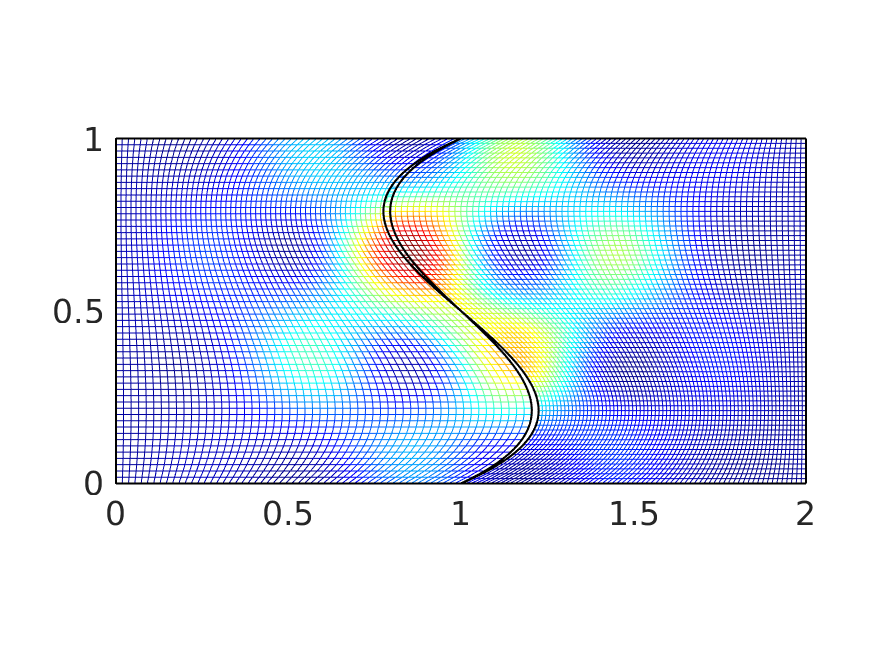}
\end{center}
\caption{\emph{Test case \#2.} At left, the parameterization of the two patches
with $7\times 7$ elements in $\Omega^{(1)}$ and  $9\times 9$ elements in
$\Omega^{(2)}$. The maximum distance between $\Gamma_1$ and $\Gamma_2$ is 
 $d_\Gamma=0.0197$.
At right the numerical solution computed by INTERNODES (with
RL-RBF interpolation) with $p^{(1)}=4$ and $p^{(2)}=3$.
The number of grid points used in the plot is about
four times the number of degrees of freedom }
\label{fig:test2_patches}
\end{figure}

The computational domain $\Omega$ is split into two patches as shown in Fig.
\ref{fig:test2_patches}, left, where we have considered
a sinusoidal physical interface $\Gamma_{12}=\{(x,y)\in {\mathbb R}^2: \
x=g(y)=1+0.2\sin(2\pi y), \ 0\leq y\leq 1\}$.
The interfaces $\Gamma_1$ and $\Gamma_2$ are built as
B-spline interpolation of $\Gamma_{12}$ and they are non-watertight.
The size of gaps and overlaps that are generated by the approximation 
(we denote by $d_\Gamma$ the maximum distance between the two interfaces)
depends on the parameterization of the patches.
To face the non-watertight interfaces we implement INTERNODES
with the RL-RBF interpolation
operators defined in Sect. \ref{sec:rbf_interpolation}.

We analyze the accuracy of INTERNODES by measuring the 
broken-norm error (\ref{err_broken}) versus the mesh size $h=\max_k h_k$ in two
situations: \emph{i)} with fixed $d_\Gamma$, \emph{ii)} with variable $d_\Gamma$.

\emph{i) Fixed $d_\Gamma$.} We
fix two different polynomial degrees $p^{(1)}$ and 
$p^{(2)}$, thus for $k=1,2$ we determine the interfaces $\Gamma_k$ as follows:
\begin{enumerate}[noitemsep]
\item let $\Xi^{(k)}=\{\underbrace{0\ldots 0}_{p^{(k)}},
\underbrace{1\ldots 1}_{p^{(k)}}\}$ be the knots set on a single element,
\item let $\widehat\varphi_{i,p^{(k)}}^{(k)}$  denote the univariate B-spline basis functions of
degree $p^{(k)}$ associated with $\Xi^{(k)}$,
\item let $\widehat t_{j,G}^{(k)}$, for $j=1,\ldots,p^{(k)}+1$ the Greville abscissae 
associated with $\Xi^{(k)}$,
\item let ${\bf Q}_j^{(k)}=(g(\widehat t_{j,G}^{(k)}),\widehat t_{j,G}^{(k)})$, for $j=1,\ldots,
p^{(k)}+1$, 
\item compute the points ${\bf P}_i^{(k)}$, for $i=1,\ldots p^{(k)}+1$,
by solving the linear system
(see, e.g. \cite[Sect.  9.2.1]{piegl_nurbs_book})
\begin{equation}
\sum_{i=1}^{p^{(k)}+1}\widehat\varphi_{i,p^{(k)}}^{(k)}(\widehat t_{j,G}^{(k)}){\bf P}_i^{(k)}=
{\bf Q}_j^{(k)}, \qquad j=1,\ldots, p^{(k)}+1.
\end{equation}
\end{enumerate}

The points ${\bf P}_i^{(k)}$ are in fact the   control points of the curve
that defines $\Gamma_k$. We have chosen the associated weights equal to 1
(but different weights can be chosen as well providing NURBS 
instead of B-spline parametrization), 
the parameterization of 
the patch $\Omega^{(k)}$ can be defined as ruled surface between
$\Gamma_k$ and the opposite straight side. Finally, uniform
$h-$refinement is
adopted to reach the desired final discretization, while the polynomial degrees
remain fixed.

\emph{ii) Variable $d_\Gamma$}. In this case we choose equal polynomial degrees
$p^{(1)}=p^{(2)}$. Moreover, the initial set $\Xi^{(k)}$ includes all the knots
corresponding to the desired uniform $h-$refinement.
Then we proceed as in case \emph{i)}, by omitting the final $h-$refinement.

\begin{figure}[b!]
\begin{center}
\includegraphics[width=0.37\textwidth]{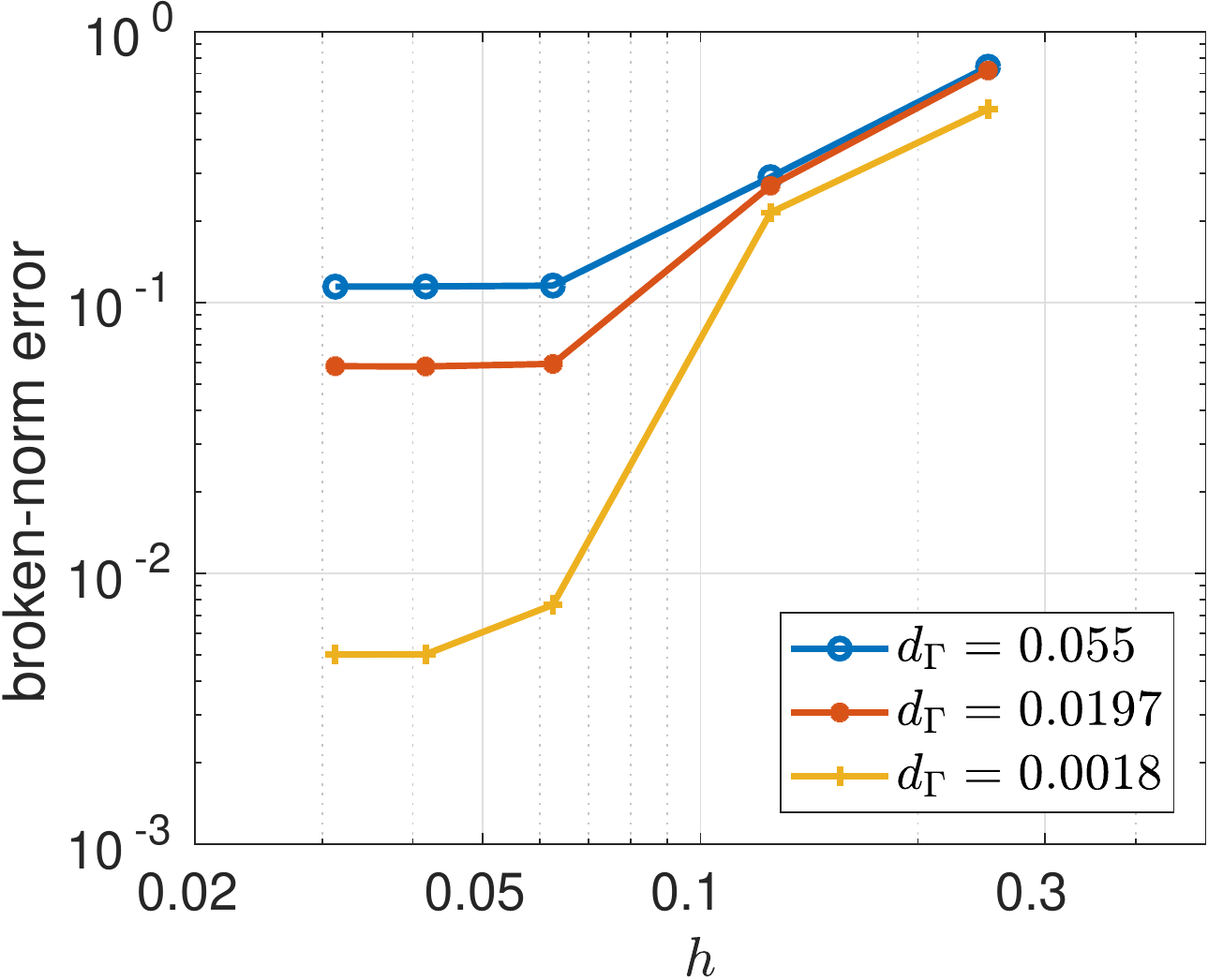}\hskip 1.cm
\includegraphics[width=0.4\textwidth]{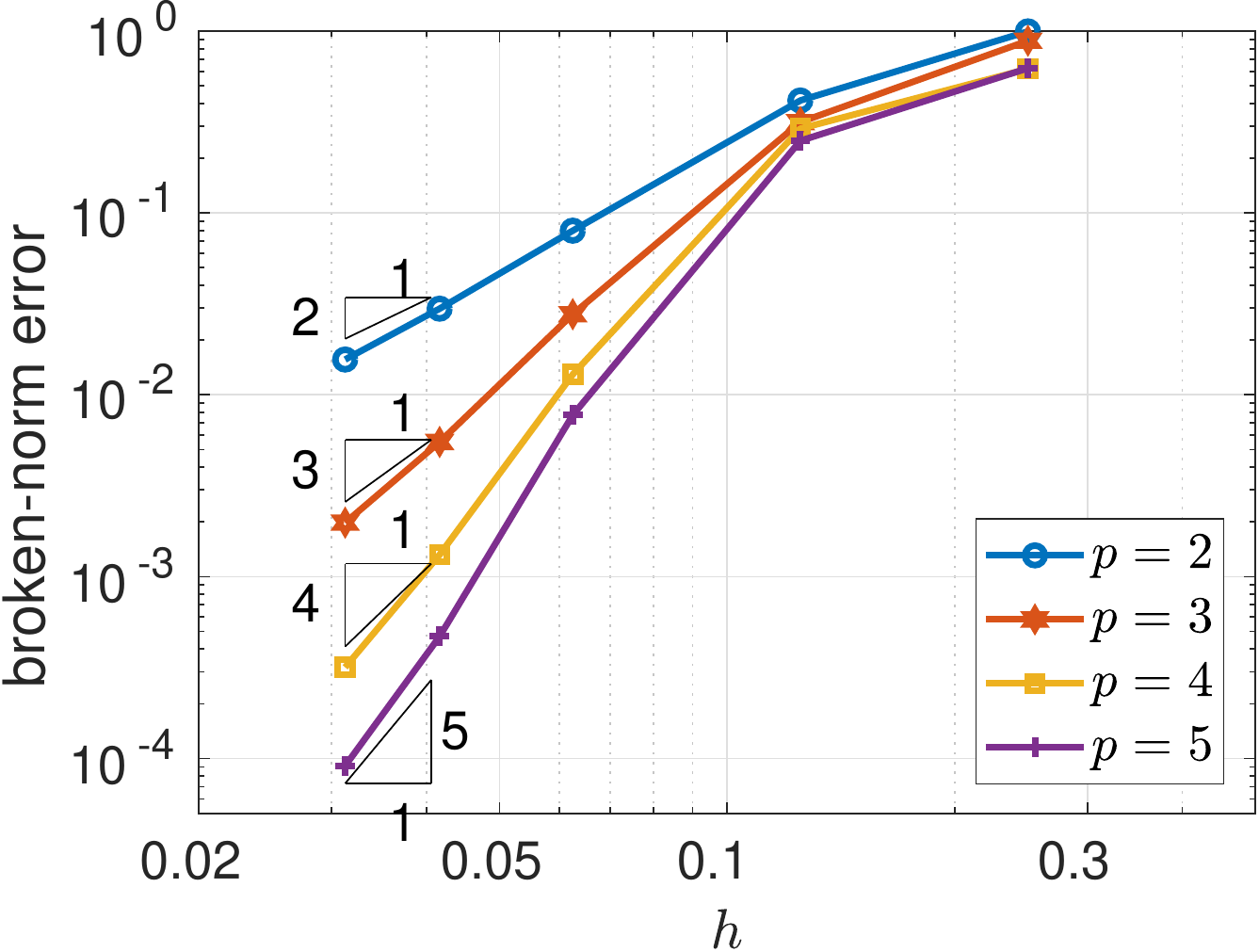}
\end{center}
\caption{\emph{Test case \#2.} At left, the errors for the case \emph{i)}:
$p^{(1)}=5$ and $p^{(2)}=3$ when $d_\Gamma=0.055$; 
$p^{(1)}=4$ and $p^{(2)}=3$ when $d_\Gamma=0.0197$; 
$p^{(1)}=6$ and $p^{(2)}=5$ when $d_\Gamma=0.0018$.
At right, the errors for the
case \emph{ii)} with $p=p^{(1)}=p^{(2)}$ and $d_\Gamma={\cal O}(h^p)$}
\label{fig:test2_err}
\end{figure}

For the case \emph{i)} we have chosen three couples of polynomial degrees:
 the first one with 
$p^{(1)}=5$ and $p^{(2)}=3$ that provides $d_\Gamma=0.055$; the second one
with $p^{(1)}=4$ and $p^{(2)}=3$ that provides $d_\Gamma=0.0197$; the third one
with $p^{(1)}=6$ and $p^{(2)}=5$ that provides $d_\Gamma=0.0018$.
In the left picture of Fig. \ref{fig:test2_err} we plot the broken-norm errors
versus the mesh size $h=\max_k h_k$; we have taken
$(\overline n-1)\times (\overline n-1)$ elements in $\Omega^{(1)}$ and
$(\overline n+1)\times (\overline n+1)$ elements in $\Omega^{(1)}$.
The plateaux in the error curves are clearly due to the presence of 
non-watertight interfaces and,
smaller $d_\Gamma$, lower the plateau.

For the case \emph{ii)} we have chosen $p=p^{(1)}=p^{(2)}\in\{2,\ldots,5\}$ and
different values of $h=\max_k h_k$. In this case the maximum size
$d_\Gamma$ of gaps and overlaps decreases as $h^p$ (in fact  it corresponds to
 the B-spline interpolation error of the curve $g(y)$).
In the right picture of Fig. \ref{fig:test2_err} 
the corresponding broken-norm errors are shown,
they decrease like $h^p$ when $h\to0$, as in the case of watertight
interfaces.

We have solved the Schur complement system (\ref{schur}) by the preconditioned 
Bi-CGStab (PBi-CGStab) method
as mentioned in Sect. \ref{sec:schur}, with  preconditioner $P$ given by
the local Schur complement matrix $S_{\Gamma_1}$.
It is well known (see, e.g., \cite{qv_ddm}) that,
in the framework of conforming Finite
Element/Spectral Element discretizations, 
such a preconditioner is optimal,
that is the condition number of
$(S_{\Gamma_1})^{-1}S$ is bounded independently of the mesh size and the
polynomial degree. Our numerical results show that this preconditioner looks
optimal also for non-matching NURBS parametrizations.

In Table \ref{tab:test2_iter} we report the number of iterations required by
the PBi-CGStab to solve
the Schur complement system (\ref{schur}) up to a
tolerance $\epsilon=10^{-10}$. 
In the right picture of 
Fig. \ref{fig:test2_patches} the INTERNODES solution obtained 
with RL-RBF interpolation is shown. We have set
$p^{(1)}=5$, $p^{(2)}=3$, $7\times 7$ elements in $\Omega^{(1)}$ and $9\times
9$ elements in $\Omega^{(2)}$.

\begin{table}
\begin{center}
\begin{tabular}[t]{l|ccc}
$\overline n$ & $d_\Gamma=0.055$ & $d_\Gamma=0.0197$ & $d_\Gamma=0.0018$\\
\hline 
$4 $ & 6 & 5 & 7 \\
$8 $ & 7 & 6 & 7 \\
$16$ & 8 & 8 & 7 \\
$24$ & 9 & 8 & 7 \\
$32$ & 9 & 9 & 7 \\
\end{tabular}
\quad
\begin{tabular}[t]{l|cccc}
$\overline n$ & $p=2$ & $p=3$ & $p=4$ & $p=5$ \\
\hline 
$4 $ & 3 & 4 & 5 & 6\\
$8 $ & 6 & 6 & 7 & 7\\
$16$ & 6 & 7 & 6 & 7 \\
$24$ & 7 & 7 & 7 & 7\\
$32$ & 7 & 7 & 7 & 7\\
\end{tabular}
\end{center}
\caption{\emph{Test case \#2.}
The number of Bi-CGStab iterations required to
solve the Schur complement system (\ref{schur}). The stop tolerance is
$\epsilon=10^{-10}$. We have set $(\overline n-1)\times (\overline
n-1)$ elements in $\Omega^{(1)}$ and $(\overline n+1)\times (\overline
n+1)$ elements in $\Omega^{(2)}$. The left table refers to case \emph{i)}, while
the right one to the case \emph{ii)}  }
\label{tab:test2_iter}
\end{table}

\section{More general second order elliptic PDEs}\label{sec:general_L}

INTERNODES methods can be applied to solve general 
elliptic second order PDEs, where the differential operator  is
\begin{equation}\label{def:L}
Lu=-\nabla\cdot(\nu\nabla u)+{\bf b}\cdot \nabla u+\alpha u,
\end{equation}
with
$\nu\in L^\infty(\Omega)$ such that there exists $\nu_0>0$ and 
$\nu\geq \nu_0$ a.e. in $\Omega$; 
${\bf b}=(b_1,\ldots,b_d)$, with $b_i\in L^\infty(\Omega)$;
$\alpha\in L^\infty(\Omega)$ with $\alpha\geq 0$ a.e. in $\Omega$.

Given $f\in L^2(\Omega)$ and $g\in H^{1/2}(\partial\Omega)$, and 
under the assumption that $\alpha-\frac{1}{2}\nabla\cdot{\bf b}\ge 0$ a.e. in 
$\Omega$, the problem to find $u\in H^1(\Omega)$ such that
\begin{eqnarray}\label{elliptic_problem}
\left\{\begin{array}{ll}
Lu=f & \mbox{ in }\Omega\\
u=g & \mbox{ on }\partial\Omega
\end{array}\right.
\end{eqnarray}
admits a unique solution that is stable w.r.t. the data $f$ and $g$.

In such a case, while the interface condition 
(\ref{transmission_strong})$_3$ enforcing the continuity of
the traces across $\Gamma_{12}$ does not change, the interface
conditions (\ref{transmission_strong})$_4$ 
involving the normal derivatives must be replaced by
\begin{equation}\label{ic_conormal}
 \nu_1\frac{\partial u^{(1)}}{\partial {\bf n}_1}+
 \nu_2\frac{\partial u^{(2)}}{\partial {\bf n}_2} =0 
\qquad \mbox{ on }\Gamma_{12},
\end{equation}
where  $\nu_k=\nu_{|\Omega^{(k)}}$.

\medskip

When Neumann boundary conditions are assigned on a subset
$\partial\Omega_N$ of the boundary $\partial\Omega$, the
definition of the set $G_k$ (that is used to define the real 
values $r_i^{(k)}$, see  (\ref{flux})) becomes
$
G_k=\partial\Omega^{(k)}\setminus(\Gamma_k\cup\partial\Omega_N),
$
(see \cite[formula (44)]{gq_camc}).


\section{INTERNODES for decompositions with $M\geq 2$
patches}\label{sec:internodesM}

Let now 
$\Omega^{(k)}$, with $k=1,\ldots,M$, denote a family of disjoint patches
of $\Omega\subset{\mathbb R}^d$, with $d=2,3$, s.t.
$\cup_k \overline{\Omega}^{(k)} = \overline{\Omega}$.
Let us suppose that each $\Omega^{(k)}$ has Lipschitz boundary
$\partial\Omega^{(k)}$ (for $k=1,\ldots,M$).

Let $\Gamma_k=\partial\Omega^{(k)}\setminus\partial\Omega$ be
the part of the boundary of $\Omega^{(k)}$ internal to $\Omega$, and
\begin{equation*}\label{Gammakl}
\Gamma_{k\ell}=\Gamma_{\ell k}=\partial\Omega^{(k)}\cap\partial\Omega^{(\ell)}
\end{equation*}
be the interface
 between the two subdomains $\Omega^{(k)}$ and $\Omega^{(\ell)}$.
Intersections having null measure in the topology of ${\mathbb R}^{d-1}$
are considered empty.
Finally, let $L$ be the differential operator introduced in (\ref{def:L}).

The multidomain formulation of problem (\ref{problem})
reads:  look for $u^{(k)}$ for
$k=1,\ldots,M$ such that:
\begin{eqnarray}\label{strong_multidomain}
\left\{\begin{array}{ll}
Lu^{(k)}=f & \mbox{ in }\Omega^{(k)}, \quad k=1,\ldots,M\\
u^{(k)}=g  & \mbox{ on }\partial\Omega_D^{(k)},\\
u^{(k)}=u^{(\ell)}, & \mbox{ on }\Gamma_{k\ell}\neq\emptyset,
\ \ell=1,\ldots,M, \ \ell\neq k\\
\displaystyle 
\nu_k\frac{\partial u^{(k)}}{\partial{\bf n}_k}+
\nu_\ell\frac{\partial u^{(\ell)}}{\partial{\bf n}_\ell}=0
& \mbox{ on }\Gamma_{k\ell}\neq\emptyset, \ \ell=1,\ldots,M, \ \ell\neq k.
\end{array}\right.
\end{eqnarray}

\begin{figure}
\begin{center}
\includegraphics[trim={0 20 0
20},width=0.25\textwidth]{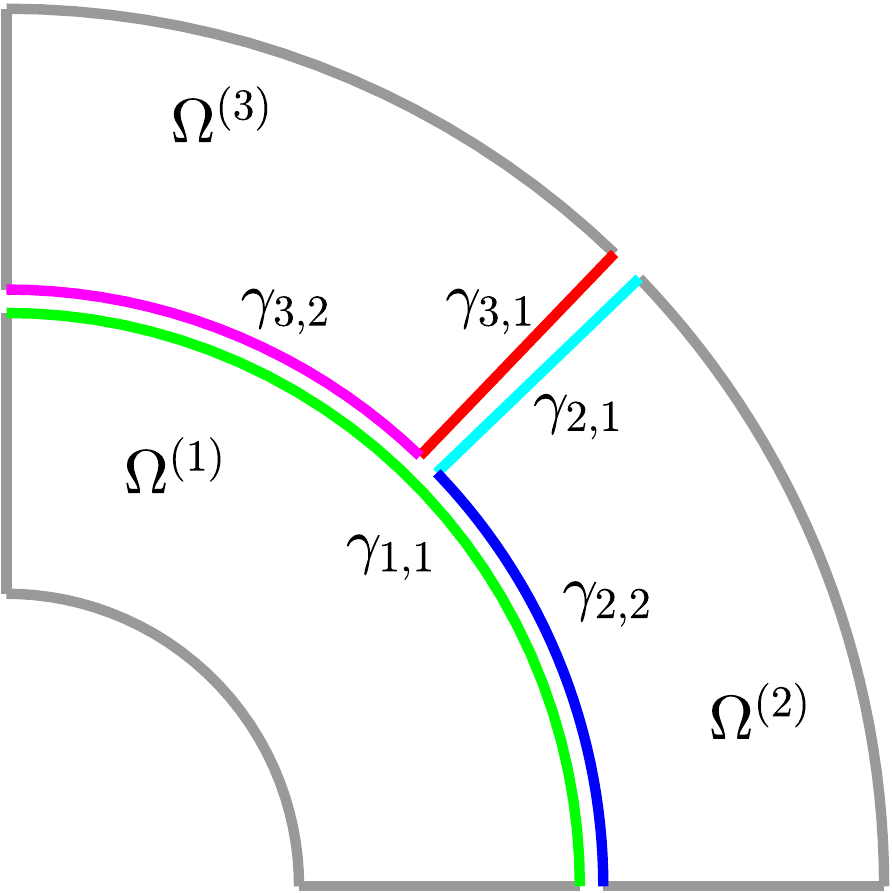}\quad
\includegraphics[trim={0 20 0
20},width=0.25\textwidth]{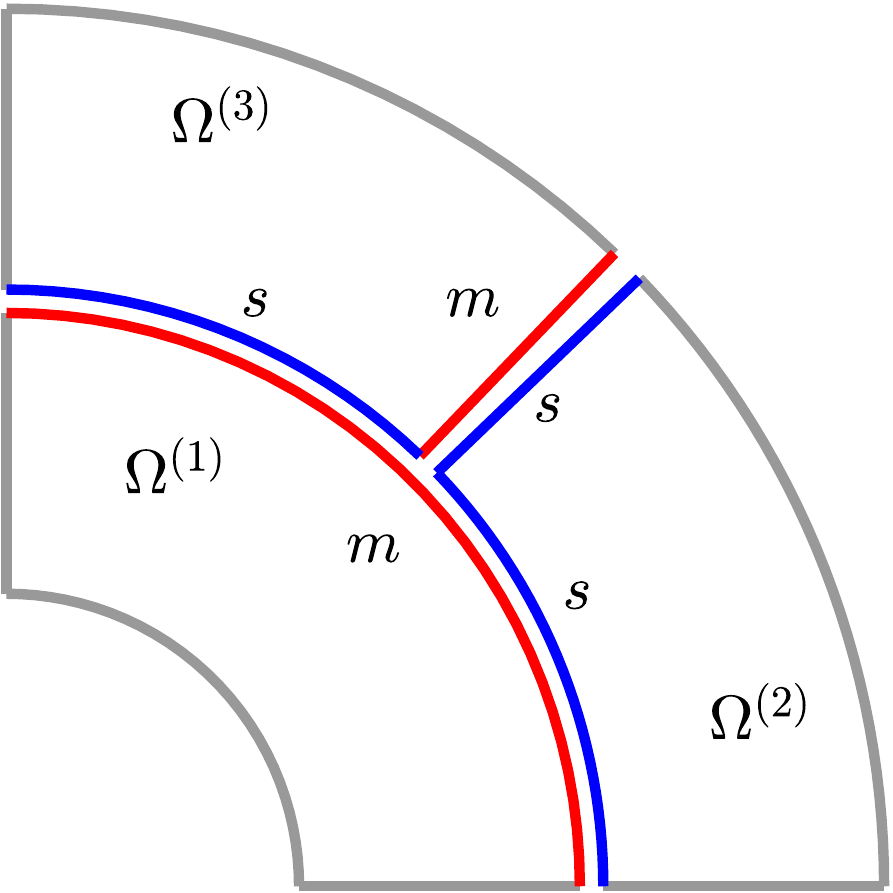}
\end{center}
\caption{The physical domain $\Omega\subset{\mathbb R}^2$
 split in 3 patches
 and the faces $\gamma_{k,\alpha}$. The interfaces $\Gamma_{k\ell}$ are:
$\Gamma_{12}=\gamma_{1,1}\cap\gamma_{2,2}$,
$\Gamma_{13}=\gamma_{1,1}\cap\gamma_{3,2}$,
$\Gamma_{23}=\gamma_{2,1}\cap\gamma_{3,1}$. At right, a possible
choice of the master/slave faces is shown}
\label{fig:3domini}
\end{figure}

\begin{figure}
\begin{center}
\includegraphics[trim={0 20 0 0},width=0.4\textwidth]{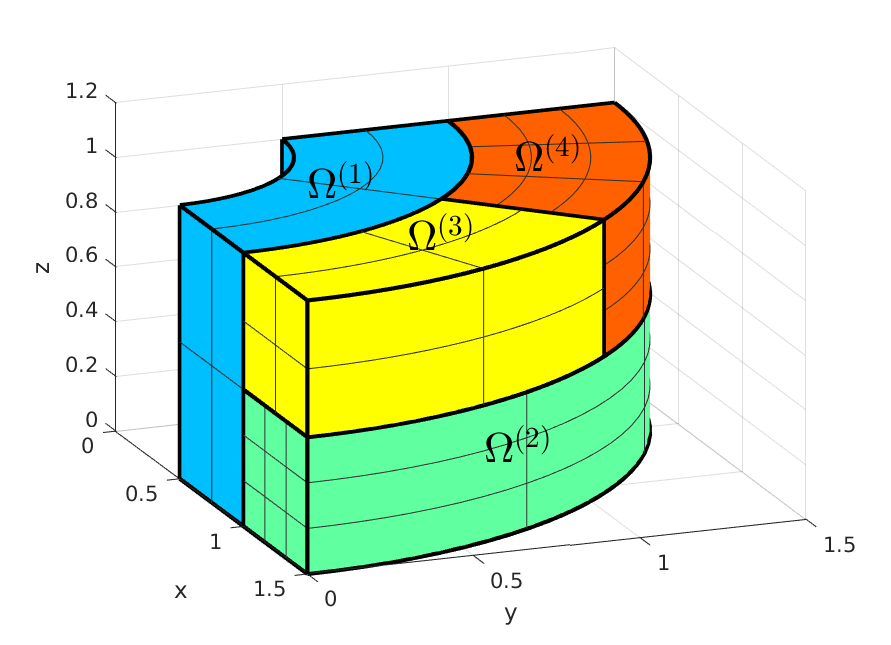}\quad
\includegraphics[trim={0 20 0 0},width=0.4\textwidth]{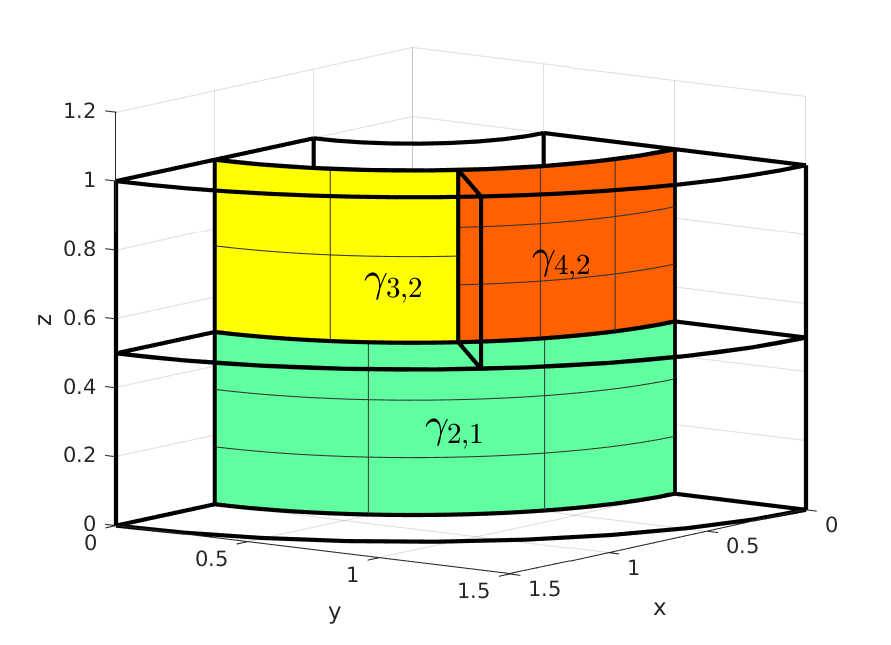}\quad

\end{center}
\caption{3D-decompositions.
$\gamma_{1,1}$ is the sole interface of $\Omega^{(1)}$ and it
has non-null intersection
with just one interface of each one of the other subdomains.
$\Omega^{(2)}$ has two interfaces:
$\gamma_{2,1}$ (that adjacent to $\Omega^{(1)}$)
and $\gamma_{2,2}$ (the top one adjacent to both $\Omega^{(3)}$ and
$\Omega^{(4)}$);
$\Omega^{(3)}$ has three interfaces:
$\gamma_{3,1}$ (that adjacent to $\Omega^{(4)}$),
$\gamma_{3,2}$  (that adjacent to $\Omega^{(1)}$)
and $\gamma_{3,3}$ (that adjacent to $\Omega^{(2)}$);
$\Omega^{(4)}$ has three interfaces:
$\gamma_{4,1}$ (that adjacent to $\Omega^{(3)}$),
$\gamma_{4,2}$  (that adjacent to $\Omega^{(1)}$)
and $\gamma_{4,3}$ (that adjacent to $\Omega^{(2)}$)}
\label{fig:4dom_3d}
\end{figure}

We split the 
internal boundary $\Gamma_k$ of $\partial\Omega^{(k)}$ in faces and
we denote by $\gamma_{k,\alpha}$  
the $\alpha$th face of $\Gamma_k$ (see Fig. \ref{fig:3domini} and Fig.
\ref{fig:4dom_3d}),
the first sub-index $k$ identifies the domain, while the second one $\alpha$ is
the index of the face of $\Gamma_k$.
As for the case of two subdomains, we assume that each interface
$\gamma_{k,\ell}$ is sufficiently regular (i.e. of class $C^{1,1}$)
 to allow the conormal derivative of $u_k$ on
$\gamma_{k,\ell}$ to be well defined.

\begin{remark}
{\rm We assume that $\gamma_{k,\alpha}$ includes its boundary.}
\end{remark}

For example, in the multipatch configuration 
of Fig. \ref{fig:3domini}, we have
$\Gamma_1=\gamma_{1,1}$,
$\Gamma_2=\gamma_{2,1}\cup\gamma_{2,2}$ and $\Gamma_3=
\gamma_{3,1}\cup\gamma_{3,2}$.

Moreover, for any face $\gamma_{k,\alpha}$ we define the set
\begin{equation}\label{S_kalpha}
{\cal A}_{k,\alpha}=\{(\ell,\beta):\
\gamma_{\ell,\beta}\cap\gamma_{k,\alpha}\neq \emptyset\}
\end{equation}
of the faces (of the other patches) that are adjacent to $\gamma_{k,\alpha}$.
In the multipatch configuration of
Fig. \ref{fig:3domini}, we have ${\cal A}_{1,1}=\{(2,2),(3,2)\}$, 
${\cal A}_{2,2}=\{(1,1)\}$,  ${\cal A}_{2,1}=\{(3,1)\}$ and so on.

%

Between $\gamma_{k,\alpha}$ and $\gamma_{\ell,\beta}$, one is tagged
as \emph{master} and the other as \emph{slave}
and we define the \emph{master skeleton}
\begin{equation}\label{global_interface}
\Gamma=\bigcup_{(k,\alpha)} \gamma_{k,\alpha} \quad\mbox{ with }
\gamma_{k,\alpha} \mbox{ master.}
\end{equation}
In the mortar community $\Gamma$ is named \emph{mortar interface}.

A-priori there is no constraint in tagging a face as either master
or slave.
In the example of Fig. \ref{fig:3domini} right, we could tag as master
the face $\gamma_{1,1}$ (in which case
$\gamma_{2,2}$ and $\gamma_{3,2}$ will be both slave), or other way
around.

In the patch $\Omega^{(k)}$ we define a NURBS space ${\cal N}^{(k)}_{h_k}$ as
defined in (\ref{spaces_Nk})  and the corresponding
finite dimensional spaces $V^{(k)}_{h_k}$ (see (\ref{spaces1})) that
are totally independent of the discretizations inside the
adjacent patches. For each face $\gamma_{k,\alpha}\subset \Gamma_k$ we define 
the trace space 
$$Y_{h_k}^{(k,\alpha)}=\{\lambda=v|_{\gamma_{k,\alpha}}, \ 
v\in {\cal N}^{(k)}_{h_k}\}$$
whose dimension is denoted by $n^{(k,\alpha)}$.

The INTERNODES method for $M>2$ patches reads as follows.
For $k=1,\ldots,M$, let $\widetilde g_{h_k}^{(k)}\in{\cal N}^{(k)}_{h_k}$ be
a suitable approximation of $\widetilde g^{(k)}$,  we look for
$u^{(k)}_{h_k}\in{\cal N}^{(k)}_{h_k}$ such that $(u_{h_k}^{(k)}-\widetilde 
g_{h_k}^{(k)})\in V^{(k)}_{h_k}$ and
\begin{eqnarray}\label{internodes_weakM}
\left\{
\begin{array}{ll}
a^{(k)}(u_{h_k}^{(k)},v_{h_k}^{(k)}) = {\cal F}^{(k)}(v_{h_k}^{(k)})
&\ \forall v_{h_k}^{(k)}\in V^{(k)}_{0,h_k},\\[3mm]
\mbox{for any }(\ell,\beta): \gamma_{\ell,\beta} \mbox{ is slave}\\
\displaystyle\qquad u^{(\ell)}_{h_\ell}=\sum_{(k,\alpha)\in {\cal
A}_{\ell,\beta}}
\Pi_{(\ell,\beta)(k,\alpha)}u^{(k)}_{h_k} & 
\mbox{ on }\gamma_{\ell,\beta}\\[3mm]
\mbox{for any }(k,\alpha): \gamma_{k,\alpha} \mbox{ is master}\\
\qquad \displaystyle\langle \nu_k\frac{\partial u^{(k)}_{h_k}}{\partial{\bf n}_k}
+\sum_{(l,\beta)\in {\cal A}_{k,\alpha}}
\widetilde\Pi_{(k,\alpha)(\ell,\beta)}(\nu_\ell\frac{\partial
u^{(\ell)}_{h_\ell}}
{\partial{\bf n}_\ell})
,\eta\rangle = 0 &\
 \forall \eta \in Y^{(k,\alpha)}_{h_k},
\end{array}\right.
\end{eqnarray}
where:
\begin{itemize}[noitemsep]
\item $\Pi_{(k,\alpha)(\ell,\beta)}$ and $\Pi_{(\ell,\beta)(k,\alpha)}$ are the 
interpolation operators used to 
transfer information from one side to the other of
$\gamma_{k,\alpha}\cap\gamma_{\ell,\beta}\neq\emptyset$, more precisely, 
$\Pi_{(k,\alpha)(\ell,\beta)}$ moves from $\gamma_{\ell,\beta}$ to 
$\gamma_{k,\alpha}$, while
$\Pi_{(\ell,\beta)(k,\alpha)}$ moves from $\gamma_{k,\alpha}$ to
$\gamma_{\ell,\beta}$;
\item
${\cal J}_{k,\alpha}$ (${\cal J}_{\ell,\beta}$, resp.)
denotes the canonical isomorphism between 
$Y_{h_k}^{(k,\alpha)}$ ($Y_{h_k}^{(\ell,\beta)}$, resp.) and its dual space;
\item
$\widetilde\Pi_{(k,\alpha)(\ell,\beta)}={\cal J}_{k,\alpha}
\Pi_{(k,\alpha)(\ell,\beta)}{\cal J}_{\ell,\beta}^{-1}$.
\end{itemize}

When the cardinality of ${\cal A}_{k,\alpha}$ is equal to one, i.e. 
there is only one face adjacent to $\gamma_{k,\alpha}$, then
the summations in (\ref{internodes_weakM})$_{2,3}$ disappear and we recover
the interface conditions (\ref{internodes_weak})$_{2,3}$. In such a case 
 the definition of the interpolation operators 
 is as in either Sect. \ref{sec:nurbs_interpolation} or 
Sect. \ref{sec:rbf_interpolation}.

\begin{figure}
\begin{center}
\includegraphics[trim={0 30 0
10},width=0.3\textwidth]{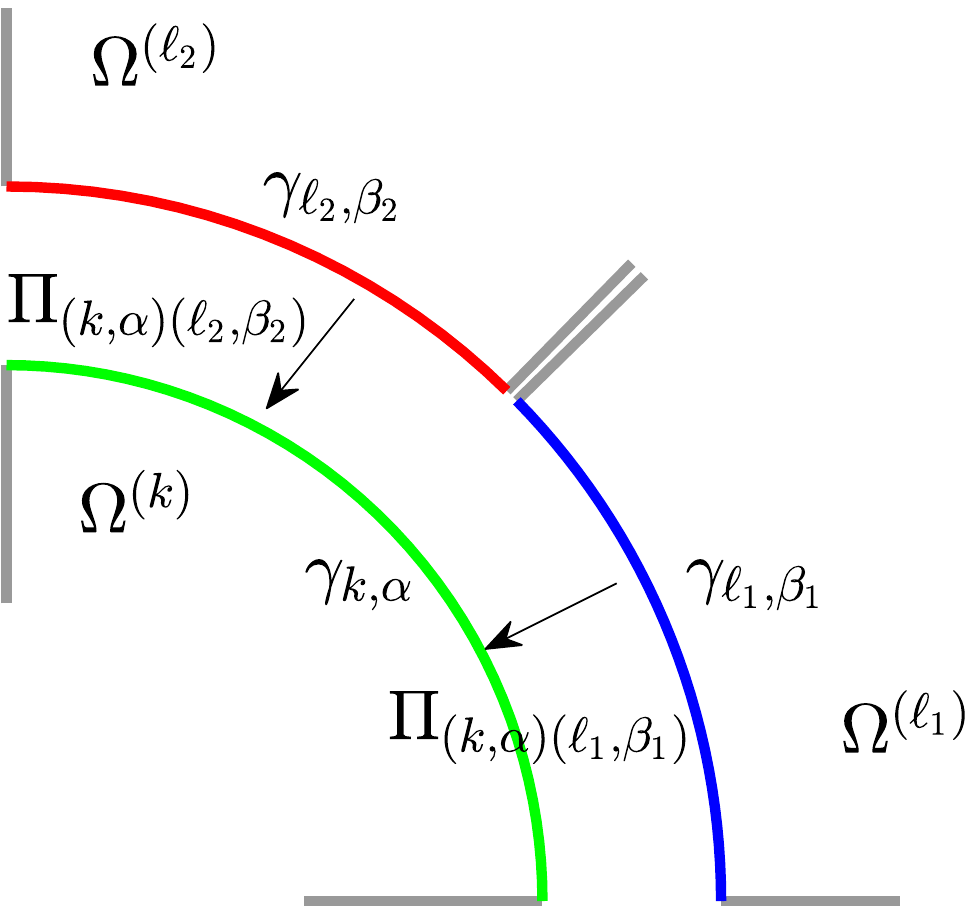}
\qquad
\includegraphics[trim={0 30 0
10},width=0.3\textwidth]{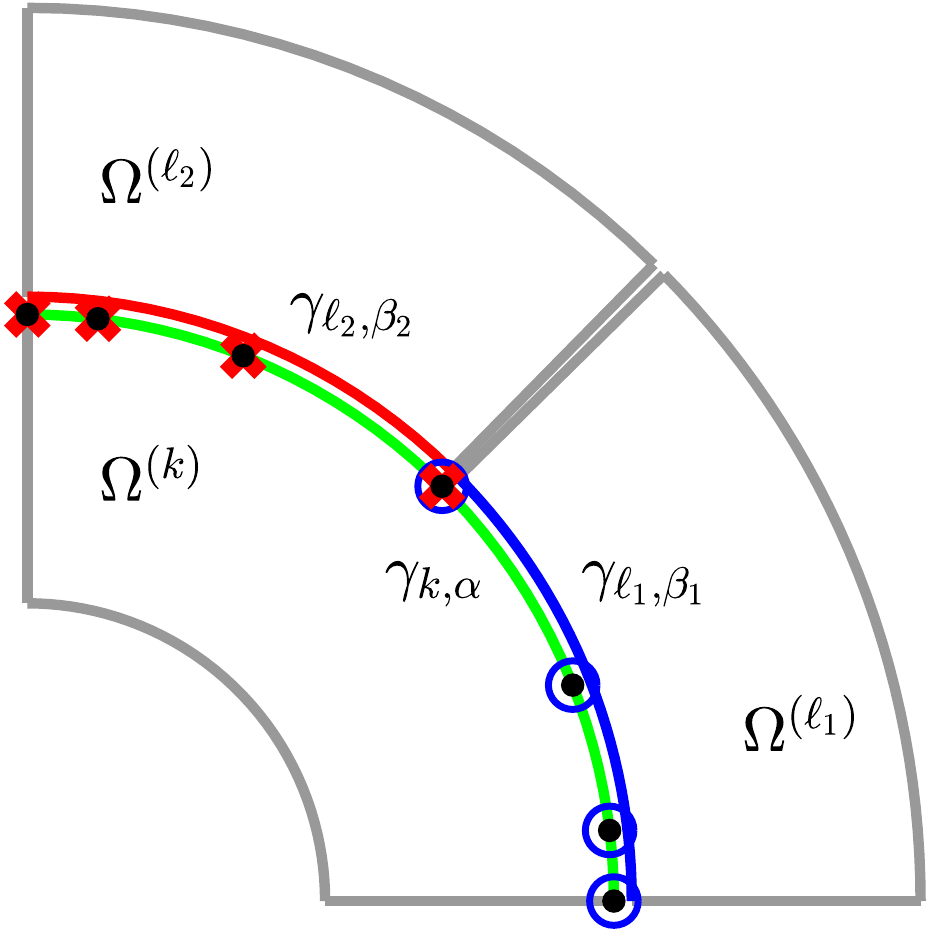}
\end{center}
\caption{At left, interpolation from
$\left(\gamma_{\ell_1,\beta_1}\cup\gamma_{\ell_2,\beta_2}\right)$ to
$\gamma_{k,\alpha}$. Here ${\cal
A}_{k,\alpha}=\{(\ell_1,\beta_1),(\ell_2,\beta_2)\}$. At right, 
the black dots are the Greville nodes ${\bf
x}_{i,G}^{(k,\alpha)}\in\gamma_{k,\alpha}$. The indices of the nodes
surrounded by blue rings belong to ${\cal G}_{(k,\alpha)(\ell_1,\beta_1)}$,
while the indices of 
the nodes marked with a red cross belong to 
${\cal G}_{(k,\alpha)(\ell_2,\beta_2)}$ (see formula (\ref{Gkl})). There is
one Greville node that belongs to both the sets ${\cal
G}_{(k,\alpha)(\ell_1,\beta_1)}$ and ${\cal G}_{(k,\alpha)(\ell_2,\beta_2)}$,
the corresponding value $(U_{k,\alpha})_i$ is equal to 1/2. For all the other
Greville nodes, $(U_{k,\alpha})_i$ is equal to 1}
\label{fig:multi_interp}
\end{figure}

When instead the cardinality of ${\cal A}_{k,\alpha}$ is
greater than one, that is the face
$\gamma_{k,\alpha}$ interfaces with at least two adjacent faces (as, 
e.g. for  the face $\gamma_{k,\alpha}$ 
in Fig. \ref{fig:multi_interp}), and we want to interpolate from 
$\displaystyle\bigcup_{(\ell,\beta)\in {\cal
A}_{k,\alpha}}\gamma_{\ell,\beta}$
to $\gamma_{k,\alpha},$ then we have to slightly modify the
 definition of the interpolation operators. In Sect. \ref{sec:interpM} we
will show how to extend the definition of the interpolation matrices 
(\ref{matrix_P21}) and (\ref{matrix_P12}) when the interfaces are watertight.
If two patches are non-watertight, we assume that the cardinality of the set 
${\cal A}_{k,\alpha}$ associated with the corresponding interfaces 
is equal to one, thus the interpolation matrices are built as in
(\ref{matrix_P21_RBF}) and (\ref{matrix_P12_RBF}).

Finally, in Sect. \ref{sec:normalderM} we will precise how to generalize formula
(\ref{residual}) for the computation of normal derivatives.

\subsection{Extension of the interpolation matrices (\ref{matrix_P21}) and
(\ref{matrix_P12}) for watertight interfaces. }\label{sec:interpM}

We denote by ${\bf F}^{(k,\alpha)}:{\mathbb R}^{d-1}\to {\mathbb R}^d$ the
restriction of ${\bf F}^{(k)}$ to the face $\widehat\gamma_{k,\alpha}=
({\bf F}^{(k)})^{-1}(\gamma_{k,\alpha})$.
Let  ${\bf x}_{i,G}^{(k,\alpha)}$, for $i=1,\ldots, n^{(k,\alpha)}$,  be
the Greville nodes associated with the patch $\Omega^{(k)}$ that belong to
$\gamma_{k,\alpha}$ and $\widehat{\bf x}_{i,G}^{(k,\alpha)}=
({\bf F}^{(k,\alpha)})^{-1}({\bf x}_{i,G}^{(k,\alpha)})$
their counter-image in the face $\widehat\gamma_{k,\alpha}$ 
of the parameter domain.
For any face $\gamma_{k,\alpha}$ and for any $(\ell,\beta)\in {\cal
A}_{k,\alpha}$
we define the set
\begin{equation}\label{Gkl}
{\cal G}_{(k,\alpha)(\ell,\beta)}=\{i=1,\ldots,n^{(k,\alpha)}|\ 
{\bf x}_{i,G}^{(k,\alpha)}\in \gamma_{\ell,\beta}\}
\end{equation}
of the indices of the Greville nodes associated with the domain 
$\Omega^{(k)}$
belonging to
$\gamma_{k,\alpha}$ that lay on $\gamma_{\ell,\beta}$ too (see Fig.
\ref{fig:multi_interp}, right).

Notice that ${\cal G}_{(k,\alpha)(\ell,\beta)}$ and ${\cal
G}_{(\ell,\beta)(k,\alpha)}$  denote two different sets.

Finally,  we define a sort of \emph{partition of unity function}, with the aim
of interpolating correctly the data coming from two contiguous faces adjacent to
 $\gamma_{k,\alpha}$.

For any $i=1,\ldots,n^{(k,\alpha)}$, we set 
\begin{equation}\label{pu}
(U_{k,\alpha})_i=1/\mbox{card}\{(\ell,\beta)\in{\cal A}_{k,\alpha}|\ i\in {\cal
G}_{(k,\alpha)(\ell,\beta)}\},
\end{equation}
that is $(U_{k,\alpha})_i$ is the inverse of the
number of faces adjacent to $\gamma_{k,\alpha}$  which the 
Greville node ${\bf x}_{i,G}^{(k,\alpha)}$ lays on 
(up to the boundary). Here, card$A$ denotes the cardinality of
the set $A$.

Given $\lambda\in Y^{(\ell,\beta)}_{h_\ell}$ , we define
the interpolation operator
$\Pi_{(k,\alpha)(\ell,\beta)}:Y^{(\ell,\beta)}_{h_\ell}\to
Y^{(k,\alpha)}_{h_k}$ by imposing the following interpolation conditions
for $i=1,\ldots, n^{(k,\alpha)}$:
\begin{eqnarray}\label{Pikl2}
\begin{array}{lll}
(\Pi_{(k,\alpha)(\ell,\beta)}\lambda)({\bf x}_{i,G}^{(k,\alpha)})
=\left\{\begin{array}{ll}
(U_{k,\alpha})_i\ \lambda({\bf
x}_{i,G}^{(k,\alpha)}) & \mbox{ if }i\in{\cal G}_{(k,\alpha)(\ell,\beta)},\\
0 & \mbox{ if }i\not\in{\cal G}_{(k,\alpha)(\ell,\beta)}
\end{array}\right.
%
\end{array}
\end{eqnarray}

Notice that
 $\Pi_{(k,\alpha)(\ell,\beta)}\lambda$ is defined on the whole face 
$\gamma_{k,\alpha}$, even when $\gamma_{\ell,\beta}\subsetneqq
\gamma_{k,\alpha}$.

\begin{figure}
\begin{center}
\includegraphics[width=0.3\textwidth]{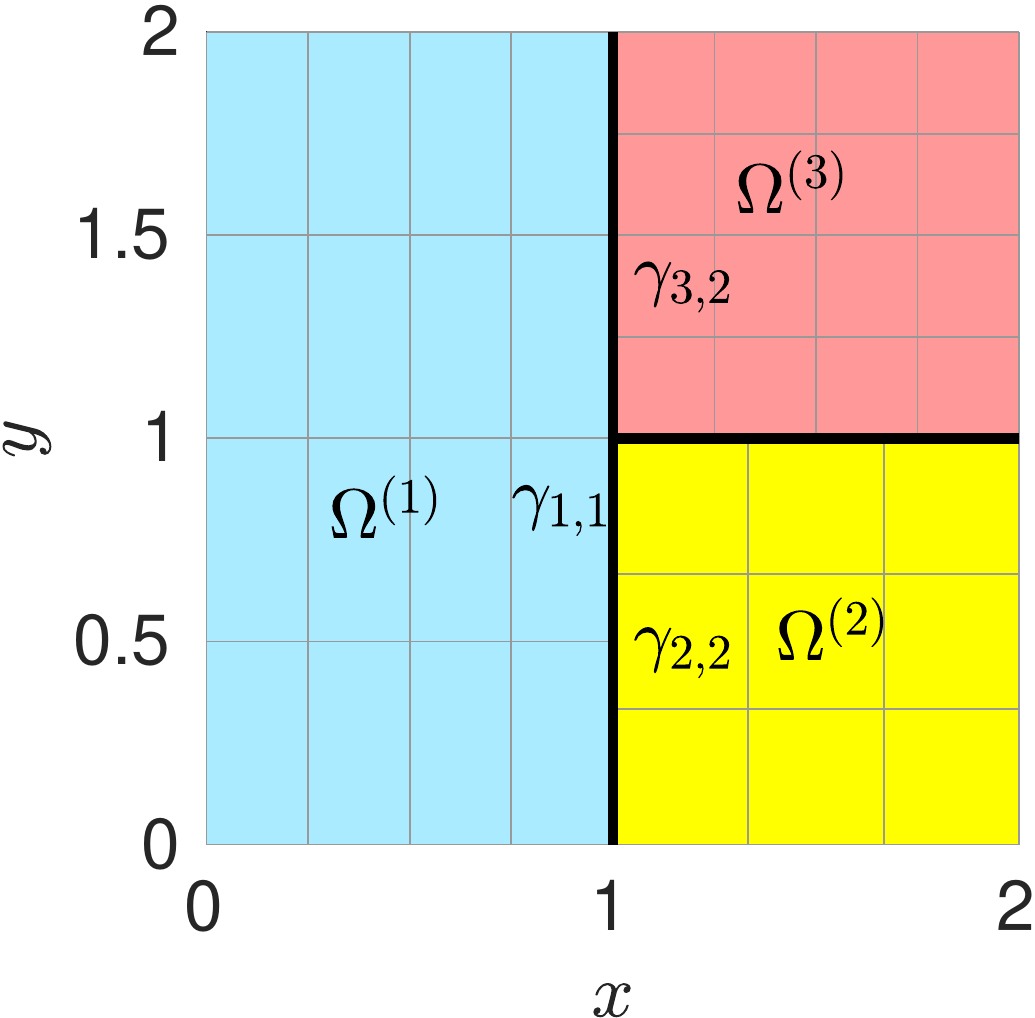}\qquad
\includegraphics[width=0.35\textwidth]{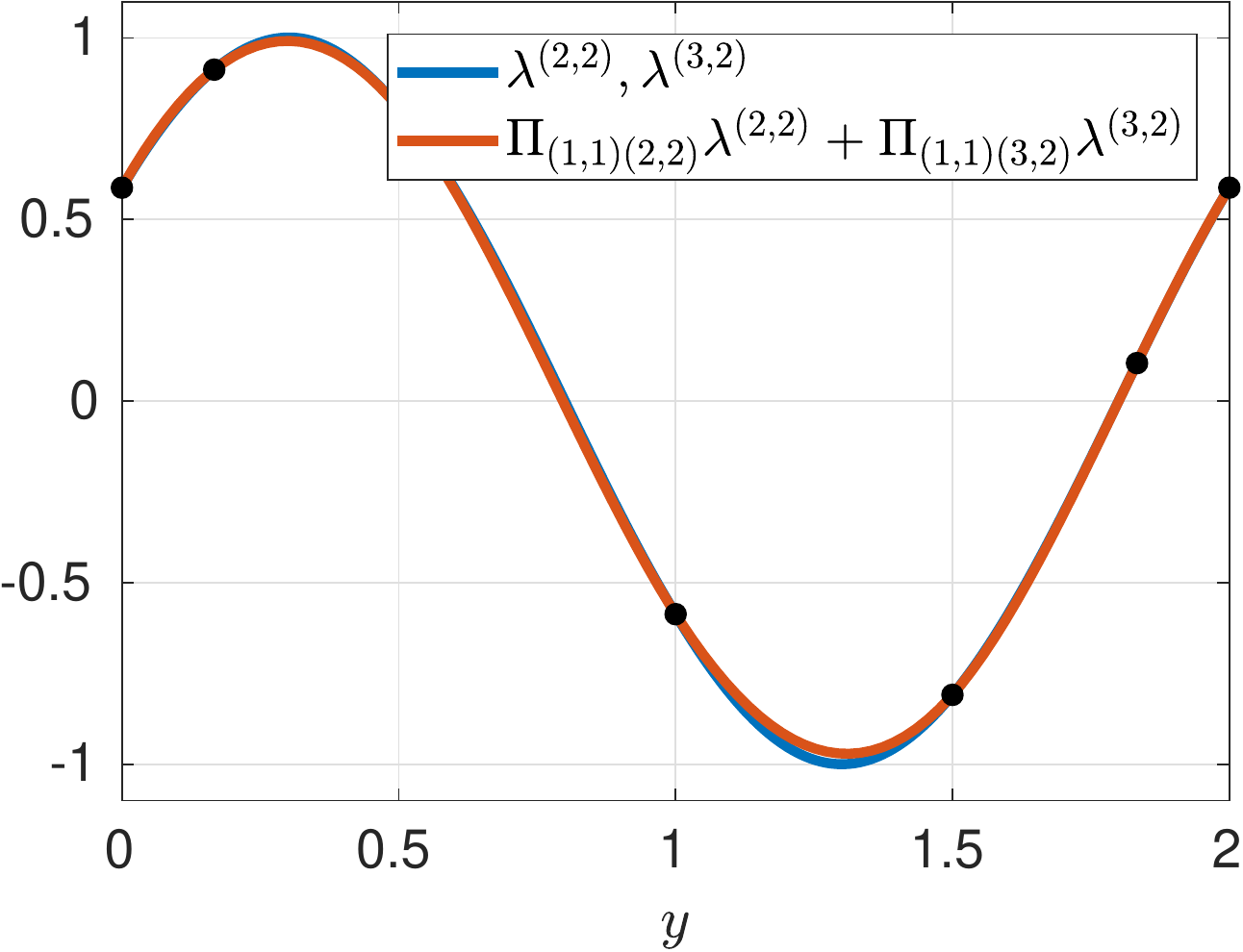}\\[2mm]
\includegraphics[trim={0 20 0
0},width=0.35\textwidth]{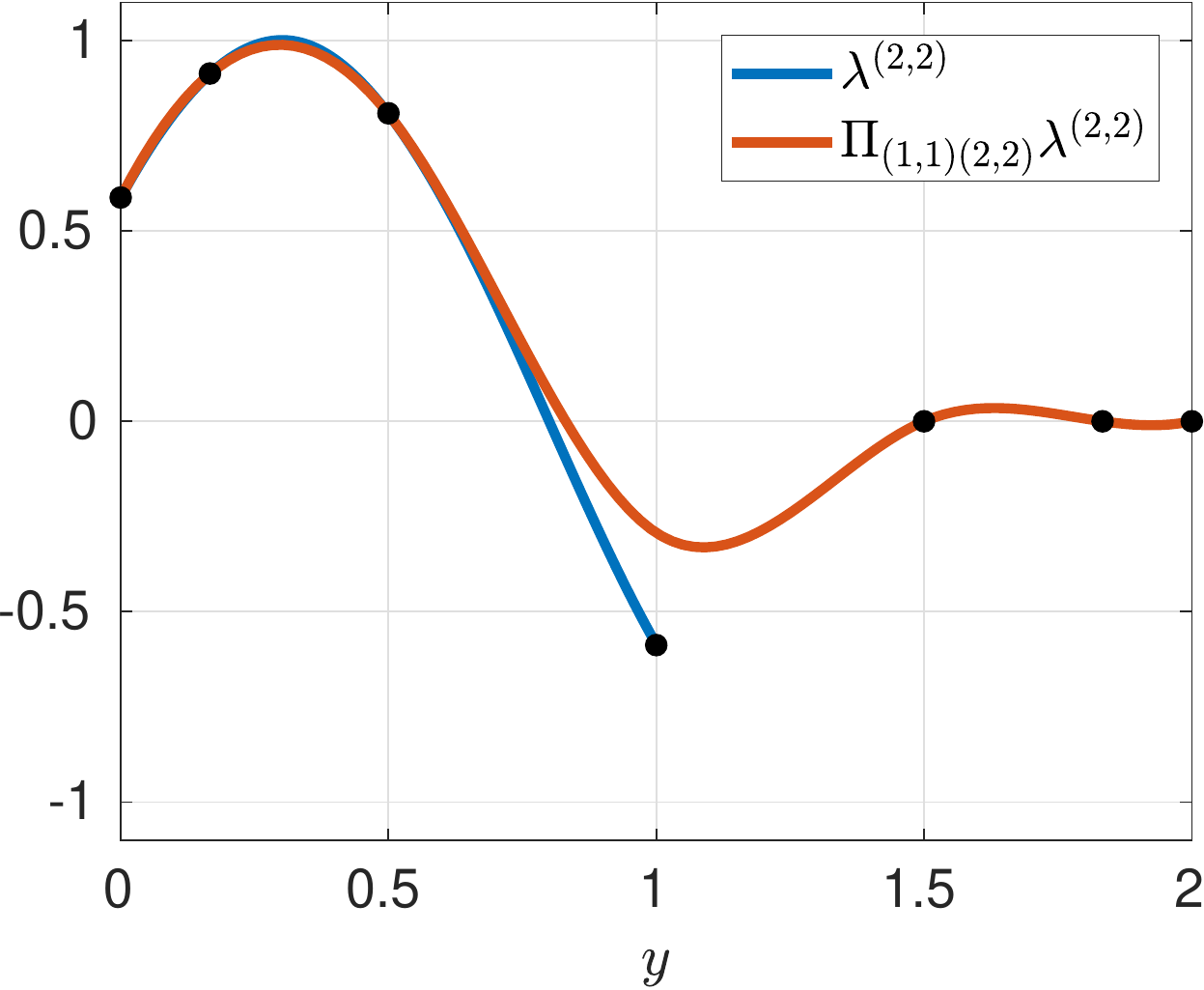}\quad
\includegraphics[trim={0 20 0
0},width=0.35\textwidth]{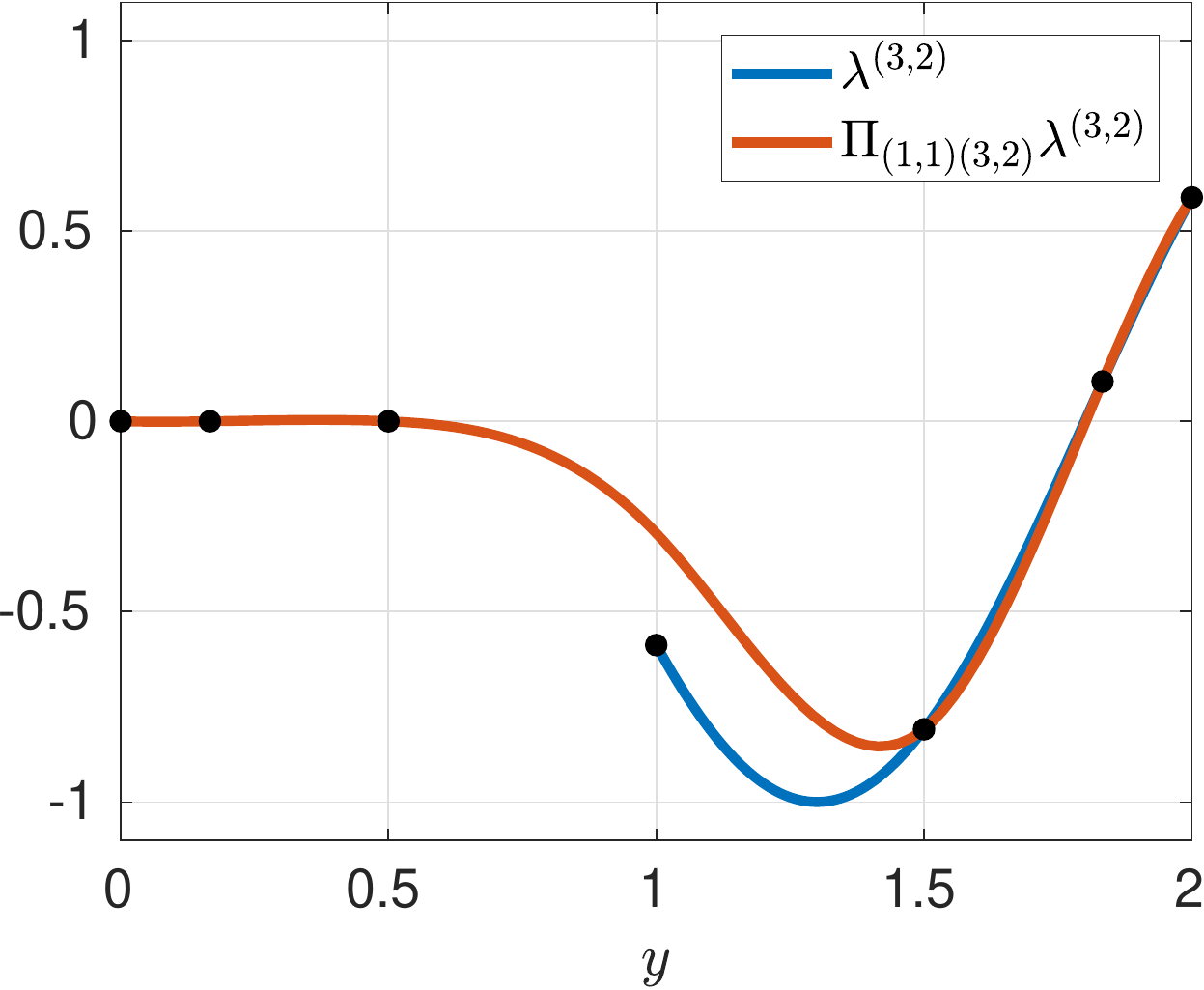}
\end{center}
\caption{Example of how the interpolation operators
$\Pi_{(k,\alpha)(\ell,\beta)}$ work. 
 The black dots represent the interpolated values at the Greville nodes}
\label{fig:domain123}
\end{figure}

Let us consider the  multipatch geometry shown in the top left picture of Fig.
\ref{fig:domain123}, we interpolate a trace function from
$\gamma_{2,2}\cup\gamma_{3,2}$ to $\gamma_{1,1}$. 
In the bottom pictures of Fig. \ref{fig:domain123} we show how the interpolation
operators $\Pi_{(1,1)(2,2)}$ (from $\gamma_{2,2}$ to $\gamma_{1,1}$) and 
 $\Pi_{(1,1)(3,2)}$ (from $\gamma_{3,2}$ to $\gamma_{1,1}$) work.
The point whose coordinates are $(x,y)=(1,1)$ is a Greville point for
$\Omega^{(1)}$ belonging to $\gamma_{1,1}$ and it 
lays on both $\gamma_{2,2}$ and $\gamma_{3,2}$ too (the faces include
their boundary), thus the corresponding
weight defined in (\ref{pu}) is equal to 1/2. 
Notice that $(\Pi_{(1,1)(2,2)}\lambda)$ takes null value at the 
Greville nodes of $\gamma_{1,1}$ 
not laying on $\gamma_{2,2}$. Analogous considerations hold
for $(\Pi_{(1,1)(3,2)}\lambda)$.

The sum $\Pi_{(1,1)(2,2)}\lambda^{(2,2)}+\Pi_{(1,1)(3,2)}\lambda^{(3,2)}$
interpolates the piece-wise function  $\lambda$ such that
$\lambda_{|\gamma_{2,2}}=\lambda^{(2,2)}$ and 
$\lambda_{|\gamma_{3,2}}=\lambda^{(3,2)}$ (see the right picture of Fig.
\ref{fig:domain123}).

Let $\mu_j^{(k,\alpha)}$,  with 
$j=1,\ldots,n^{(k,\alpha)}$, be
the NURBS basis functions of the trace space
$Y_{h_k}^{(k,\alpha)}$  and $\mu_j^{(k,\alpha)}$ the corresponding basis
function in the parameter domain, i.e. such that $\mu_j^{(k,\alpha)}=
\widehat\mu_j^{(k,\alpha)}\circ({\bf F}^{(k,\alpha)})^{-1}$.
By setting
\begin{eqnarray}
\begin{array}{lll}
(G_{(k,\alpha)(k,\alpha)})_{ij}&=\widehat\mu_j^{(k,\alpha)}(\widehat{\bf
x}_{i,G}^{(k,\alpha)}), & i,j=1,\ldots,n^{(k,\alpha)},\\[2mm]
(G_{(k,\alpha)(\ell,\beta)})_{ij}&=\left\{\begin{array}{ll}
\widehat\mu_j^{(\ell,\beta)}
(({\bf F}^{(\ell,\beta)})^{-1}({\bf x}_{i,G}^{(k,\alpha)})) 
& \mbox{ if }i\in{\cal G}_{(k,\alpha)(\ell,\beta)}\\
0 & \mbox{ if }i\not\in{\cal G}_{(k,\alpha)(\ell,\beta)}
\end{array}\right.,
& i=1,\ldots,n^{(k,\alpha)},\ j=1,\ldots,n^{(\ell,\beta)}
\end{array}
\end{eqnarray}
the matrix associated with $\Pi_{(k,\alpha)(\ell,\beta)}$ is
\begin{equation}\label{P_kalb}
P_{(k,\alpha)(\ell,\beta)}=G_{(k,\alpha)(k,\alpha)}^{-1}\mbox{diag}(U_{k,\alpha})
G_{(k,\alpha)(\ell,\beta)}.
\end{equation}

Similarly, we define the interpolation operator
$\Pi_{(\ell,\beta)(k,\alpha)}:Y^{(k,\alpha)}_{h_k}\to
Y^{(\ell,\beta)}_{h_\ell}$ and the corresponding matrix
\begin{equation}
P_{(\ell,\beta)(k,\alpha)}=G_{(\ell,\beta)(\ell,\beta)}^{-1}\mbox{diag}
(U_{\ell,\beta})
G_{(\ell,\beta)(k,\alpha)}.
\end{equation}

Denoting by $\boldsymbol\lambda_{k,\alpha}$ ($\boldsymbol\lambda_{\ell,\beta}$,
resp.) the array whose entries are
the degrees of freedom of the function $u^{(k)}_{h_k}$ 
($u^{(\ell)}_{h_\ell}$, resp.) associated with
the face $\gamma_{k,\alpha}$ ($\gamma_{\ell,\beta}$, resp.),
the algebraic implementation of the interface conditions
(\ref{internodes_weakM})$_{2}$ reads:
\begin{eqnarray}\label{summa1}
\begin{array}{ll}
\mbox{for any }(\ell,\beta): \gamma_{\ell,\beta} \mbox{ is slave,}
\displaystyle\qquad \boldsymbol\lambda_{\ell,\beta}=\sum_{(k,\alpha)\in
{\cal A}_{\ell,\beta}} P_{(\ell,\beta),(k,\alpha)}
\boldsymbol\lambda_{k,\alpha}.
\end{array}
\end{eqnarray}

\subsection{Transferring normal derivatives}\label{sec:normalderM}

For each face (either slave or master)
$\gamma_{k,\alpha}\subset\Gamma_k=\partial\Omega^{(k)}\setminus\partial\Omega$,
we define the set
$G_{k,\alpha}=\partial\Omega^{(k)}\setminus\gamma_{k,\alpha}$
(or $G_{k,\alpha}=\partial\Omega^{(k)}\setminus(\gamma_{k,\alpha}\cup
\partial\Omega^{(k)}_N)$ when mixed boundary conditions are given on
$\partial\Omega$)
and the real values (similar to (\ref{residual}))
\begin{equation}\label{residualM}
r^{(k,\alpha)}_i=\int_{\gamma_{k,\alpha}} \nu_k\frac{\partial u^{(k)}_{h_k}}{\partial {\bf n}_k}
\,\overline{\cal L}^{(k)}\mu_i^{(k,\alpha)}\,d\Gamma, \quad
i=1,\ldots,n^{(k,\alpha)}.
\end{equation}

As  done for the configuration with only two patches, we can compute 
$r^{(k,\alpha)}_i$ by exploiting the weak form of the differential equations
inside the patches, i.e.
\begin{equation}\label{fluxM}
r^{(k,\alpha)}_i=a^{(k)}(u^{(k)}_{h_k},
\overline{\cal L}^{(k)}\mu_i^{(k,\alpha)})-{\cal F}^{(k)}
(\overline{\cal L}^{(k)}\mu_i^{(k,\alpha)})
-\int_{G_{k,\alpha}} \nu_k\frac{\partial u^{(k)}_{h_k}}{\partial {\bf n}_k}
\,\overline{\cal L}^{(k)}\mu_i^{(k,\alpha)}\,d\Gamma, \quad
i=1,\ldots,n^{(k,\alpha)}.
\end{equation}

Following the notations of Sect. \ref{sec:alg},
for any $k=1,\ldots,M$  and for any face $\gamma_{k,\alpha}\subset \Gamma_k$,
we define the sets ${\cal I}_{\overline\gamma_{k,\alpha}}$, 
${\cal I}_{\gamma_{k,\alpha}}$, and  ${\cal I}_{\partial\gamma_{k,\alpha}}$.
To evaluate the last integral of (\ref{fluxM})
we define the matrix $C^{(k,\alpha)}$ (of size $N^{(k)}\times N^{(k)}$)
whose non-null entries are
\begin{eqnarray}\label{matrix_corrM}
\begin{array}{lll}
C^{(k,\alpha)}_{ij}&
 \displaystyle =-\int_{G_{k,\alpha}} \nu_k\frac{\partial\varphi_j ^{(k)}}
{\partial{\bf n}_k}
\varphi_i^{(k)}, & \mbox{ for }i\in{\cal I}_{\partial\gamma_{k,\alpha}},
\ j\in{\cal I}_{\overline\Omega^{(k)}}.
\end{array}
\end{eqnarray}
Then we set
\begin{equation}\label{AC}
\widehat A^{(\overline\gamma_{k,\alpha},\overline\Omega^{(k)})}= 
A^{(\overline\gamma_{k,\alpha},\overline\Omega^{(k)})}
+C^{(\overline\gamma_{k,\alpha},\overline\Omega^{(k)})}
\end{equation}
and ${\bf f}_{k,\alpha}={\bf f}({\cal I}_{\overline\gamma_{k,\alpha}})$,
and we compute
$$ {\bf r}_{k,\alpha}=\widehat
A^{(\overline\gamma_{k,\alpha},\overline\Omega^{(k)})} {\bf u}^{(k)}
-{\bf f}_{k,\alpha},$$
that will contain the values $r^{(k,\alpha)}_i$ defined in (\ref{fluxM}).

Finally, for any face $\gamma_{k,\alpha}$ we define the mass matrix 
\begin{equation}\label{local_massM}
(M_{(k,\alpha)})_{ij}=(\mu_j^{(k)},\mu_i^{(k)})_{L^2(\gamma_{k,\alpha})}, 
\qquad i,j=1,\ldots, n^{(k,\alpha)}.
\end{equation}

The algebraic implementation of the interface conditions
(\ref{internodes_weakM})$_{3}$ reads:
\begin{eqnarray}\label{summa2}
\begin{array}{ll}
\mbox{for any }(k,\alpha): \gamma_{k,\alpha} \mbox{ is master,}
\qquad \displaystyle{\bf r}_{k,\alpha}+\sum_{(\ell,\beta)\in{\cal A}_{k,\alpha}}
M_{(k,\alpha)}P_{(k,\alpha),(\ell,\beta)}M_{(\ell,\beta)}^{-1}
{\bf r}_{\ell,\beta}={\bf 0}.
\end{array}
\end{eqnarray}

\subsection{Comparison between INTERNODES and mortar methods}\label{sec:mortar}

The mortar method is a well-established technology for the coupling of
non-conforming discretizations. It suits also for problems where the
non-conformities are intrinsic, such as the contact problems. Typically, in
mortar methods the continuity between interfaces is imposed weakly, so that a
procedure for integrating the product between basis functions of 
the master and those of the slave sides is required. 
In order to obtain optimal convergence rates, it is
critical to define  quadrature rules that compute very accurately, if not ever
exactly, the cross mass matrix at the interface \cite{brivadis_mortar_quad}. 
Adapting the quadrature points to the elements of the 
master side does not results in an optimal quadrature rule for the basis
functions of the slave side, and vice-versa. 
A common strategy for overcoming this issue in the case of watertight
interfaces consists in 
creating a third "intersection mesh" between the master and slave meshes, in
which the integration elements can be seen as belonging to both sides and thus
the quadrature points defined over them allow for high accurate
 integration of both master and slave basis functions.
Nevertheless, the generation of intersection meshes is far from trivial (from a theoretical
and practical point of view) and may easily result in a number of elements that
is of orders of magnitude greater than those of 
the two original meshes.  Moreover, in
the case of non-watertight interfaces the task acquires even greater complexity
due to the need to project the quadrature nodes from one interface to the
other.

On the contrary, INTERNODES 
requires local interface mass matrices involving only basis
functions from one side of the interface (we refer to matrices
(\ref{local_mass})), and they can be assembled by
taking into account only the parametrization of the relative patch.

\medskip

Most often mortar methods are formulated as
a saddle point problem by introducing an extra field, the Lagrange
multiplier.
This yields an inf-sup compatibility condition to be
fulfilled in order to ensure well-posedness, and such condition 
affects the choice of the
polynomial degrees of the NURBS spaces.
\cite{brivadis_mortar_15}. On the contrary, INTERNODES does not introduce extra
fields, so that it does not care of this problem and no constraints on the
polynomial degrees are required.

When mortar methods are formulated as a single
field problem (like (\ref{alg_2dom_nc_3x3})), the corresponding algebraic system
 is symmetric, provided the original differential operator
is self-adjoint. 
This property is not fulfilled
by INTERNODES due to the two a-priori different intergrid operators
$\Pi_{12}$ and $\Pi_{21}$.
Nevertheless, the local stiffness
matrices continue to be symmetric and positive definite and we can solve the
local differential subsystems by ad-hoc algebraic methods.

In \cite[Sec. 6]{dfgq}, the authors have compared both the eigenvalues and the iterative 
condition number  of the matrix ${\mathbb A}$ with those of the analogous
matrix built with mortar methods instead of INTERNODES. 
The eigenvalues of the
mortar matrix are all real positive (since the global matrix results
 symmetric positive
definite), whereas those of the INTERNODES matrix feature tiny imaginary parts
that vanish as the step size does. Moreover the iterative condition number of
the two matrices behaves right the same way. Since the eigenvalues of the Schur
complement matrix are strictly connected with those of the original matrix, 
we conclude that also the eigenvalues of the corresponding Schur complement
matrices behave similarly.

\medskip
We notice that using only one intergrid interpolation operator (as, e.g., 
$\Pi_{21}$ jointly with its transpose $\Pi_{21}^t$) would not
guarantee an accurate non-conforming method. The use of only one interpolation
operator would  yield
the so-called \emph{pointwise matching} discussed, e.g., in
\cite{bmp94,bbdmkmp}.

\medskip
About the accuracy,  INTERNODES and mortar methods behaves exactly in the same
way: the $H^1-$broken norm of the error decays optimally when the
mesh-size tends to zero, 
i.e. the error produced by these two coupling methods is
proportional to that of the discretization used inside the local subdomains and
it depends on the regularity of the exact solution of the differential
problem. 
The proof of this result for INTERNODES is under
consideration in the IGA context, while 
in the Finite Element framework the
theoretical analysis has been established in \cite{gq_internodes}.

%
\section{The iterative algorithm to solve
(\ref{internodes_weakM}).}\label{sec:schurM}

We extend to each patch with index $k=1,\ldots, M$ 
the notations on the matrices and on 
the arrays introduced in Sect. \ref{sec:alg}
and we split the degrees of freedom (the unknown coefficients of each
$u_{h_k}^{(k)}$ with respect to the NURBS basis functions) in:
\begin{enumerate}
\item ${\bf u}^{(k)}_0$: the degrees of freedom internal to $\Omega^{(k)}$, 
\item ${\bf g}^{(k)}$: the degrees of freedom associated with the Dirichlet 
boundary $\partial\Omega^{(k)}_D$, if it is not empty,
\item ${\bf u}_\Gamma$: the (not replicated) degrees of freedom associated with
 the master skeleton $\Gamma$ defined in (\ref{global_interface}) deprived of
the degrees of freedom associated with $\Gamma\cap\partial\Omega^{(k)}_D$.
\end{enumerate}

Notice that only the degrees of freedom associated with the vertices of 
the patches  are interpolatory and they are the only degrees of freedom which
we must be careful to not replicate  inside ${\bf u}_\Gamma$.
In this way we automatically
enforce the continuity of the solution at such interpolatory points.

\begin{figure}
\begin{center}
\includegraphics[width=0.35\textwidth]{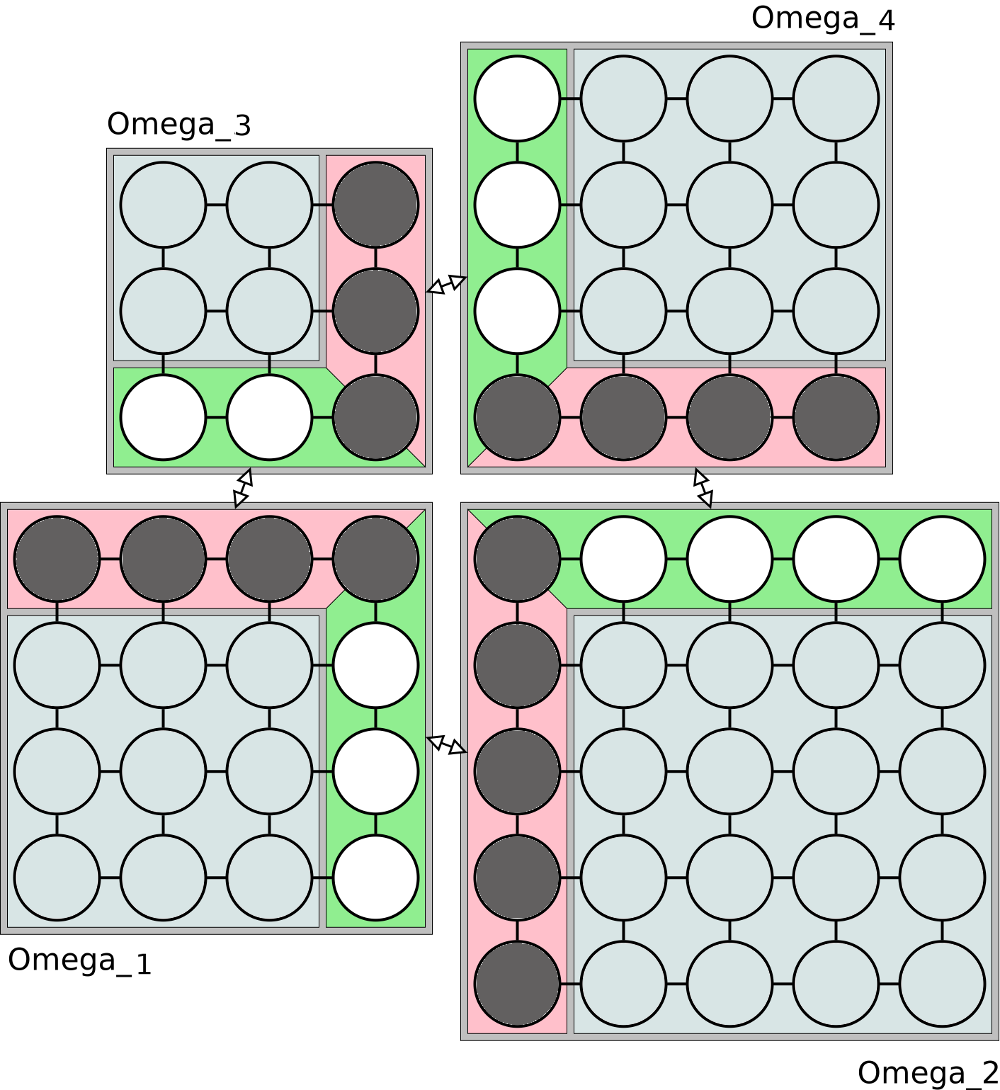}\qquad
\includegraphics[width=0.35\textwidth]{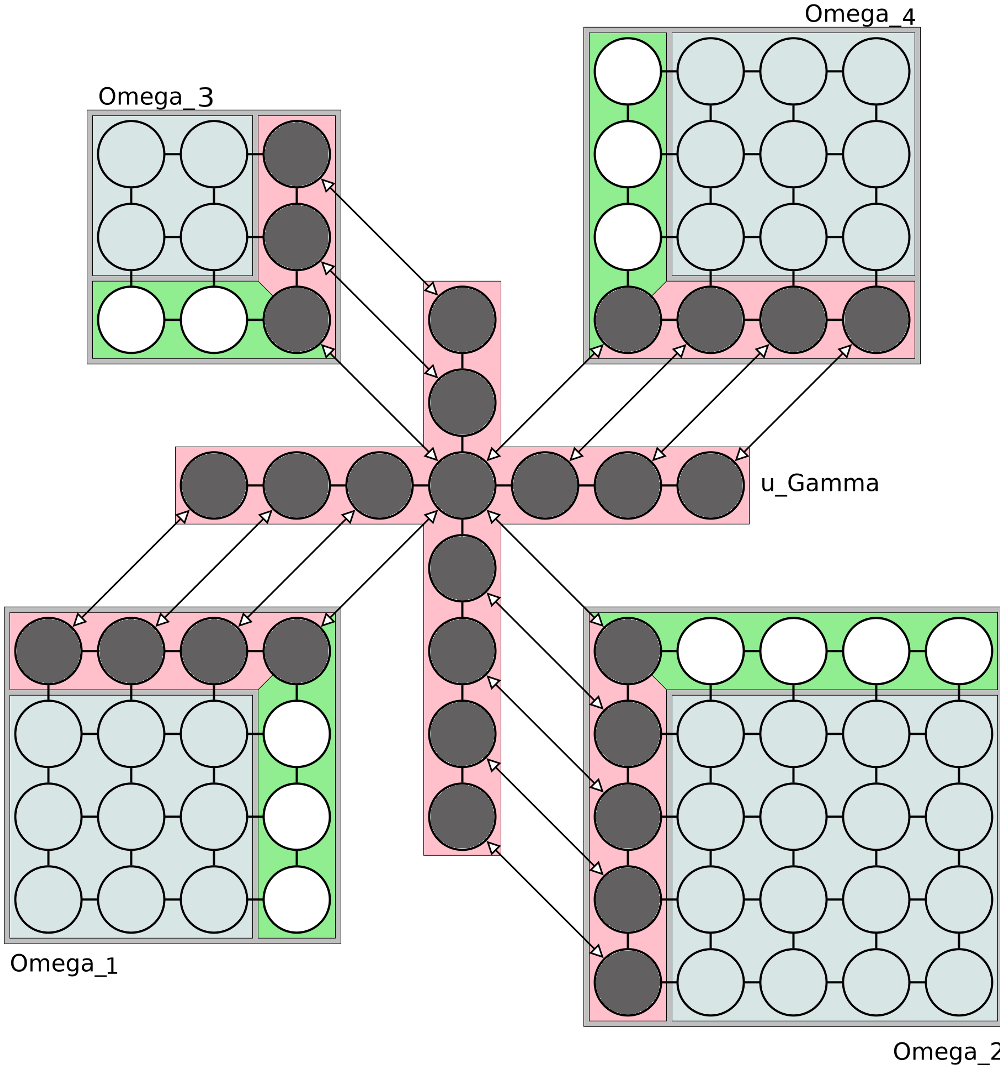}
\end{center}
\caption{At left, the global layout for a situation with 4 patches, the circles
denote the degrees of freedom.
The dark-gray circles are those on the master faces, 
while the white ones are on
the slave faces. 
At right, the degrees of freedom on the master skeleton and
the correspondence between the master skeleton $\Gamma$ and the master faces
$\gamma_{k,\alpha}$
}
\label{fig:ex1}
\end{figure}

For example, when we consider a decomposition like that
sketched in Fig. \ref{fig:ex1}, the vertex shared by all the four patches
belongs to four master faces (the dark-gray ones), but only one occurrence of it
must be considered in the array ${\bf u}_\Gamma$. 

Notice that in the decompositions depicted in Fig. \ref{fig:3domini} the
vertex of $\gamma_{3,1}$ (internal to $\Omega$)
 lays on the face $\gamma_{1,1}$, but 
it is not a degree of freedom for the patch $\Omega^{(1)}$
(the solely interpolation degrees of freedom of $\Omega^{(k)}$ are at the
vertices of the patch itself).

If we eliminate the internal degrees of freedom ${\bf u}^{(k)}_0$, we 
obtain the Schur complement system with
respect to ${\bf u}_\Gamma$:
\begin{equation}\label{SchurM}
S{\bf u}_\Gamma={\bf b}
\end{equation}
and we can solve it by a Krylov method (e.g. Bi-CGStab, GMRES and so on),
for which it is sufficient to provide an algorithm (see
Algorithm \ref{alg:Slambda}) that, given an array
$\boldsymbol\lambda$ of the same size of ${\bf u}_\Gamma$, 
computes $\boldsymbol\psi=S\boldsymbol\lambda$.

Since the array ${\bf u}_\Gamma$ does not contain the degrees 
of freedom associated
with $\Gamma\cap\partial\Omega^{(k)}_D$, but at the same time
the interpolation matrices work on the degrees of freedom associated with the 
faces up to their boundary, for practical purposes it is convenient to
extend ${\bf u}_\Gamma$  to  an array $\widetilde {\bf u}_\Gamma$
that includes also the null
degrees of freedom associated with $\Gamma\cap\partial\Omega^{(k)}_D$.

We denote by $n_\Gamma$ the size of ${\bf u}_\Gamma$, by $\widetilde
n_\Gamma$ the size of $\widetilde {\bf u}_\Gamma$ and we define 
the \emph{restriction matrix} $R_D$ of size $\widetilde n_\Gamma \times
n_\Gamma$ (whose entries are 0 or 1) that, with any array $\widetilde{\bf u}$
defined on $\Gamma$, associates its restriction to
$\Gamma\setminus\partial\Omega$ such that ${\bf u}_\Gamma=R_D\widetilde{\bf
u}_\Gamma$, and 
the \emph{prolongation} (or extension-by-zero) matrix of size
$n_\Gamma \times \widetilde n_\Gamma$ such that
$\widetilde{\bf u}_\Gamma=R_D^T{\bf u}_\Gamma$.

Similarly, for any couple $(k,\alpha)$, 
we define the \emph{restriction matrix} 
$R_{k,\alpha}$ of size $n^{(k,\alpha)}\times \widetilde n_\Gamma$
that implements the restriction of $\widetilde {\bf u}_\Gamma$ from $\Gamma$ 
to the degrees of freedom associated with $\gamma_{k,\alpha}$ 
(up to the boundary), i.e. 
${\bf u}_{k,\alpha}=R_{k,\alpha}\widetilde {\bf u}_\Gamma$. Consequently, 
$R_{k,\alpha}^T$ implements the prolongation of ${\bf u}_{k,\alpha}$ to
$\widetilde {\bf u}_\Gamma$.

\begin{figure}
\begin{center}
\includegraphics[width=0.4\textwidth]{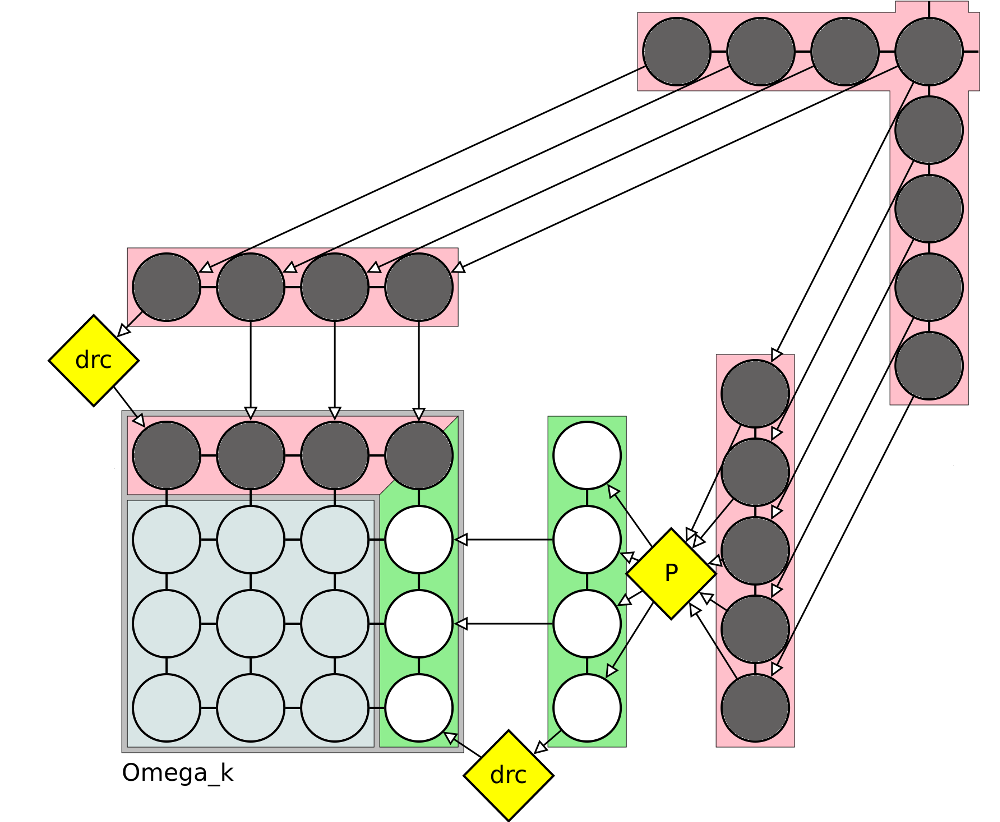}\quad
\includegraphics[width=0.4\textwidth]{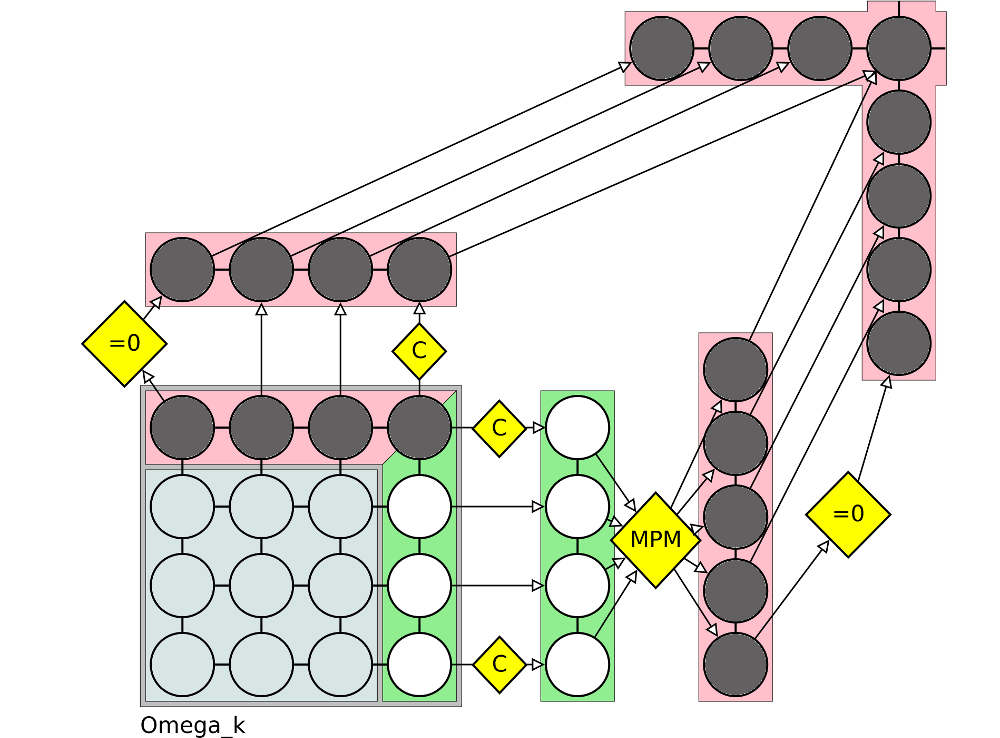}
\end{center}
\caption{At left, the interpolation of the trace from the master faces to
the slave faces, $P$ stands for the interpolation operator
$\Pi_{(k,\alpha)(\ell,\beta)}$, ``drc'' stands for Dirichlet
condition. At right,  the interpolation of the normal derivatives 
from slave faces  to master faces, $MPM$ stands for the interpolation operator 
$\widetilde \Pi_{(\ell,\beta)(k,\alpha)}$. ``C'' identifies the degrees of
freedom on which the
correction matrix $C^{(k)}$ defined in (\ref{matrix_corrM}) acts
}
\label{fig:local-global}
\end{figure}

In Fig. \ref{fig:ex1} and \ref{fig:local-global} we sketch a decomposition with 4
patches, we identify master (dark-gray) and slave (white)
 faces and we explain how the
interpolation operators act.

Algorithm \ref{algo:internodes1} contains the instructions to
initialize INTERNODES, Algorithm \ref{alg:b} 
computes the right hand side ${\bf b}$ of (\ref{SchurM}), while 
Algorithm \ref{alg:Slambda} implements the matrix vector product 
$\boldsymbol\psi=S\boldsymbol\lambda$ that can be used at each iteration 
of the iterative method called to solve (\ref{SchurM}). Once that ${\bf
u}_\Gamma$ has been computed, we can recover the solution ${\bf u}^{(k)}$ in
every patch by applying Algorithm \ref{alg:final}.

\begin{algorithm}[h!]
{\fontsize{10}{11}\selectfont{
\begin{algorithmic}
\ForAll{patch $k=1,\ldots,M$}
\State{build the local stiffness matrices $A^{(k)}$\;}
\State{build the arrays ${\bf
f}^{(k)}$, ${\bf g}_k$, ${\bf g}_{\partial\Gamma_k}$\;}
\State{build the matrices
$C^{(k)}$ (see  (\ref{matrix_corrM})) and
$\widehat{A}^{(\overline\gamma_{k,\alpha},\overline\Omega^{(k)})}$ (see
(\ref{AC})) \;}
\State{build the Greville nodes in $\Omega^{(k)}$\;}
\ForAll{face $\alpha$ of $\Gamma_k$}
\State{build the local interface mass matrices $M_{(k,\alpha)}$ (formula
(\ref{local_massM}))\;}
\ForAll{face $(\ell,\beta)\in{\cal A}_{(k,\alpha)}$}
\State{build the interpolation matrices $P_{(k,\alpha)(\ell,\beta)}$ (either
formulas
(\ref{P_kalb}) or (\ref{matrix_P21_RBF})--(\ref{matrix_P12_RBF}))\;}
\EndFor
\EndFor
\EndFor
\end{algorithmic}
\caption{Initialization of INTERNODES}
\label{algo:internodes1}
}}
\end{algorithm}

\begin{algorithm}[h!]
{\fontsize{10}{11}\selectfont{
\begin{algorithmic}
\State{\emph{\% Distribute the Dirichlet dof}}
\ForAll{patch $k=1,\ldots,M$}
\State{${\bf t}^{(k)}={\bf 0}$ (array of the dof associated with
 $\partial\Omega^{(k)}$)\;}
\State{${\bf t}^{(k)}|_{\partial\Omega_D^{(k)}}={\bf
g}|_{\partial\Omega_D^{(k)}}$\;}
\EndFor
\State{\emph{\% Interpolate from master faces to slave faces}}
\ForAll{patch $\ell=1,\ldots,M$}
\ForAll{slave face $\beta$ of $\Gamma_\ell$}
\State{$\displaystyle 
{\bf t}^{(\ell)}|_{(\ell,\beta)} = \sum_{(k,\alpha)\in{\cal A}_{(\ell,\beta)}}
P_{(\ell,\beta)(k,\alpha)}{\bf t}^{(k)}|_{(k,\alpha)} $\;}
\EndFor
\EndFor
\State{\emph{\% Solve local independent subproblems and compute ${\bf
r}_{k,\alpha}$ face by face}}
\ForAll{patch $k=1,\ldots,M$}
\State{solve $ A^{(k,k)}{\bf u}^{(k)}_0={\bf f}^{(k)}_0
-A^{(k,\partial\Omega^{(k)})}{\bf
t}^{(k)}$\;}
\State{assemble ${\bf u}^{(k)}=[{\bf u}^{(k)}_0,{\bf t}^{(k)}]$\;}
\ForAll{face $\alpha$ of $\Gamma_k$}
\State{${\bf r}_{k,\alpha}=
\widehat A^{(\overline\gamma_{k,\alpha},\overline\Omega^{(k)})}{\bf u}^{(k)}-{\bf
f}_{k,\alpha}$\;}
\EndFor
\EndFor
\State{\emph{\% Interpolate the derivatives from slave to master faces 
and assemble from local faces to global master skeleton}}
\State{$\widetilde{\bf b}={\bf 0}$ \;}
\ForAll{patch $k=1,\ldots,M$}
\ForAll{master face $\alpha$ of $\Gamma_k$}
\State{$\widetilde{\bf b}=\widetilde{\bf b}+R_{(k,\alpha)}^T\left(
\displaystyle {\bf r}_{k,\alpha}+\sum_{(\ell,\beta)\in{\cal A}_{(k,\alpha)}}
M_{(k,\alpha)}P_{(k,\alpha)(\ell,\beta)}M_{(\ell,\beta)}^{-1}{\bf
r}_{\ell,\beta}\right)$\;}
\EndFor
\EndFor
\State{\emph{\% Restrict $\widetilde{\bf b}$ to $\Gamma\setminus
\partial\Omega$}}
\State{${\bf b}=R_D\widetilde{\bf b}$\;}
\end{algorithmic}
\caption{Computation of the right hand side
${\bf b}$ of (\ref{SchurM})}
\label{alg:b}
}}
\end{algorithm}


\begin{algorithm}[h!]
{\fontsize{10}{11}\selectfont{
\begin{algorithmic}
\State{\emph{\% Expand $\boldsymbol\lambda$ from
$\Gamma\setminus\partial\Omega$ to $\Gamma$}\;}
\State{$\widetilde{\boldsymbol\lambda}=R_D^T\boldsymbol\lambda$ \;}
\State{\emph{\% Distribute the trace from the 
global master skeleton $\Gamma$ to the local master faces}}
\ForAll{patch $k=1,\ldots,M$}
\State{${\bf t}^{(k)}={\bf 0}$ (array of the dof associated with
 $\partial\Omega^{(k)}$)\;}
\ForAll{master face $\alpha$ of $\Gamma_k$}
\State{${\bf t}^{(k)}|_{(k,\alpha)} =
R_{(k,\alpha)}\widetilde{\boldsymbol\lambda}$\;}
\EndFor
\EndFor
\State{\emph{\% Interpolate from master faces to slave faces}}
\ForAll{patch $\ell=1,\ldots,M$}
\ForAll{slave face $\beta$ of $\Gamma_\ell$}
\State{$\displaystyle 
{\bf t}^{(\ell)}|_{(\ell,\beta)} = \sum_{(k,\alpha)\in{\cal A}_{(\ell,\beta)}}
P_{(\ell,\beta)(k,\alpha)}{\bf t}^{(k)}|_{(k,\alpha)} $\;}
\EndFor
\EndFor
\State{\emph{\% Solve local independent subproblems and compute ${\bf
r}_{k,\alpha}$ face by face}}
\ForAll{patch $k=1,\ldots,M$}
\State{solve $ A^{(k,k)}{\bf u}^{(k)}_0=-A^{(k,\partial\Omega^{(k)})}{\bf
t}^{(k)}$\;}
\State{assemble ${\bf u}^{(k)}=[{\bf u}^{(k)}_0,{\bf t}^{(k)}]$\;}
\ForAll{face $\alpha$ of $\Gamma_k$}
\State{compute ${\bf r}_{k,\alpha}=
\widehat A^{(\overline\gamma_{k,\alpha},\overline\Omega^{(k)})}{\bf u}^{(k)}$\;}
\EndFor
\EndFor
\State{$\widetilde{\boldsymbol\psi}={\bf 0}$ (same size as
$\widetilde{\boldsymbol\lambda}$)\;}
\State{\emph{\% Interpolate the derivatives from slave to master faces
and assemble from local faces to global master skeleton}}
\ForAll{patch $k=1,\ldots,M$}
\ForAll{master face $\alpha$ of $\Gamma_k$}
\State{$\widetilde{\boldsymbol\psi}=\widetilde{\boldsymbol\psi}+R_{(k,\alpha)}^T\left(
\displaystyle {\bf r}_{k,\alpha}+\sum_{(\ell,\beta)\in{\cal A}_{(k,\alpha)}}
M_{(k,\alpha)}P_{(k,\alpha)(\ell,\beta)}M_{(\ell,\beta)}^{-1}{\bf
r}_{\ell,\beta}\right)$\;}
\EndFor
\EndFor
\State{\emph{\% Restrict $\boldsymbol\psi$ to $\Gamma\setminus
\partial\Omega$}}
\State{$\boldsymbol\psi=R_D\widetilde{\boldsymbol\psi}$\;}
\end{algorithmic}
\caption{Given $\boldsymbol\lambda$, computation of
$\boldsymbol\psi=S\boldsymbol\lambda$. This is 
the matrix-vector product needed to solve
(\ref{SchurM}) by Krylov methods}
\label{alg:Slambda}
}}
\end{algorithm}

\begin{algorithm}[h!]
{\fontsize{10}{11}\selectfont{
\begin{algorithmic}
\State{\emph{\% Expand ${\bf u}_\Gamma$ from
$\Gamma\setminus\partial\Omega$ to $\Gamma$}\;}
\State{$\widetilde{\bf u}_\Gamma=R_D^T{\bf u}_\Gamma$ \;}
\State{\emph{\% Distribute the Dirichlet dof and the trace from $\Gamma$ to 
the master faces $\gamma_{k,\alpha}$}}
\ForAll{patch $k=1,\ldots,M$}
\State{${\bf t}^{(k)}={\bf 0}$ (array of the dof associated with
 $\partial\Omega^{(k)}$)\;}
\ForAll{master face $\alpha$ of $\Gamma_k$}
\State{${\bf t}^{(k)}|_{(k,\alpha)} =
R_{(k,\alpha)}\widetilde{\bf u}_\Gamma$\;}
\EndFor
\State{${\bf t}^{(k)}|_{\partial\Omega_D^{(k)}}={\bf
g}|_{\partial\Omega_D^{(k)}}$\;}
\EndFor
\State{\emph{\% Interpolate from master faces to slave faces}}
\ForAll{patch $\ell=1,\ldots,M$}
\ForAll{slave face $\beta$ of $\Gamma_\ell$}
\State{$\displaystyle 
{\bf t}^{(\ell)}|_{(\ell,\beta)} = \sum_{(k,\alpha)\in{\cal A}_{(\ell,\beta)}}
P_{(\ell,\beta)(k,\alpha)}{\bf t}^{(k)}|_{(k,\alpha)} $\;}
\EndFor
\EndFor
\State{\emph{\% Solve local independent subproblems}}
\ForAll{patch $k=1,\ldots,M$}
\State{solve $ A^{(k,k)}{\bf u}^{(k)}_0={\bf f}^{(k)}_0
-A^{(k,\partial\Omega^{(k)})}{\bf
t}^{(k)}$\;}
\State{assemble ${\bf u}^{(k)}=[{\bf u}^{(k)}_0,{\bf t}^{(k)}]$\;}
\EndFor
\end{algorithmic}
\caption{Given ${\bf u}_\Gamma$, computation of
 the local solutions ${\bf u}^{(k)}$, for $k=1,\ldots,M$}
\label{alg:final}
}}
\end{algorithm}

\section{Numerical results for $M>2$ patches}\label{sec:numres_mdom}


\subsection{Test case \#3}\label{sec:numres_mdom_2d}

\begin{figure}[h!]
\begin{center}
\includegraphics[width=0.32\textwidth]{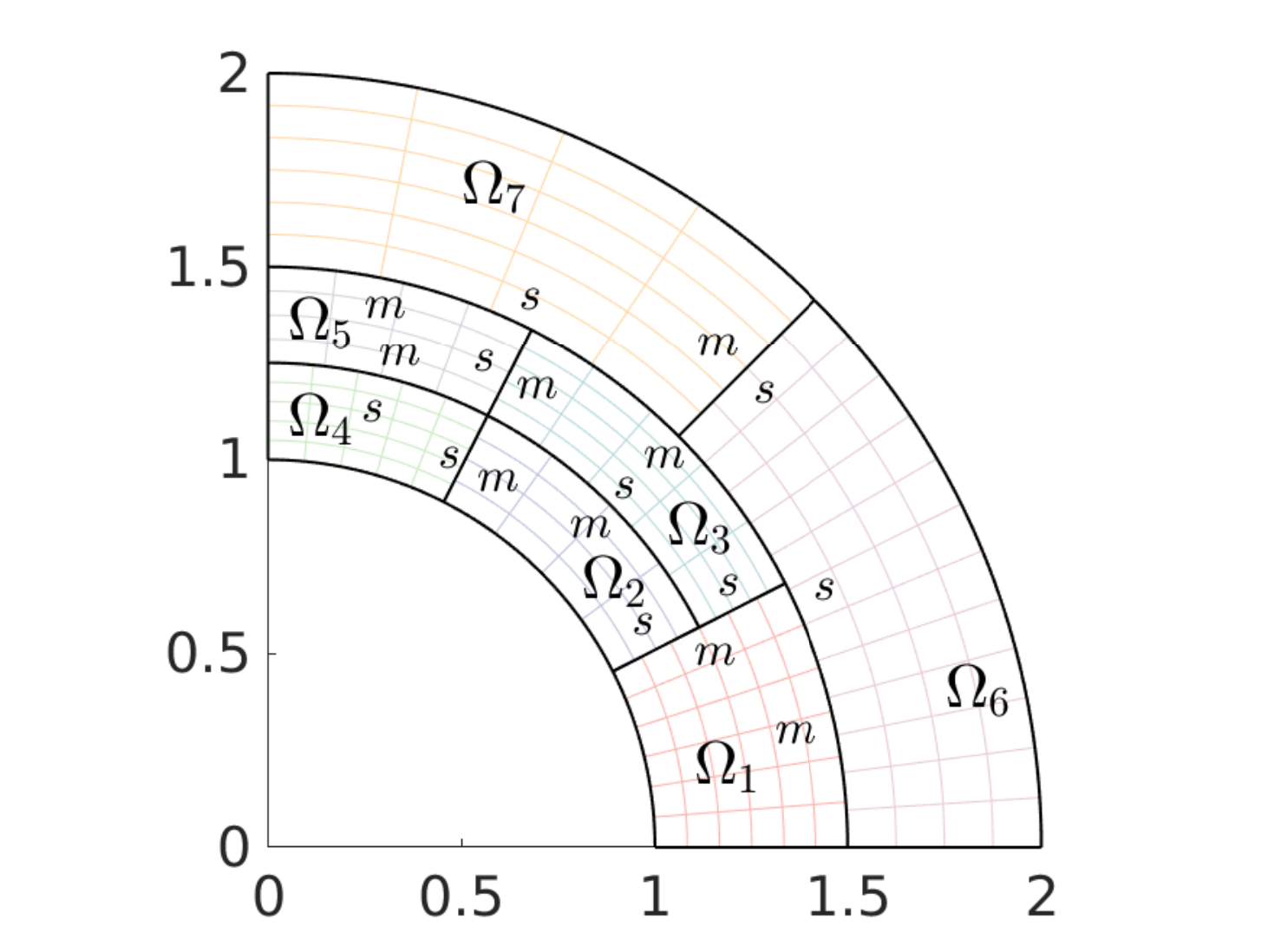}\quad
\includegraphics[width=0.3\textwidth]{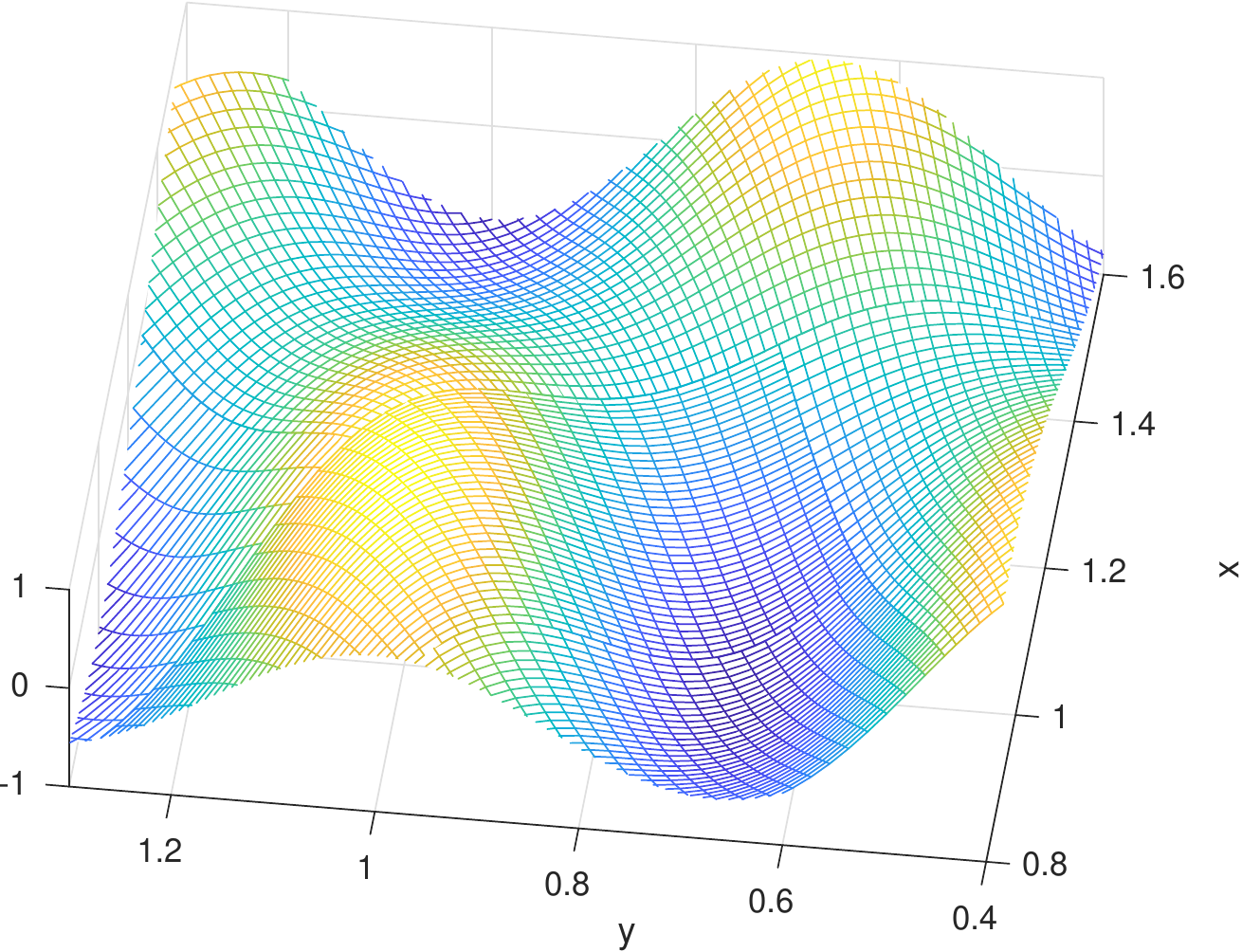}\quad
\includegraphics[width=0.3\textwidth]{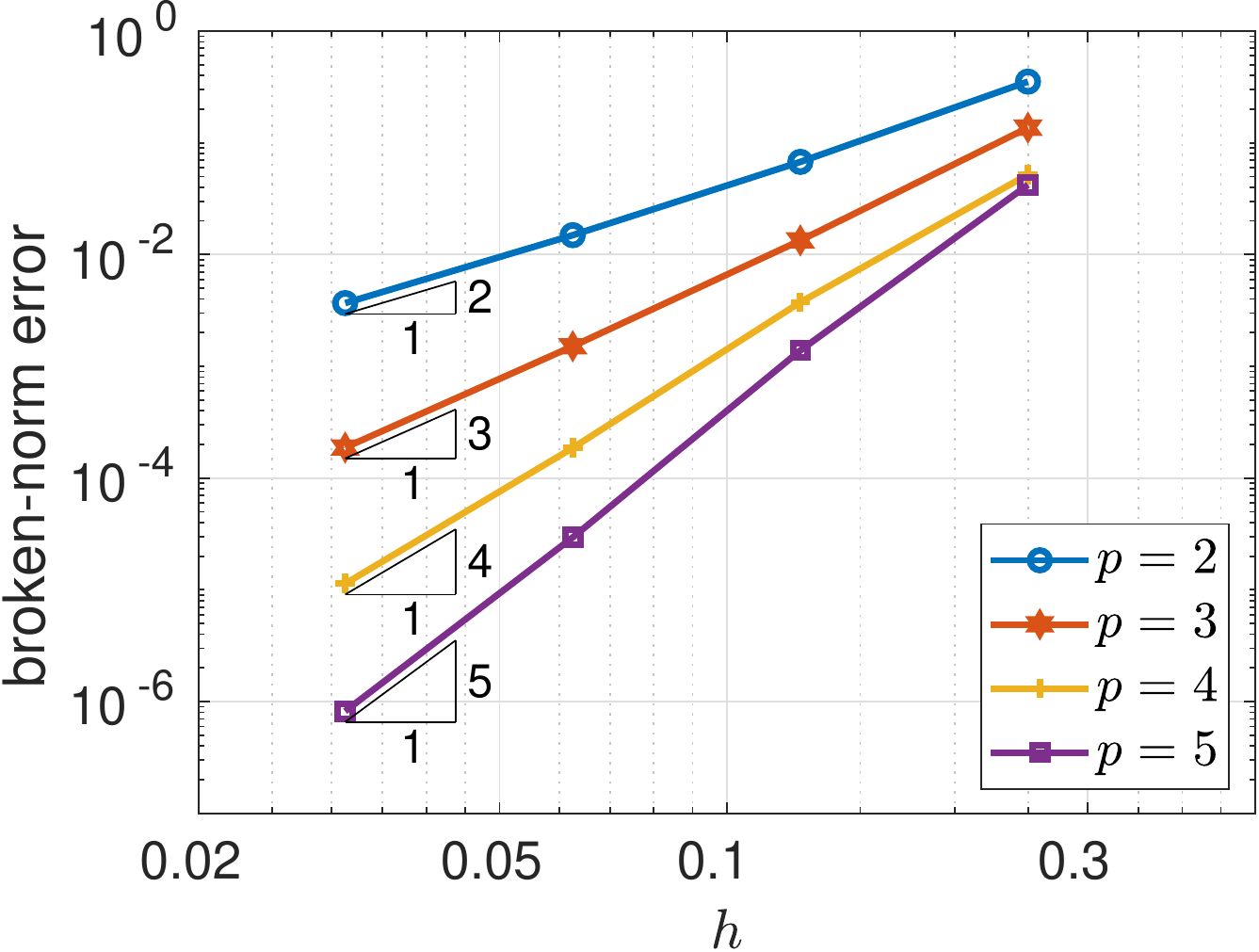}
\end{center}
\caption{\emph{Test case \#3}.
At left, the multipatch configuration; in the middle, a zoom on 
the numerical solution computed with $p=4$ and $\overline n=8$;
at right, the broken-norm error versus the mesh size $h=\max_k h_k$}
\label{fig:7dom}
\end{figure}

In this test case we consider watertight interfaces, but with non-matching
parametrizations.
Let us consider again the differential problem (\ref{problem}) in 
$\Omega=\{(x,y)\in
{\mathbb R}^2:\ x\geq0, \ y\geq 0, \ 1\leq x^2+y^2\leq 4\}$
 with  $\alpha=0$, and
$f$ and  $g$ such that the exact solution is
$u(x,y)=\sin(1.5\pi x)\sin(3\pi y)$. Now we split $\Omega$ in 7 patches and we
tag the master/slave interfaces as
shown in the left picture of Fig. \ref{fig:7dom}. We fix the polynomial degree
equal to $p=2,\ldots,5$ in all the patches, while the number of elements inside
the patches is defined as in the following table, with $\overline n=4,8,16,32$:
\begin{center}\begin{tabular}{cc}
patch & number of elements\\
\hline
$\Omega^{(1)}$ & $(\overline{n}+2)\times \overline{n}$\\
$\Omega^{(2)}$ and $\Omega^{(5)}$ & $\overline{n}\times \overline{n}$\\
$\Omega^{(3)}$ and $\Omega^{(4)}$ & $(\overline{n}+1)\times (\overline{n}+1)$\\
$\Omega^{(6)}$ & $\overline{n}\times 3\overline{n}$\\
$\Omega^{(7)}$ & $(\overline{n}+2)\times \overline{n}$\\
\end{tabular}
\end{center}


As in the Test case \#1, the first parametric coordinate is mapped onto
the radial coordinate of the physical domain, and the second coordinate
onto the angular one. The patches are sectors of rings, the control points are
defined as described in \cite[Sect. 2.4.1.1]{chb_iga_book}, in particular the
patches $\Omega^{(k)}$ with $k=1,\ldots,5$ have $30^\circ$ degrees arcs, while 
$\Omega^{(6)}$ and $\Omega^{(7)}$ have $45^\circ$ degrees arcs. Once the
control points and the weights have been computed, the NURBS
parametrization of the patches is defined as in (\ref{map}).
 In the right picture of Fig. \ref{fig:7dom} the broken-norm
errors (\ref{err_broken}) are shown versus $h=\max_k h_k$ for $p=2,\ldots,5$.
As in the case of two subdomains, INTERNODES exhibits optimal convergence order
with respect to the mesh size $h$.

\subsection{Test case \#4. Jumping coefficients, the Kellogg's test case}

We solve the elliptic problem $-\nabla\cdot(\nu\nabla u)=0$ in
$\Omega=(-1,1)^2$
with Dirichlet boundary conditions on $\partial\Omega$ 
and piece-wise constant coefficient $\nu$ such that
the exact solution is 
the so-called Kellogg's function (see, e.g.,
\cite{morin_nochetto_siebert,gq_camc,gq_dd24}).
This is a very challenging problem whose solution features low regularity.
We refer to \cite{Langer_meshgrading} for a more in-depth analysis of the
problem in the framework of isogeometric analysis.

The Kellogg's solution can be written
in terms of the polar coordinates $r$ and $\theta$ as
$u(r,\theta)=r^\gamma \mu(\theta)$,
where  $\gamma\in(0,2)$ is a given parameter, while
$\mu(\theta)$ is a $2\pi-$periodic continuous function defined like follows:
\begin{eqnarray}\label{Kellogg_mu}
\mu(\theta)=\left\{\begin{array}{ll}
\cos((\pi/2-\sigma)\gamma)\cos((\theta-\pi/2+\rho)\gamma) & 0\leq \theta \leq
\pi/2\\
\cos(\rho\gamma)\cos((\theta-\pi+\sigma)\gamma) & \pi/2\leq \theta\leq \pi\\
\cos(\sigma\gamma)\cos((\theta-\pi-\rho)\gamma) & \pi\leq \theta\leq 3\pi/2\\
\cos((\pi/2-\rho)\gamma)\cos((\theta-3\pi/2-\sigma)\gamma) &
3\pi/2\leq\theta\leq 2\pi.
\end{array}
\right.
\end{eqnarray}
The parameters $\sigma$, $\rho$, $\gamma$ and the coefficient $R$ (that
is involved in the definition of $\alpha$)
 must satisfy the following
non-linear system:
\begin{eqnarray}\label{kellogg_conditions}
\left\{\begin{array}{l}
R=-\tan((\pi/2-\sigma)\gamma)\cot(\rho\gamma)\\
\frac{1}{R}=-\tan(\rho\gamma)\cot(\sigma\gamma)\\
R=-\tan(\sigma\gamma)\cot((\pi/2-\rho)\gamma)\\
0<\gamma<2\\
\max\{0,\pi\gamma-\pi\}<2\gamma\rho<\min\{\gamma\pi,\pi\}\\
\max\{0,\pi-\gamma\pi\}<-2\gamma\sigma<\min\{\pi,2\pi-\gamma\pi\}.
\end{array}\right.
\end{eqnarray}

We set $\nu=R>0$ in
the first and the third quadrants, and $\nu=1$ in the second and in the
fourth ones.  

The case $\gamma=1$ is trivial since the solution is a plane.
When $\gamma\neq 1$, then $u\in H^{1+\gamma-\varepsilon}(\Omega)$ for any
$\varepsilon>0$, in particular
the solution features low regularity at the origin
and on the axes.

We look for the approximation of the Kellogg's solution
 by applying INTERNODES to the 4-subdomains
decomposition induced by the discontinuity of $\nu$.

For $\bar n\in\{5,10,\ldots,30\}$ we
 consider $(2\bar n+1)\times (2\bar n+1) $ equal elements in both 
$\Omega^{(1)}$ and $\Omega^{(3)}$ and
$(\bar n -1)\times(\bar n -1)$ equal elements
in both $\Omega^{(2)}$ and $\Omega^{(4)}$, 
so that the parametrizations
(defined as in (\ref{map}) do not
match on any couple of interfaces.
To analyze the errors we take the same
 polynomial degree $p$ along each direction and in all patches.

\begin{figure}
\begin{center}
\includegraphics[trim={0 20 0
0},width=0.25\textwidth]{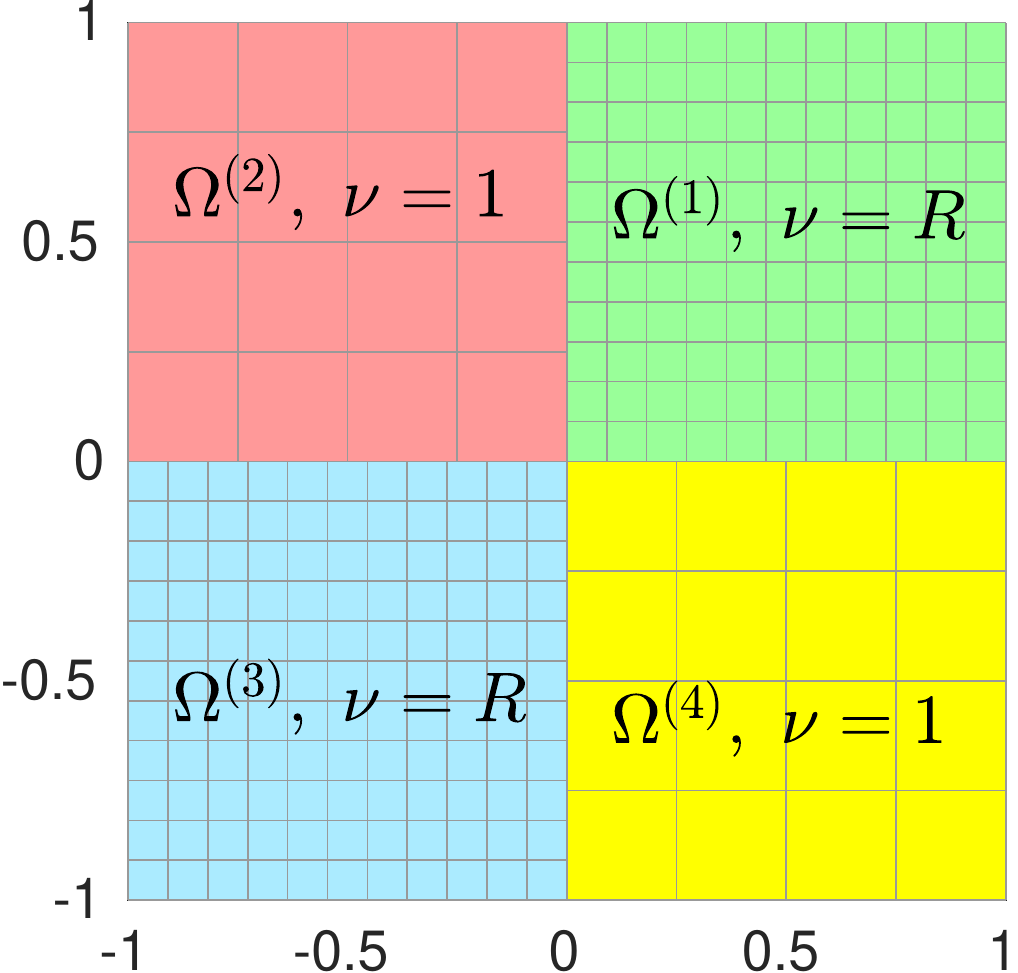}\quad
\includegraphics[trim={0 20 0
0},width=0.25\textwidth]{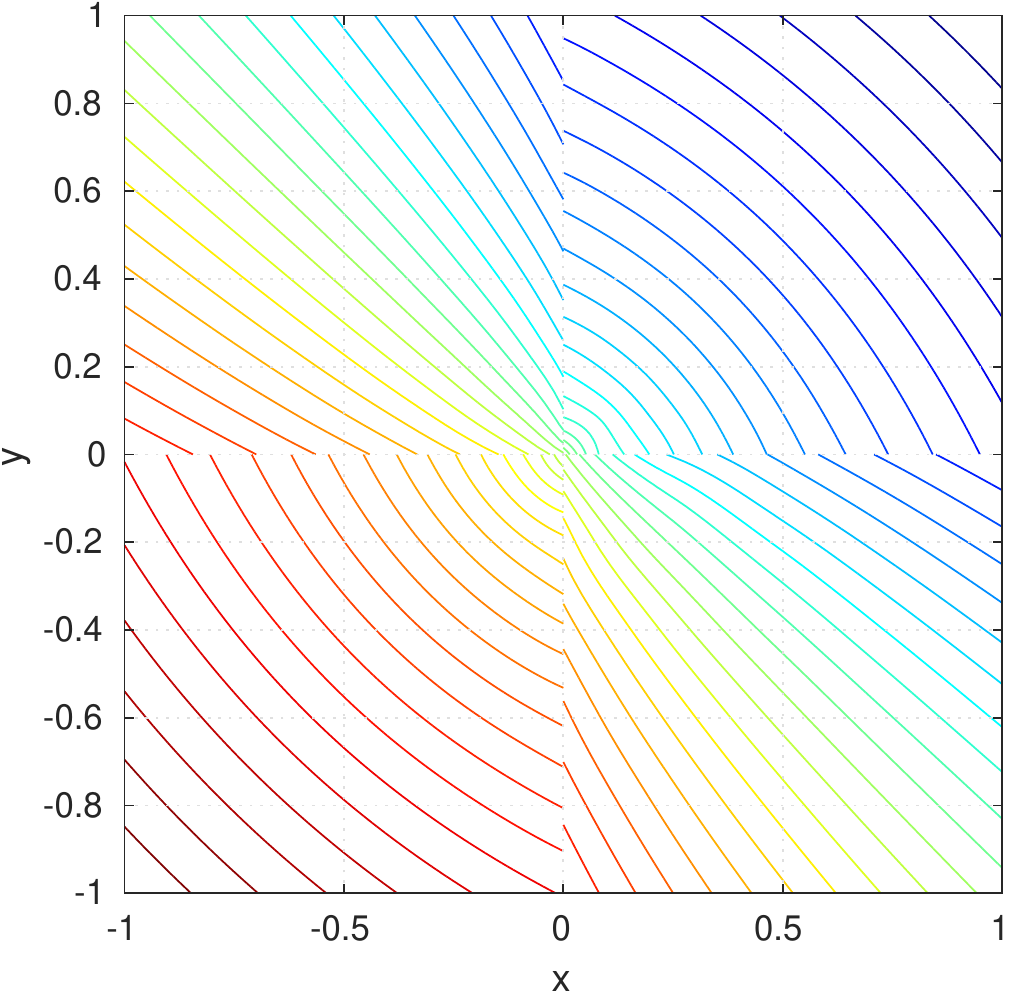}\quad
\includegraphics[trim={0 20 0
0},width=0.4\textwidth]{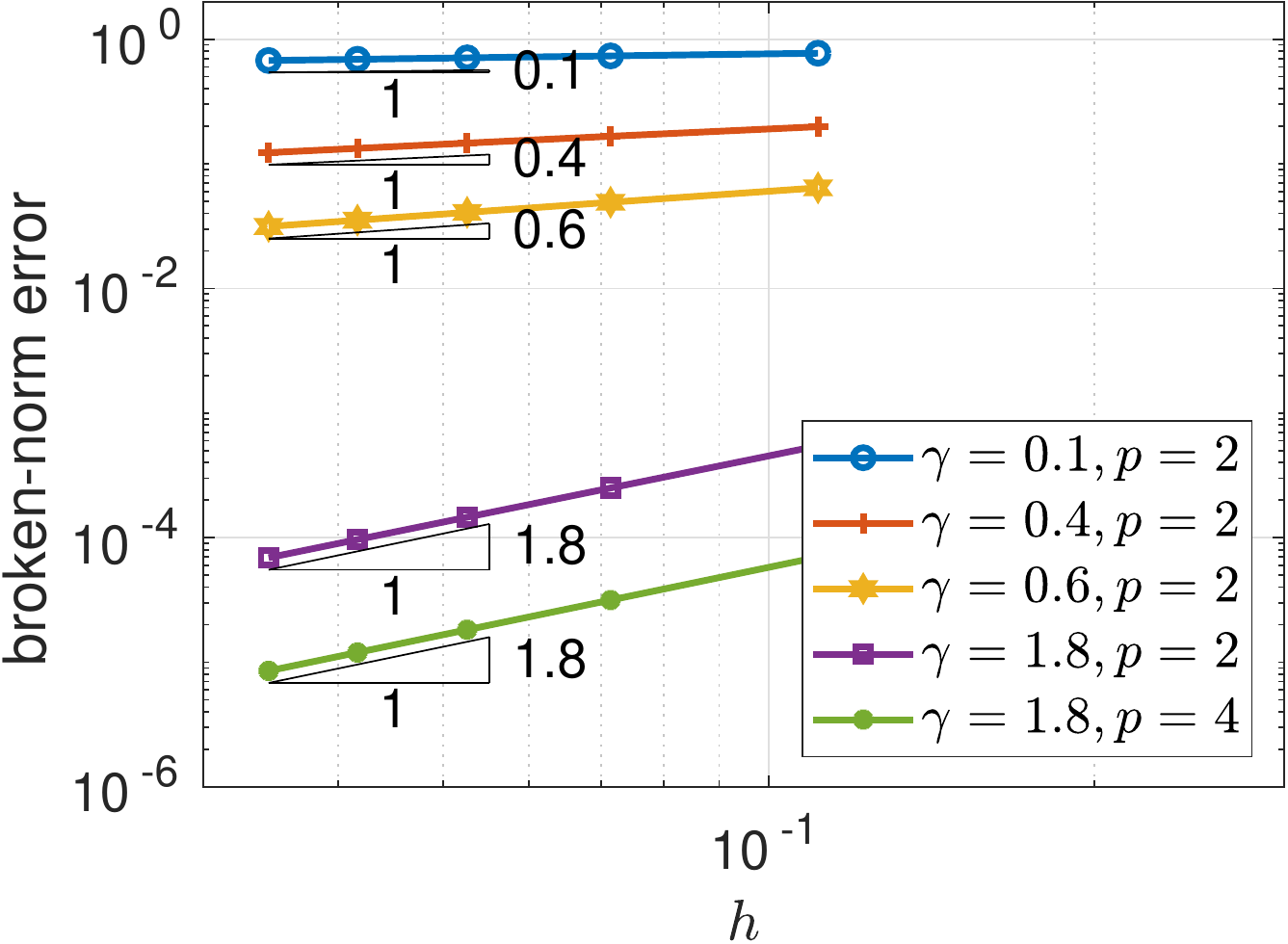}
\end{center}
\caption{\emph{Test case \#4. Kellogg's solution}.
At left, the geometry for the Kellogg's test case. At middle, the
numerical solution with $\gamma=0.6$. At right the errors versus the
discretization parameter $h_1=1/n_{el}^{(1)}$}
\label{fig:kellogg}
\end{figure}

In the left picture of 
Fig. \ref{fig:kellogg} the multipatch configuration is shown;
 in the middle picture of the same Figure the numerical
solution corresponding to $\gamma=0.6$ is plotted, it is computed by setting
the polynomial degree $p=2$ in each patch, $11\times 11$ equal elements 
in $\Omega^{(1)}$ and $\Omega^{(3)}$ and
$4\times 4$ equal elements
in both $\Omega^{(2)}$ and $\Omega^{(4)}$.

We have considered four different values of $\gamma$:
\begin{enumerate}[noitemsep]
\item $\gamma=0.1$ and
$R\simeq 161.45$,  so that the corresponding Kellogg's solution $u$ belongs to
$H^{1.1-\varepsilon}(\Omega)$ (for any $\varepsilon>0$),
\item $\gamma=0.4$ and
$R\simeq 9.47$,  they provide
$u\in H^{1.4-\varepsilon}(\Omega)$,
\item $\gamma=0.6$ and
$R\simeq 3.85$, they provide
$u\in H^{1.6-\varepsilon}(\Omega)$,
\item
$\gamma=1.8$ and
 $R\simeq 2.5\cdot 10^{-2}$, 
they provide  $u\in H^{2.8-\varepsilon}(\Omega)$.
\end{enumerate}

The broken-norm errors  (\ref{err_broken})
versus the mesh size $h=\max_k h_k$ are shown in the right picture of Fig.
\ref{fig:kellogg}. 
They behave like 
$h^{\min(s-1,p)}$ when $h\to 0$, where
$s=1+\gamma-\varepsilon$ is the Sobolev regularity of the Kellogg's
solution. We conclude that
the INTERNODES solution is converging to the exact
one when $h\to 0$ with the best possible convergence rate dictated by the
regularity of the Kellogg's solution.

We briefly comment on the algebraic solution of  the Schur complement
system (\ref{schur}). The number of Bi-CGStab iterations needed to solve 
(\ref{schur}) depends on the regularity of the solution (dictated by the
parameter $\gamma$), on the mesh-sizes $h_k$, on the polynomial degrees
$p$, and on the choice of the master/slave edges.
In Table \ref{tab:test4_iter}, left,
 we show the number \#it of Bi-CGStab iterations
needed to converge up a tolerance $\epsilon=10^{-10}$ for three different
configurations of master/slave interfaces, for $\gamma\in\{0.6, \ 1.8\}$,
 and for some values of $h_k$ and $p$. The master/slave configuration that
performs better is that featuring all the master edges inside the patches with
the coarser meshes (i.e. $\Omega^{(2)}$ and $\Omega^{(4)}$). The dependence of
\#it on $p$ reflects the way the condition numbers of stiffness IGA matrices
depend on $p$.

In the case that we can distinguish between
master and slave patches (we say that a patch is master (slave resp.)
if all its edges are tagged as master (slave, resp.)), 
to reduce the number of iterations we implement a
 preconditioner $P$ as follows (it is a generalization of the Dirichlet/Neumann
preconditioner for substructuring domain decomposition methods \cite{qv_ddm}).
We denote  by $R_{\Gamma_k}$ the restriction matrix
from the global master skeleton $\Gamma$ to the internal
boundary $\Gamma_k$ of $\Omega^{(k)}$ and by $U$ a diagonal matrix
whose element $U_{ii}$ is the number of master patches which the $i$th 
 degree of freedom
of $\Gamma$ belongs to. Then we define $P$ such that:
\begin{equation}\label{preco}
P^{-1}=U^{-1}\sum_{\tiny\begin{array}{c}k\\ master\end{array}}
R^T_{\Gamma_k}S_{\Gamma_k}^{-1}R_{\Gamma_k}.
\end{equation}
The matrices $S_{\Gamma_k}^{-1}$ are never built explicitly; to compute
$S_{\Gamma_k}^{-1}{\bf v}$ we must 
solve a differential problem in $\Omega^{(k)}$ of the same
nature of (\ref{problem}), but with Neumann (instead of Dirichlet) 
data on $\Gamma_k$ (see \cite{qv_ddm}).

The  PBi-CGStab iterations with $P$ defined as in (\ref{preco}) are
shown in Table \ref{tab:test4_iter}, right, they are independent of the
discretization parameters $h$ and $p$, but they depend on both the master/slave
configuration and the regularity of the Kellogg's solution.

A more in-depth analysis of suitable preconditioners for system (\ref{schur})
will be subject of future work.

\begin{table}
\begin{center}
\begin{tabular}[t]{c|ccc|ccc|ccc}
 & \multicolumn{3}{c|}{$\gamma=0.6$} 
 & \multicolumn{3}{c|}{$\gamma=1.8$} 
 & \multicolumn{3}{c}{$\gamma=1.8$} \\
 & \multicolumn{3}{c|}{$p=2$} 
 & \multicolumn{3}{c|}{$p=2$} 
 & \multicolumn{3}{c}{$p=4$} \\
$\overline n$ & (a) & (b) & (c) & (a) & (b) & (c) &(a) & (b) & (c) \\
\hline
10 &  17 & 12 & 16 & 43 & 13 & 42 &  86 & 23 & 79 \\
15 &  21 & 14 & 19 & 47 & 15 & 44 & 103 & 25 & 97 \\
20 &  25 & 17 & 22 & 56 & 17 & 54 & 118 & 23 & 99 \\
25 &  30 & 18 & 25 & 62 & 20 & 51 & 115 & 25 & 105 \\
30 &  33 & 22 & 28 & 61 & 22 & 52 & 110 & 27 &124 \\
\end{tabular}
\qquad \qquad \begin{tabular}[t]{c|cc|cc|cc}
 & \multicolumn{2}{c|}{$\gamma=0.6$}
 & \multicolumn{2}{c|}{$\gamma=1.8$}
 & \multicolumn{2}{c}{$\gamma=1.8$} \\
 & \multicolumn{2}{c|}{$p=2$}
 & \multicolumn{2}{c|}{$p=2$}
 & \multicolumn{2}{c}{$p=4$} \\
$\overline n$ & (a) & (b) & (a) & (b) &(a) & (b) \\
\hline
10 &  7 & 10 & 21 & 5 & 21 & 5 \\
15 &  8 & 11 & 21 & 5 & 24 & 5 \\
20 &  7 & 11 & 23 & 5 & 22 & 5 \\
25 &  9 & 11 & 23 & 5 & 23 & 5 \\
30 &  9 & 11 & 21 & 5 & 25 & 5 \\
\end{tabular}
\end{center}
\caption{\emph{Test case \# 4.} At left,
the number of Bi-CGStab iterations needed to solve
(\ref{schur}). (a) $\Omega^{(1)}$ and $\Omega^{(3)}$ master patches;
(b) $\Omega^{(2)}$ and $\Omega^{(4)}$ master patches;
(c) one master edge and one slave edge for each patch. At right, 
number of PBi-CGStab iterations needed to solve
(\ref{schur}). The preconditioner $P$ is defined in (\ref{preco})
}
\label{tab:test4_iter}
\end{table}

Finally, in Table \ref{tab:kellog_iter} the number of PBi-CGStab iterations required to
solve the Schur complement system (\ref{schur}) is shown for 
all the considered values of $\gamma$, here
 we have tagged as master the interfaces of the patches $\Omega^{(2)}$
and $\Omega^{(4)}$ and slave the others.

\begin{table}
\begin{center}
\begin{tabular}{ccc|ccccc}
$\overline n$ &\multicolumn{2}{c|}{ dof$(\Omega)$} &$\gamma=0.1$ & $\gamma=0.4$ & $\gamma=0.6$ & \multicolumn{2}{c}{$\gamma=1.8$} \\
&$p=2$ & $p=4$& $p=2$ & $p=2$& $p=2$ & $p=2$ & $p=4$\\
\hline
$10$ &1300 & 1588 & 11 & 10 & 10 & 5 & 5 \\
$15$ &2690 & 3098 & 11 & 11 & 11 & 5 & 5 \\
$20$ &4580 & 5108 & 12 & 12 & 11 & 5 & 5 \\
$25$ &6970 & 7618 & 12 & 11 & 11 & 5 & 5 \\
$30$ &9860 & 10628 & 12 & 11 & 11 & 5 & 5 \\
\end{tabular}
\end{center}
\caption{\emph{Test case \#4.}
Comparison of the number of PBi-CGStab iterations required to
solve the Schur complement system (\ref{schur}) for all the values of $\gamma$
considered in Fig. \ref{fig:kellogg}.  $\Omega^{(2)}$
and $\Omega^{(4)}$ master patches}
\label{tab:kellog_iter}
\end{table}

\begin{remark}
The  preconditioner defined in (\ref{preco}) can be extended to
more general decompositions provided that
all the interfaces of a single patch can be tagged either slave or master. 
Moreover, in the case that a patch has empty intersection with
$\partial\Omega$ and the corresponding local Schur complement is singular,
the strategy to add the mass matrix to the stiffness one can be adopted
 (see, e.g., \cite{qv_ddm}).
The design of suitable preconditioners for more general 
configurations with both master and slave edges in a single patch
requires further work.
\end{remark}

\subsection{Test case \#5. Non-watertight interfaces}

\begin{figure}
\begin{center}
\includegraphics[width=0.26\textwidth]{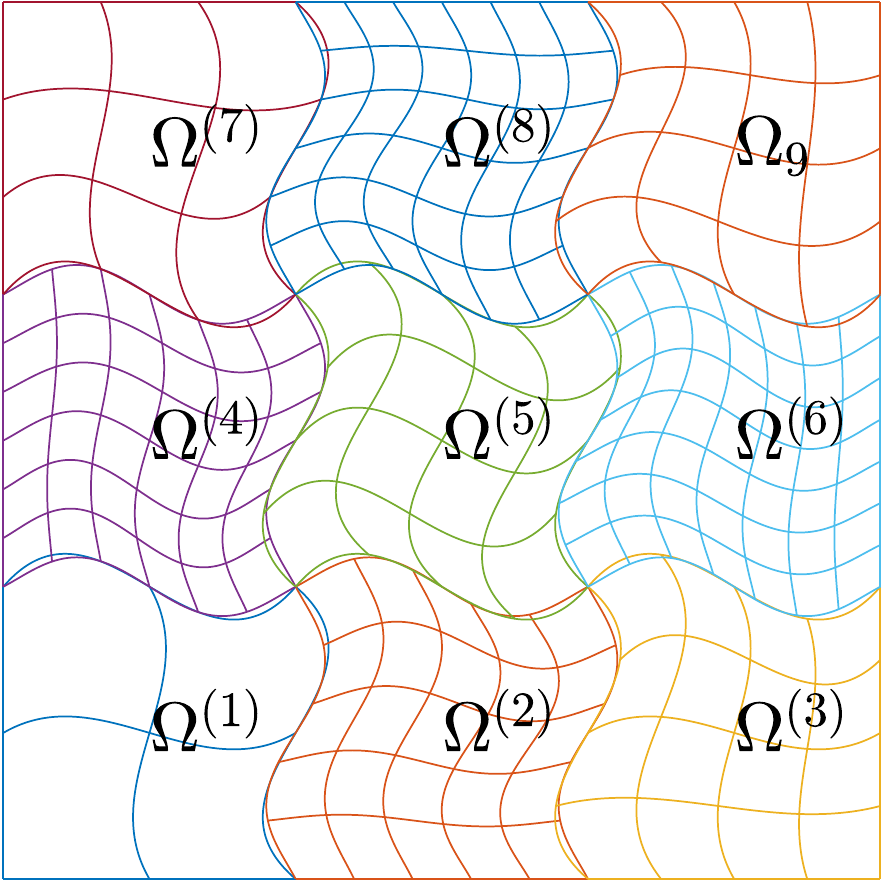}\qquad 
\includegraphics[trim={0 20 0 0},width=0.4\textwidth]{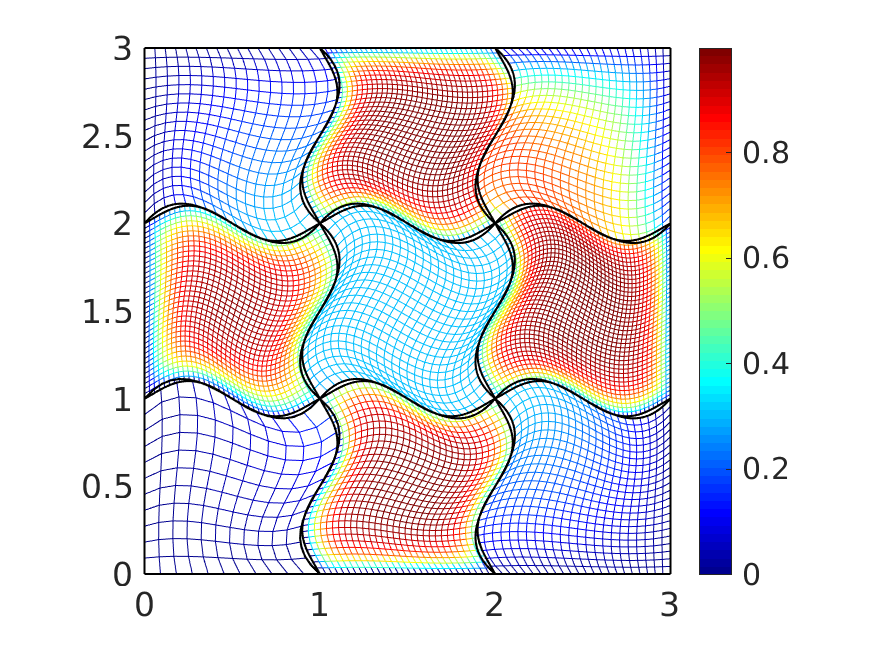}
\end{center}
\caption{\emph{Test case \#5.}
At left, the geometry; at right the 
 numerical solution computed by INTERNODES and RL-RBF.
The number of grid points used in the right plot is about
twice the number of degrees of freedom}
\label{fig:test5}
\end{figure}

We solve the elliptic problem $-\nabla\cdot(\nu\nabla u)+u=0$ in
$\Omega=(0,3)^2$
with homogeneous Dirichlet boundary conditions on $\partial\Omega$
 and piece-wise constant coefficient $\nu$.
The computational domain is split into nine patches as shown in Fig.
\ref{fig:test5}. The interfaces are non-watertight and they 
are obtained by B-spline interpolation  of sinusoidal curves,
 as already described  in Sect. \ref{sec:test2}. The maximum size of gaps and
overlaps has been set $d_\Gamma\simeq 0.0276$.

Then we define a piecewise function $\nu$ that assumes the values $\{
10,0.005,1,0.01,100,0.005,1,.005,0.1\}$ inside the patches 
(ordered as in the left picture of Fig. \ref{fig:test5}) (when two patches
overlap we consider two different values of $\nu$, depending on the patch
(see, e.g., \cite{hofer_langer_2019})).

We have considered uniform knot sets in all the patches, the number of elements
can be read from  the left picture of Fig. \ref{fig:test5}.
The numerical solution, computed with INTERNODES and RL-RBF interpolation
is shown in Fig. \ref{fig:test5}, right.

The interfaces of the odd patches have been tagged as master, those of the even
patches as slave and we have solved the Schur complement system (\ref{schur})
by the preconditioned Bi-CGStab iterations, with preconditioner $P$ build as
mentioned for the Test case \#4.
25 iterations have been required by 
the PBi-CGStab method to converge up to a tolerance of $\epsilon=10^{-10}$.

\subsection{Test case \#6. $C^\infty-$ solution on a 3D geometry}

Let us consider the differential problem (\ref{problem}) with
$\alpha=1$ in the domain $\Omega=\{(x,y,z)\in{\mathbb R}^3:\ 
0.25\leq (x^2+y^2)\leq 2.25,\ 0\leq z\leq 1\}$.
The functions $f$ and $g$ are such that the
exact solution is $u(x,y,z)=\sin(\pi x)\sin(\pi y)\cos(2\pi z)$.

The domain $\Omega$ is split into four patches like in
Fig. \ref{fig:4dom_3d} featuring non-matching NURBS 
parametrizations at the
interfaces, even if they are watertight. More precisely, in cylindrical
coordinates the patches are
$ \Omega^{(1)}=\{(r,\theta,z):\ 0.5\leq r\leq 1,\ 0\leq \theta \leq \pi/2,\ 
0\leq z\leq 1\}$,
$ \Omega^{(2)}=\{(r,\theta,z):\ 1\leq r\leq 1.5,\ 0\leq \theta \leq \pi/2,\ 
0\leq z\leq 0.5\}$,
$ \Omega^{(3)}=\{(r,\theta,z):\ 1\leq r\leq 1.5,\ 0\leq \theta \leq \pi/4,\ 
0.5\leq z\leq 1\}$,
$ \Omega^{(4)}=\{(r,\theta,z):\ 1\leq r\leq 1.5,\ \pi/4\leq \theta \leq \pi/2,\ 
0.5\leq z\leq 1\}$. 
The control
points and the weights  of the sector of rings are defined accordingly
to \cite[Sect.
2.4]{chb_iga_book}, then the $xy-$surfaces have been 
extruded along the $z$ direction.

\begin{figure}
\begin{center}
\includegraphics[trim={0 20 0
0},width=0.4\textwidth]{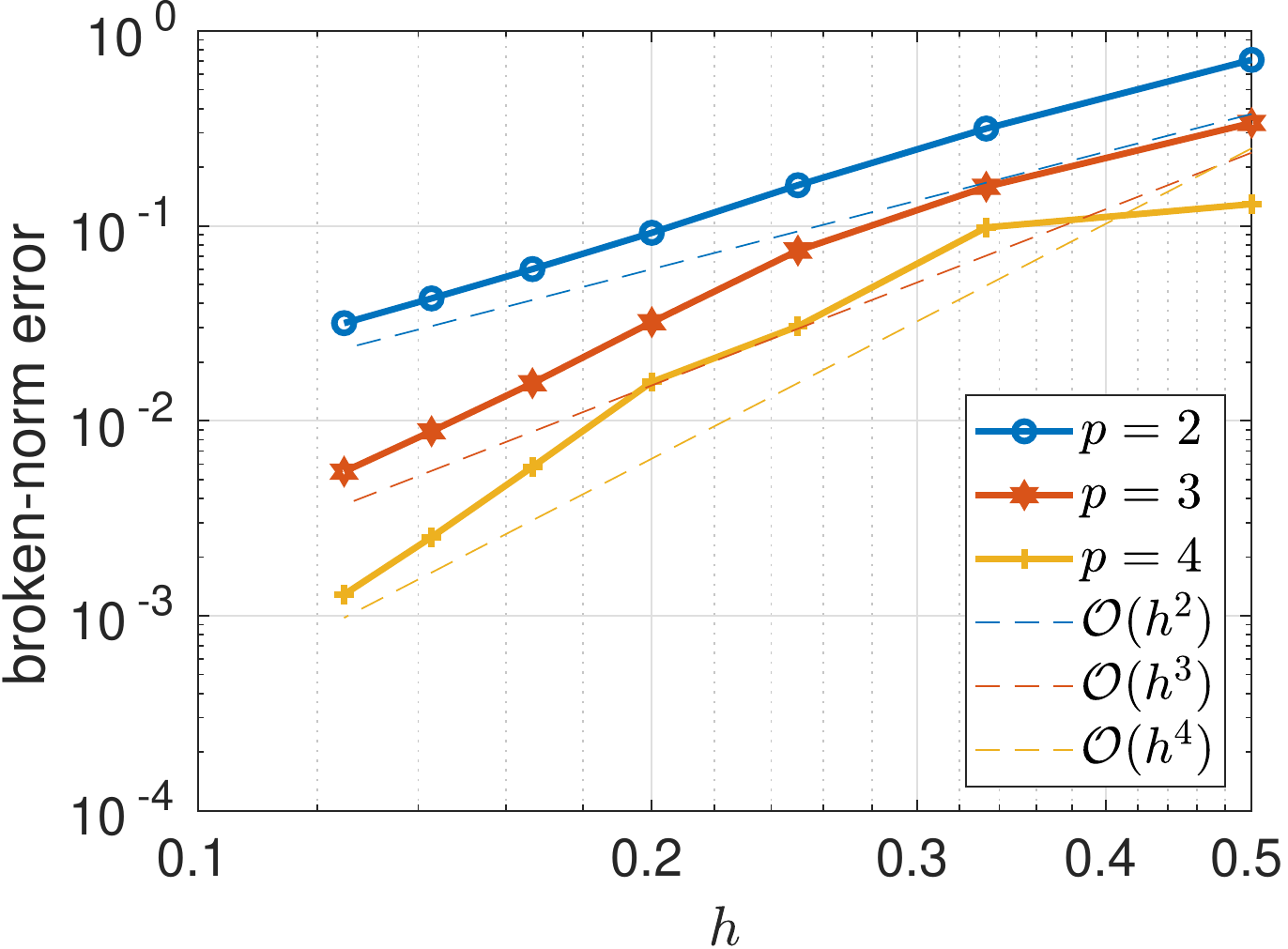}\quad
\includegraphics[trim={0 20 0 0},width=0.44\textwidth]{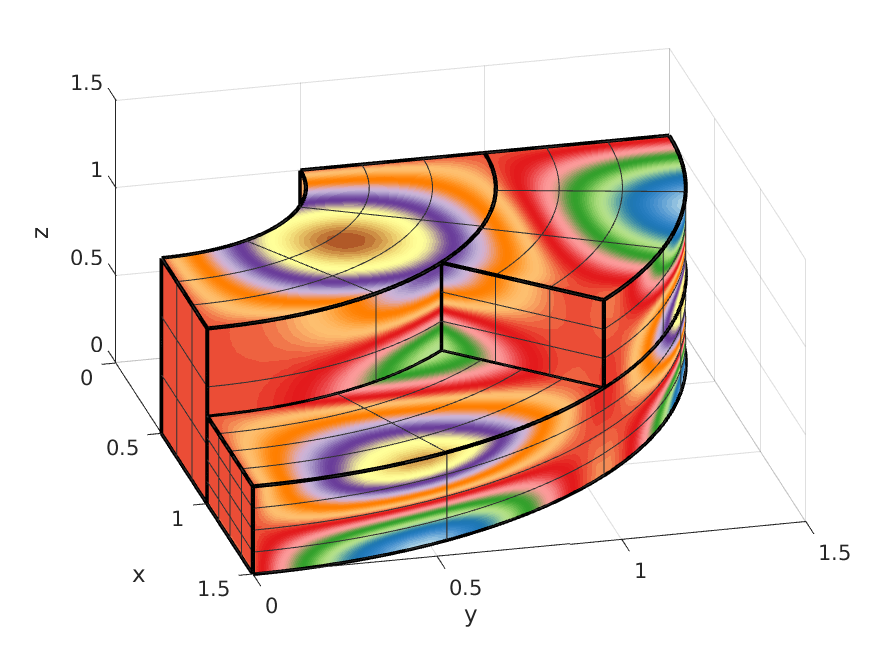}
\end{center}
\caption{\emph{Test case \#6}.
At left, the broken-norm errors versus the
discretization parameter $h$, $p$ is the same in all patches. At right, the
plot of the numerical solution in the patches
 $\Omega^{(1)}\cup \Omega^{(2)}\cup 
\Omega^{(4)}$, computed with: $p^{(1)}=p^{(3)}=4$, 
 $p^{(2)}=p^{(4)}=3$, $3\times 3\times 3$ elements 
in  both $\Omega^{(1)}$ and $\Omega^{(4)}$ and
$4\times 4\times 4$ elements
in  both $\Omega^{(2)}$ and $\Omega^{(3)}$.
 We have removed the patch 
$\Omega^{(3)}$ to have a look on the solution inside $\Omega$}
\label{fig:ring3d}
\end{figure}

The master skeleton is $\Gamma=\gamma_{1,1}\cup 
\gamma_{3,1}\cup\gamma_{3,3}\cup\gamma_{4,3}$ (see the caption of Fig.
\ref{fig:ring3d} for the numbering of the faces).

For $\overline n=2,\ldots,8$ we have considered
 $\overline n\times \overline n\times \overline n$ equal elements 
in both $\Omega^{(1)}$ and $\Omega^{(4)}$, and
$(\overline n+1)\times (\overline n+1)\times (\overline n+1)$ equal elements
in both $\Omega^{(2)}$ and $\Omega^{(3)}$.
The discretizations on the two
sides of any interface are totally unrelated.

To analyze the behaviour of the broken-norm
error with respect to the mesh size, we have considered the same
polynomial degree $p$ along each parametric direction and in all patches.

In the left picture of Fig. \ref{fig:ring3d} we show the broken-norm error 
(\ref{err_broken}) versus the mesh size $h=\max_k h_k$.
We observe that
$\|u-u_h\|_*\simeq {\cal O}(h^p)$ when $h\to 0$ and 
we conclude that
the INTERNODES solution is converging to the exact
one when $h\to 0$ with the best possible convergence rate dictated by the
NURBS-discretization inside the  patches
(see, e.g., \cite[Thm. 3.4 and Cor. 4.16]{bbsv}). 
In the right picture of Fig. \ref{fig:ring3d} the numerical solution
computed with $p^{(1)}=p^{(3)}=4$,
 $p^{(2)}=p^{(4)}=3$, $3\times 3\times 3$ elements
in  both $\Omega^{(1)}$ and $\Omega^{(4)}$ and
$4\times 4\times 4$ elements
in  both $\Omega^{(2)}$ and $\Omega^{(3)}$ is shown.

\subsection{Test case \#7. 3D re-entrant corner}\label{sec:reentrant}

\begin{figure}
\begin{center}
\includegraphics[trim={0 20 0 0},width=0.42\textwidth]{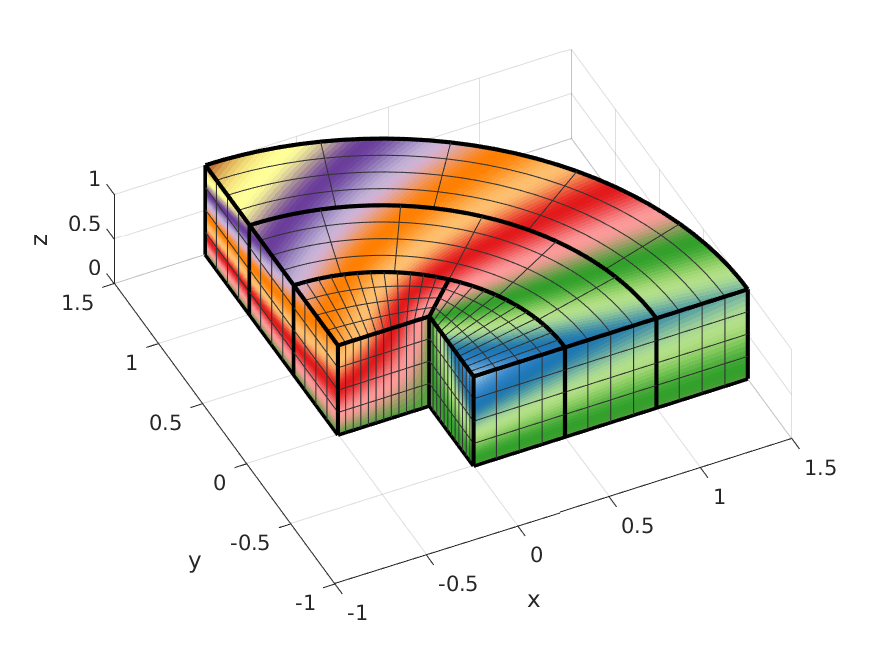}\quad
\includegraphics[trim={0 20 0
0},width=0.42\textwidth]{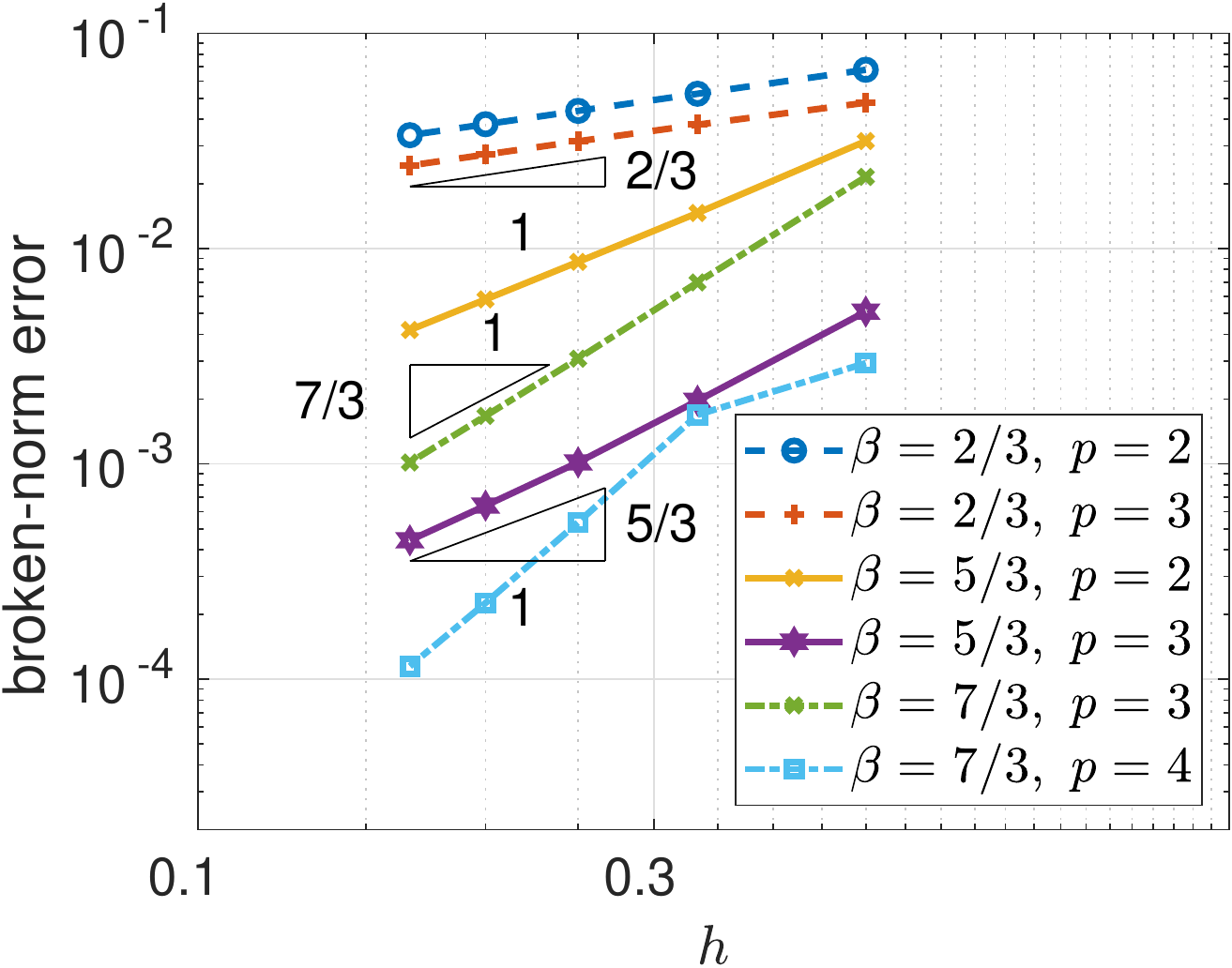}
\end{center}
\caption{\emph{Test case \#7.}
At left, the multipatch configuration. At right, 
the broken-norm errors versus the
discretization parameter $h$, for three values of $\beta$.
The polynomial degree $p$ is the same in all patches }
\label{fig:reentrant}
\end{figure}

Let us consider the differential problem (\ref{problem}) in a non-convex
domain with a re-entrant corner, as shown in Fig.
\ref{fig:reentrant},
more precisely $\Omega=\{(x,y,z)\in{\mathbb R}^3: (x+0.5)^2+(y+0.5)^2\leq 4, \ 
x\geq -0.5,\ y\geq -0.5,\ 0\leq z\leq
1\}\setminus([-0.5,0)\times[-0.5,0)\times[0,1])$.
Then we set $\alpha=0$ , and
$f$ and  $g$ such that the exact solution in cylindrical coordinates reads
$u(r,\theta,z)=r^{\beta}\sin(\beta \theta)\sin(\beta z)$, with $\beta>0$.
The solution features low regularity in a neighborhood of the $z-$axis, in
particular it holds $u\in H^{1+\beta}(\Omega)$.

The computational domain $\Omega$ is split into four patches as shown in Fig.
\ref{fig:reentrant}, all the interfaces are watertight, but with non-matching
parametrizations:
 for $\overline n=2,\ldots,6$ we have discretized the patches as
follows: $\overline n\times 3\overline n\times \overline n$ elements in 
$\Omega^{(1)}$, 
$3\overline n\times \overline n\times \overline n$ elements in 
$\Omega^{(2)}$, $\overline n\times (\overline n+2)\times \overline n$ elements in 
$\Omega^{(3)}$, and $\overline n\times(\overline n+1)\times \overline n$ elements
in $\Omega^{(4)}$.  The first, second, and third  parametric coordinates
 are mapped onto the radial, the angular, and the vertical
physical coordinates, respectively.

In the right picture of Fig. \ref{fig:reentrant} we show the
broken-norm errors for $\beta=2/3$, $\beta=5/3$ and $\beta=7/3$ for the
INTERNODES solution,
that is converging to the exact
one when $h\to 0$ with the best possible convergence rate dictated by the
regularity of the test solution.

\section*{Acknowledgments}
We are very grateful to Luca Ded\`e for his valuable advice.
The research of the first author was partially supported 
by GNCS-INDAM. This support is gratefully acknowledged.

\bibliographystyle{plain}
\bibliography{bibl}

\end{document}